\documentclass[a4paper,12pt]{amsart}

\usepackage{amsmath}
\usepackage{amssymb,amsbsy,amsmath,amsfonts,amssymb,amscd}
\usepackage{latexsym}
\usepackage{graphics}
\usepackage{color}
\input xy
\xyoption{all}

\newcommand\sC{{\mathcal C}}
\newcommand\sT{{\mathcal T}}
\newcommand\sD{{\mathcal D}}

\newcommand\sA{{\mathcal A}}
\newcommand\sF{{\mathcal F}}

\newcommand\sL{{\mathcal L}}

\newcommand\sB{{\mathcal B}}

\newcommand\sH{{\mathcal H}}
\newcommand\LL{{\mathbb L}}
                 \newcommand\sM{{\mathcal M}}

\newcommand\om{\omega}
\newcommand\la{\lambda}

\newcommand\al{\alpha}
\newcommand\be{\beta}
\newcommand\e{\epsilon}
\newcommand\s{\sigma}

\newcommand\Ga{\Gamma}
\newcommand\De{\Delta}
\newcommand\ga{\gamma}
\newcommand\de{\delta}
\newcommand\fie{\varphi}

\DeclareMathOperator{\Mat}{Mat}

\DeclareMathOperator{\Def}{Def}

\def\Bbb{\bf}

\newcommand{\CC}{\ensuremath{\mathbb{C}}}
\newcommand{\RR}{\ensuremath{\mathbb{R}}}
\newcommand{\ZZ}{\ensuremath{\mathbb{Z}}}
\newcommand{\QQ}{\ensuremath{\mathbb{Q}}}

\newcommand{\sS}{\ensuremath{\mathcal{S}}}
\newcommand{\sP}{\ensuremath{\mathcal{P}}}

\newcommand{\NN}{\ensuremath{\mathbb{N}}}
\newcommand{\hol}{\ensuremath{\mathcal{O}}}

\newcommand{\HH}{\ensuremath{\mathbb{H}}}
\newcommand{\BB}{\ensuremath{\mathbb{B}}}
\newcommand{\PP}{\ensuremath{\mathbb{P}}}

\newcommand{\HHH}{\ensuremath{\mathcal{H}}}

\newcommand{\ra}{\ensuremath{\rightarrow}}

\def\eea{\end{eqnarray*}}
\def\bea{\begin{eqnarray*}}

\def\C{{\Bbb C}}

\DeclareMathOperator{\Id}{Id}
\DeclareMathOperator{\Aut}{Aut}
\DeclareMathOperator{\Out}{Out}

\newcommand{\Proof}{{\it Proof. }}

\newcommand\dual{\mathrel{\raise3pt\hbox{$\underline{\mathrm{\thinspace d
\thinspace}}$}}}
\newcommand\qe{\ifhmode\unskip\nobreak\fi\quad $\Box$}       

\def\BOX{\hfill\lower.5\baselineskip\hbox{$\Box$}}

\newtheorem{theorem}[equation]{Theorem}
\newtheorem{theo}[equation]{Theorem}
\newtheorem{remark}[equation]{Remark}
\newenvironment{rem}{\begin{remark}\rm}{\end{remark}}

\newtheorem{defin}[equation]{Definition}
\newenvironment{definition}{\begin{defin}\rm}{\end{defin}}

\newtheorem{question}[equation]{Question}
\newtheorem{prop}[equation]{Proposition}
\newtheorem{cor}[equation]{Corollary}
\newtheorem{lemma}[equation]{Lemma}
\newtheorem{example}[equation]{Example}

\newtheorem{conj}[equation]{Conjecture}

\DeclareMathOperator{\im}{Im}

\def\C{{\Bbb C}}

\newcommand{\restr}[1]{{\raisebox{-0.3\height}{$\mid_{#1}$}}}

\begin{document}

\title[Kodaira and moduli theory]{ Kodaira fibrations and beyond: methods for  moduli theory}
\author{ F. Catanese}
\address {Lehrstuhl Mathematik VIII\\
Mathematisches Institut der Universit\"at Bayreuth\\
NW II,  Universit\"atsstr. 30\\
95447 Bayreuth}
\email{fabrizio.catanese@uni-bayreuth.de}

\thanks{
The present work took place in the realm of the 
 ERC Advanced grant n. 340258, `TADMICAMT' }

\date{\today}

\maketitle
\maketitle
{\em  Contribution to the 16-th Takagi lectures\\   ` In Celebration of the 100th Anniversary of Kodaira's Birth' }

\begin{abstract}
Kodaira fibred surfaces are a remarkable example of  projective classifying spaces,
and there are still many intriguing open questions concerning them, especially the slope question.
The topological characterization of Kodaira fibrations is emblematic of the use of topological methods
in the study of moduli spaces of surfaces and higher dimensional complex algebraic varieties,
and their compactifications. 
Our tour through algebraic surfaces and their moduli (with results valid also for higher dimensional  varieties) 
deals with fibrations, questions on monodromy and factorizations in the mapping class group,
old and new results on Variation of Hodge Structures, especially  a recent answer given (in joint work with Dettweiler)
to a long standing question posed by Fujita.
In the landscape of our tour, Galois coverings, deformations and rigid manifolds
(there are by the way rigid Kodaira fibrations), projective classifying spaces,
the action of the absolute Galois group on moduli spaces, stand also in the forefront. 
These questions  lead to interesting algebraic surfaces, for instance remarkable  surfaces constructed from VHS,
surfaces isogenous to a product
with automorphisms acting trivially on cohomology, hypersurfaces in  Bagnera-de Franchis
varieties, Inoue-type surfaces.

\end{abstract}
\tableofcontents

\section{Kodaira fibrations}

 It is well known that the topological Euler characteristic $e$ is multiplicative
for fibre bundles: this means that, if $ f : X \ra B$ is a fibre bundle with fibre $F$, then
$$ e(X) = e(B) e(F).$$

In 1957 Chern, Hirzebruch and Serre (\cite{chs})
  showed that the same holds true for the signature, also called index  $\sigma = b^+ - b^-$ (it is
  the index of the intersection form on the middle cohomology group) if
the fundamental group of the base $B$ acts trivially on the (rational)  cohomology of 
the fibre $F$.

In 1967 Kodaira \cite {kod67} constructed examples of fibrations
of a complex  algebraic surface over a  curve which are differentiable but not holomorphic  fibre bundles
for which \footnote{indeed, this is true for all such fibrations}the  
multiplicativity of the signature does not hold true. 
  In his honour such fibrations  are nowadays called Kodaira fibrations.
  In fact, for a compact oriented  two dimensional manifold the intersection form
is antisymmetric, hence $\sigma = 0$, whereas Kodaira fibrations have necessarily 
$\sigma > 0$.
  
  As I am now going to explain, there are many interesting properties and open questions 
  concerning Kodaira fibrations.

\subsection{Generalities on algebraic surfaces}

The signature formula of  Hirzebruch, Atiyah and Singer for a compact complex surface $S$ is
$$ \s(S)  = \frac{1}{3} (K^2_S - 2 e(S) ) =  \frac{1}{3} (c_1(S)^2 - 2 c_2(S) ) .$$ 

Here the Euler Number $e(S) $  is the alternating sum of the Betti numbers 
$$ e(S) = 1 - b_1(S) + b_2 (S) - b_3 (S) + 1 =  2 - 2 b_1(S) + b_2 (S),$$
and it equals the second Chern class $c_2(S)$ of the complex tangent bundle of $S$.
Whereas $K_S = - c_1(S)$  is, for an algebraic surface,  the Cartier divisor of a rational
section of the sheaf $\Omega^2_S$ of holomorphic differential 2-forms.

In the K\"ahler  case it was well known that the signature $\sigma (S) = b^+(S) - b^- (S)$ is determined by the Hodge numbers,
indeed $b^+(S) = 2 p_g(S) + 1,  b^+(S) + b^- (S)= b_2(S)$, where 
\begin{itemize}
\item
$p_g(S)  : = h^{2,0} (S) : = h^0(\Omega^2_S)$ is the {\bf geometric genus} of $S$,
\item
$ q(S) : = h^{0,1} (S) : = h^1(\hol_S)$ is the {\bf irregularity} of $S$,
\item
$h^{1,0} (S) : = h^0(\Omega^1_S)$ is the Albanese number (dimension of the Albanese variety).
\end{itemize}

Even if I will concentrate on algebraic or K\"ahler manifolds, I would like to point out how the beautiful
series of papers by Kodaira `On the structure of compact complex analytic surfaces' \cite{kod-scs}
(extending the Enriques classification to non algebraic complex surfaces)
begins with the following miracoulous consequence of the signature formula:
$$ (b^+(S)  - 2 p_g (S) ) + ( 2q(S)  - b_1(S) ) = 1.$$ 
Because both terms are easily shown to be non-negative, by the fact that the intersection
form is positive definite on $H^0(\Omega^2_S) \oplus \overline { H^0(\Omega^2_S) } $,
respectively that  one has  the exact sequence
$$ 0 \ra H^0(d \hol_S)  \ra H^1(S, \CC )   \ra H^1(\hol_S) ,$$
and the Dolbeault  inclusion $ \overline { H^0(d \hol_S)  } \subset H^1(\hol_S)$.

The non K\"ahler case is just the case where the  first Betti number $b_1(S)$ is odd,
$b_1(S) = 2q(S) -1,  b^+(S)  = 2 p_g (S) $; this conjecture of Kodaira was proven through a long series
of papers, culminating in \cite{siuK3}.

The so called `surface geography' problem was raised by van de Ven (\cite{vdv}), and
concerns the points of the  plane with coordinates $(\chi(S), K_S^2)$; here
$\chi(S) : = \chi (\hol_S) = 1 - q(S) + p_g (S)$ is the Euler Poincar\'e characteristic of the sheaf of holomorphic functions.

An important invariant, in the case $\chi(S) \geq 1$ (recall Castelnuovo's theorem: $\chi(S) < 0$ implies that $S$ is ruled) is the so called slope $ \nu: =  K_S^2 / \chi(S)$, whose growth is  equivalent to the one of  the Chern slope,
which is the 
ratio $\nu_C(S) : = c_1(S)^2 / c_2(S)$ of the Chern numbers, because of the Noether formula
$$ 12 \chi(S) =  c_1(S)^2 +  c_2(S) 
\Rightarrow c_1(S)^2 / c_2(S) = K_S^2 / (12 \chi(S) - K_S^2) = \frac{\nu}{12 - \nu}.$$ 

As conjectured by van de Ven, inspired by work of Bogomolov,  and proven by Miyaoka and Yau, there is the slope inequality,
$$ \nu (S) := K_S^2 / \chi(S) \leq 9 \Leftrightarrow   \nu_C (S) : =   c_1(S)^2 / c_2(S) \leq 3$$
for  surfaces with $\chi(S) \geq 1$. 

 It is commonly called the {\bf Bogomolov-Miyaoka-Yau inequality}.

Moreover, by the theorems of Aubin and Yau (\cite{Aubin}, \cite{Yau77}, \cite{Yau78}) and Miyaoka's proof \cite{miya77}, \cite{miyaoka} of ampleness of $K_S$ in the case where equality holds, follows 

\begin{theo}
Is $S$ a non ruled surface and $K_S^2 / \chi(S) =  9$, then the universal cover $\tilde{S}$ of $S$
is biholomorphic to the ball, i.e. the unit disk  $\sD_2  \subset \CC^2$.

\end{theo}

Indeed, the argument of proof about the existence of the K\"ahler Einstein metric seemed to suggest that 
if the ratio  $K_S^2 / \chi(S) $ would be rather close to $9$, then  the universal cover $\tilde{S}$ of $S$
would be diffeomorphic to Euclidean space. More than that, probably the imagination of many was struck by an
impressive result by Mostow and Siu \cite{mostow-siu}

\begin{theo}
There exists an infinite  series of rigid surfaces  $S$  with positive signature whose slope is very close to $9$ \footnote{The maximum achieved by Mostow and Siu is $8,8575$.},
and which admit a metric of nonpositive sectional curvature. In particular, by the theorem of Cartan-Hadamard,
their universal cover $\tilde{S}$ is
 diffeomorphic to Euclidean space $\RR^4$.

\end{theo}

Recently, Roulleau and Urzua \cite{R-U} disposed in the negative of a stronger form of this no name-conjecture 

\begin{theo}
The slopes of simply connected surfaces are dense in the interval $[8,9]$.

 \end{theo}
 
 They used a method introduced by Comessatti and, later, by Hirzebruch (\cite{hirzMA}), to consider abelian coverings of the projective plane $\PP^2$ 
 branched over special configurations of curves, especially configurations of lines; we shall discuss this method at a later moment.
 
 Of  course if the slope $\nu (S)  = 9$ then the universal cover $\tilde{S}$  is the ball, and the fundamental group is countable;
  while the  result of \cite{R-U} does not exclude a priori the existence of a small region around the Miyaoka-Yau
  line $K^2_S = 9 \chi (S)$ where the fundamental group is infinite and the universal cover is homeomorphic to
  a ball, it morally settles the question in the negative. A similar but even bolder conjecture, namely the so-called
  `watershed conjecture', that every surface of positive index should have infinite fundamental group had been disproved thanks
  to Miyaoka's observation that the  Galois closure $S$ of a general projection of a surface  $X \ra \PP^2$ tends to
  have positive index (see \cite{miyaoka}, \cite{m-t},  and \cite{zhije} for other examples).

\subsection{Kodaira fibrations and Kodaira's construction}

We consider in this subsection the following situation: $f : S \ra B$ is a holomorphic map  of a compact complex surface
$S$ onto a curve $B$ of genus $b$, which is a differentiable fibre bundle,  with fibres $F$ of genus $g$.
In this case, the Euler number of $S$ equals
$$  e(S) = 4 (g-1)(b-1).$$
This fibration clearly induces a morphism $\Phi : B \ra \mathfrak M_g$ into the moduli space of curves of genus $g$.

There are two cases: either $\Phi$ is constant, or, as one says, we have a nonconstant moduli fibration. The second case is exactly
the case of Kodaira fibrations. Following a generalization of Kodaira's method developed in our joint paper with Rollenske,
we are going to show how this situation can be effectively constructed, and in such a way that we can calculate the
slope $\nu(S)$ explicitly.

On the other hand, the fact that Kodaira  fibrations do in fact exist, can be shown non effectively in the following way.

\subsubsection{General Kodaira fibrations}

Let $\overline{\mathfrak M_g}^*$ be the Satake compactification of the moduli space of curves of genus $g \geq 3$: this is the closure
of $\mathfrak M_g$, embedded via the Torelli map into $\sA_g$, inside   the Satake compactification $\overline{\sA_g}^*$.

From the moduli point of view, given a moduli-stable curve $C$ of genus $g$, we associate to $C$ the product of the 
Jacobian varieties of the components in the normalization $\tilde{C}$ of $C$.
It follows that the boundary $\partial \overline{\mathfrak M_g}^*$ has dimension $3g -5$, hence codimension $2$
inside $\overline{\mathfrak M_g}^*$ . Similarly the singular set of $\overline{\mathfrak M_g}^*$ corresponds to the locus $\Sigma$  of curves
with automorphisms for $g \geq 4$, and is contained in $\Sigma$ for $g=3$; $\Sigma$  has codimension at least $2$ for $g\geq 4$,
 and codimension equal to $1$ for $g=3$ (its divisorial part is the locus of hyperelliptic curves)

By the projectivity of  $ \overline{\mathfrak M_g}^*$ we can find, for $ g \geq 4$, a smooth  linear section $B$ of $ \overline{\mathfrak M_g}^*$ having
dimension $1$ and  which avoids both 
$\partial \overline{\mathfrak M_g}^*$ and the locus $\Sigma$.

So $B \subset \mathfrak M_g $ and, since $B \cap \Sigma = \emptyset$, we have a  family of curves over $B$ with nonconstant moduli, and the
fibre curves  are all smooth.

We shall see in the sequel that a Kodaira fibration has always positive index, and indeed the Chern slope $\nu_C(S)$  lies in the open interval $(2,3)$.

 One can show that, via  Kodaira fibrations obtained as just described 
   from a general complete intersection curve $B$ 
(under the composition of the Torelli map with the Satake embedding
of the moduli space
$\mathfrak M_g$ of curves of
genus $ g \geq 4$), one obtains a  slope  which is rather small, at most around 2,18 (see \cite{cat-rollenske}).

The above argument shows, via a non explicit construction, the existence of Kodaira fibrations for any fibre genus $ g\geq 4$.
The argument can also be  adapted for the case $g=3$, while 
for $g=2$ there are no Kodaira fibrations; in fact $\mathfrak M_2$ is affine, by a theorem of Igusa \cite{igusaM2}: hence if $B$
is a complete curve, any morphism $ B \ra \mathfrak M_2$ must be constant.

For $g=3$ we take  a smooth curve $D$ which does not intersect $\partial \overline{\mathfrak M_g}^*$ and intersects $\Sigma$ transversely in smooth 
points $p_1, \dots, p_k$ corresponding to general hyperelliptic curves $C_1, \dots, C_k$ (hence such that $\Aut (C_i) = \ZZ/2$). In this case we do not have a family over $D$, 
but only over $ D \setminus \Sigma$.
Take now a double covering $ f : B \ra D$ branched over the points $p_1, \dots, p_k$ and over other points. 
Let $p'_j : = f^{-1} (p_j)$: then we observe that for each $j=1, \dots, k$ the Kuranishi family of $C_j$ has a map $\psi$  to $\mathfrak M_3$
which is  a double covering ramified over the hyperelliptic locus.
In local coordinates we may assume that the map is given by
$$ \psi : (y_1, y_2, \dots, y_5, z) \mapsto  (y_1, y_2, \dots, y_5, z^2) , i.e. , \  w = z^2$$
and the curve $D$ is locally given by  
$$  y_2=  \dots =  y_5 = 0, y_1 = k w ,$$
hence
$$ B = \{  y_2=  \dots =  y_5 = 0, y_1 = k z^2 \}$$
and the family over $B \setminus f^{-1} (\Sigma)$ extends to a family over $B$, with all the fibres smooth genus $3$ curves. 

As discussed with Terasoma after the Takagi lectures, the Kodaira construction does not work for $g=3$. There is however an explicit construction
for $g=3$ due to Zaal \cite{zaal}, which uses Prym varieties, and other explicit constructions of curves 
$B$ as above, by Gonzalez-Diez and Harvey \cite{gh}.

We do not know the answer to the following question, which might turn out to be not too difficult

\begin{question}
Given an integer $ g \geq 3$, which is the number $b(g)$, the minimum value $b \in \NN$ such that there is a Kodaira fibration with fibres of genus $g$
and base curve of genus $b$?
\end{question}

We just observe here that $ b\geq 2$, and indeed $ b \geq 3$ for $g =3,4 $. Because for Kodaira fibrations we have $ e(S) = 4 (g-1)(b-1) $,
 the positive index inequality
  $K_S^2 >   2 e (S) = 8 (g-1)(b-1) $ and the Kefeng Liu  inequality  \cite{Liu96} $K_S^2  <  9  \chi (S)$ hold;
combining with  the Noether formula $K_S^2  = 12 \chi (S) - e(S)$  one gets then
$$   3 (g-1)(b-1)  <  3 \chi (S) < 4 (g-1)(b-1),$$
by which our previous assertion follows immediately. 

Observe moreover that once you have a Kodaira fibration $ f : S \ra B$ with given fibre genus $g$, and base genus $b$, 
for each integer $n \geq 1$ we can take the  Kodaira fibration $ f_n  : S(n)  \ra B$,  where  the fibres of $S(n)$ are the unramified
coverings of the fibres $F$ of $f$ corresponding to the surjection $\pi_1(F) \ra H_1 (F , \ZZ/n)$. Then the fibre genus  of $f_n$
equals $g_n : = 1 + (g-1) n^{2g}$, while the  genus of the base curve  remains $b$. This shows that  $$ {\rm lim inf} _{g \ra \infty} (b_g )= min( b_g).$$

\subsubsection{Kodaira's construction and its generalizations}
Kodaira constructed explicit examples of these Kodaira fibrations, and we are going now  to describe how his method can be generalized.
 
The basic notion for the Kodaira type construction is the following:

\begin{defin}
A  {\bf logarithmic Kodaira fibration} is a quadruple $(X,D, f,B)$ consisting of
\begin{enumerate}
\item a smooth fibration $\psi:X \to B$ of a surface to a curve, 
  with fibres $F_t$ and
\item  a divisor $D \subset X$ such that
\begin{enumerate}
\item the projection $D \to B $ is \'etale and
\item  the fibration of pointed curves  $(F_t, F_t \cap D )$ is not isotrivial, i.e., the fibres are not all isomorphic (as pointed curves).
\end{enumerate}
\end{enumerate}
\end{defin}

Now, even if at first sight  this does not seem to help,  it really does: because now the fibration $\psi$ might have constant moduli, but  the points $F_t \cap D $
may be  moving. 

The easiest case to consider is the case where $X$ will be a product of curves $X:=B_1\times B_2$ and  $D$ shall be a divisor
such that the first projection
$D \to B_1$ is \'etale and the second projection
$D \to B_2$ is finite.

We shall now see that  in order to construct Kodaira fibrations it 
suffices to construct
log-Kodaira fibrations.

\begin{prop}
Let $(X ,D )\to B$ be a log-Kodaira fibration and  let
$f:\tilde F\to F$ be a Galois-covering of a fibre $F$,
with Galois group $G$, and branched over $D\cap F$. Then we can extend $f$
to a ramified covering of surfaces $\bar f: S \to \tilde X$  obtaining a diagram

\[\xymatrix{\tilde F\ar[d]^f\ar@{^(->}[rr] && S\ar[d]^{\bar f}\\
F\ar@{^(->}[rr] & &\tilde X:=g^*X\ar[d]\ar[dl]\\
F\ar@{=}[u]\ar@{^(->}[r]&X\ar[d]&\tilde B\ar[ld]^g\\
&B}\]
where 
\begin{itemize}
 \item $g:\tilde B\to B$ is an \'etale covering,
\item $\tilde X$ is the pullback of $X$ via $g$,
\item $\bar f$ is a ramified covering with Galois group $G$ branched 
over $\tilde D:=g^*D$ and such that $\bar f\restr{\tilde F}=f$.
\end{itemize}

\end{prop}

The idea of the proof (for which we refer to \cite{cat-rollenske}) is that the covering $\bar f\restr{\tilde F}=f$ is determined by a monodromy homomorphism
$$\mu : \pi_1(F \setminus F\cap D) \to G,$$ which extends on neighbouring fibres by the differential
local triviality of the logarithmic fibration: but then there will be an action of $\pi_1(B)$  transforming the monodromy
into another one. Since however  there is only  a finite set of possible such monodromy homomorphisms to the finite group $G$,   we get a tautological finite 
 \'etale  covering $g:\tilde B \to B$ associated to the `monodromy of the monodromy',  and then the ramified covering extends to the pull back of $X$.
 
 The problem is thus reduced to finding disjoint \'etale correspondences between the two curves $B_1, B_2$, which give  the connected components of 
 the curve $D$ we are looking for.
 
 In turn, the easiest case of an \'etale correspondence is given by the following situation: $D'$ is a curve with an automorphism group $H$,
 and there are two subgroups $H_1, H_2 < H$ acting freely on $D'$. We then set $B_j : = D' / H_j$.
 Then $D'$ admits a morphism $\psi : D' \to B_1 \times B_2$, such that the composition of $\psi$ with both projections is \'etale.
 
 If the intersection of the  two subgroups $H_0 : = H_1 \cap H_2$ is non trivial, the map $\psi$ factors through the quotient $D' / H_0$;
in any case, $\psi$ is injective if and only if there are no points $x \neq y $ such that $H_1 x = H_1 y$ and $H_2 x = H_2 y$
(i.e., there do exist $  h_1 \in H_1 \setminus \{1_H\}, h_2 \in H_2$ such that $ y = h_1 x$, $ x = h_2 h_1 x$). Equivalently, 
$$D'   \to B_1 \times B_2 \ {\rm embeds}\  D'   \Leftrightarrow  H_2 (H_1 \setminus \{1_H\})  \cap \sS = \emptyset,$$ 
 where $\sS = \{ h \in H| \exists x \ s.t. h x = x\}$ is the set of stabilizers.

\subsection{Slopes of Kodaira fibrations}

An interesting and open question, raised by Le Brun,  asks for the possible values of the Chern slope of a Kodaira fibred surface
$ f : S \ra B$.

The Chern slope 
$\nu_C (S) : =  c_1^2 (S )/ c_2(S)  = K^2_S / e(S)$ 
of a  Kodaira
fibred surface  lies in the
interval $(2,3)$, in view of the well known Arakelov inequality (that shall be discussed in a later section) and of the
improvement by Kefeng Liu (\cite{Liu96}) of the
Bogomolov-Miyaoka-Yau inequality to $ K^2_S / e(S)  < 3$.

Le Brun 
raised the question whether the slopes of Kodaira fibred surfaces 
can be effectively bounded away from 3:  is it true  that there exists $\e >0$ such that 
    for a Kodaira fibred surface $S$  we have $\nu_C (S) \leq 3 - \e?$ 
      
   The examples by Atiyah, Hirzebruch and Kodaira have  slope
not greater than   $ 2 + 1/3 = 2,33 \dots $ (see \cite{bpv}, page 221)
and, as observed already, if one considers Kodaira fibrations obtained
   from a general complete intersection curve
under the composition of the Torelli map with the Satake embedding
in the moduli space
$\mathfrak M_g$ of curves of
genus $ g \geq 3$, one obtains a smaller slope (around 2,18).

Now, given a Kodaira fibration $f : S \ra B$, and any  holomorphic map of curves $ \varphi : B' \ra B$,
one can take the pull-back of $f$, namely the fibred product $ S' : = B' \times_B S$.
The slope remains the same if $\varphi$ is \'etale, but in case where we have a ramified map,
then the slope decreases.

In fact, if we denote by $d$ the degree of $\varphi$ and by $r$ the degree of the ramification divisor on $B'$, then $b'-1= d (b-1) + r/2$,
hence $e(S') = 4 d(g-1) (b-1 + \frac{r}{2d})$, while $K_{S'}^2 = d K_S^2 + 4 r (g-1)$.

Hence $$ [\nu_C(S) - \nu_C(S')]  [ 4 (g-1)(b-1)  (b-1 + \frac{r}{2d})] = \frac{r}{2d} (K_S^2 -  8 (b-1) (g-1) ),$$
which is strictly positive as soon as $r > 0$.

We observe moreover that the slope
 $$\nu_C(S') = \frac{  K_S^2 + 4 \frac{r}{d}(g-1)} {   4 (g-1) (b-1 + \frac{r}{2d})} $$
 tends to $2$ as soon as $\frac{r}{d}$ tends to infinity.

Therefore, once one has found a given slope, it looks more like a question 
of book-keeping to show that one can realize smaller slopes.
While the hard question seems to be the one of  finding higher slopes:
for this reason we concentrate our attention on the problem of finding Kodaira fibrations with high slope.

The best known result in the direction of high slope   is the following result of \cite{cat-rollenske}:

\begin{theo}[Catanese-Rollenske] There are Kodaira fibrations with slope equal to $
2 + 2/3 = 2,66\dots $.\end{theo}

Our method of construction has been  a variant of the one used by Kodaira, namely to consider {\em double Kodaira fibred surfaces}.

The first main point is that the slope of a Kodaira fibred surface $ f : S \ra B'$ obtained as a Galois branched covering of a logarithmic
Kodaira fibred surface  $\psi : (X, D) \to B$  is determined by the logarithmic structure of $(X,D)$,
namely, given the components $D_1, \dots, D_r$ of $D$, one associates to $D_i$ the branching integer $m_i$ of $f$ along the divisor $D_i$,
and  for instance the canonical divisor $K_S$ is the pull-back of the logarithmic divisor
$$K_X + \sum_i (1 - \frac{1}{m_i} ) D_i.$$

As a consequence, if $d = |G|$, $G$ being the Galois group, then, since the curves $D_i$ are disjoint ($D \ra B$ being \'etale):
$$ K_S^2 = d (K_X + \sum_i (1 - \frac{1}{m_i} ) D_i )^2 = $$$$=d (K_X^2 + 2 \sum_i (1 - \frac{1}{m_i} ) K_X D_i + \sum_i (1 - \frac{1}{m_i} )^2 D_i ^2).$$

Similarly, since in this case $D$ is smooth, 
$$ e (S) =   d [ e(X) - \sum_i (1 - \frac{1}{m_i} ) e (D_i )] =  d [ e(X) +  \sum_i (1 - \frac{1}{m_i} ) ( K_X D_i + D_i ^2)].$$

Since $(1 - \frac{1}{m_i} )^2 < (1 - \frac{1}{m_i} )$, we see that the slope $\nu_C(S)$ has little chance to  become larger unless  $D_i^2 < 0.$ 

This is however often the case: write $K_X = p^* K_B + K_{X|B}$, and observe that $K_{D_i} = p^* K_B $. Then $D_i^2= K_B D_i - K_X D_i = - K_{X|B} D_i$,
which is negative if the relative canonical divisor $K_{X|B}$ is nef (this fact  is, under suitable assumptions, 
 a consequence of  Arakelov's  theorem, and is used for the proof of the Mordell 
conjecture over function fields).

\subsection{Double Kodaira fibrations and group theory}

In \cite{cat-rollenske} a logarithmic Kodaira fibration was defined to be {\bf very simple}  if $X = B \times B$ 
and each $D_i$ is the graph of an automorphism of $B$, where the genus $b$ of $B$ is at least two.

In this case $D_i^2 = - 2 (b-1), K_X D_i = 4 (b-1), K_X^2 = 8 (b-1)^2 , e(X) = 4 (b-1)^2$ and for the Chern slope we have
$$ \nu_C(S) = 2 + \frac{ 2 (b-1) \sum_i   (1 - \frac{1}{m_i^2} ) }{ 4 (b-1)^2 + \sum_i   (1 - \frac{1}{m_i} )2 (b-1)} =   2 + \frac{  \sum_i   (1 - \frac{1}{m_i^2} ) }{ 2 (b-1) + \sum_i   (1 - \frac{1}{m_i} )} .$$ 

Now, as the indices $m_i \to \infty$, the slope tends to 
$$ \nu_C = 2 + \frac{r}{2(b-1) + r}  = 2 + \frac{\al} { 2 + \al}= 3 - \frac{2} { 2 + \al}, \ \ \al: = \frac{r}{(b-1)}.  $$

Hence one can maximize the slope if one can maximize the ratio $ \al: = \frac{r}{(b-1)}$.

One can indeed be more clever (as in loc. cit.) and take $m_i = 3$, so that 
$$  \nu_C(S)  = 2 + \frac{8 r}{18(b-1) + 6 r}  = 2 + \frac{4 \al} { 9 + 3 \al}.   $$

Then , for $\al = 3$, we obtain $$\nu_C(S)  = 2 + \frac{12}{18}= \frac{8}{3}.$$

 However, there are limits to maximizing $\al$:  first of all, by Hurwitz' theorem, $\al < 84$,
since $|Aut(B)| \leq 84 (b-1)$ \footnote{Observe that if equality holds there are several automorphisms with fixed points.}.

This inequality is however much milder than the one given by the BMY-inequality, which implies $\al \leq 9$ (indeed 
we can show that $\al < 8$).

Now, the account for the drop from $84$ to a much lower constant  is due indeed to  a  cogent restriction, namely, 
that we want all the curves $D_i$ to be disjoint!

We already showed that we would like to find an $\al >3$, in particular the action of the group $G : = Aut(B)$ on $B$
cannot be free (otherwise $|G| \leq 2 (b-1)$!).

Therefore, we denote as above by $\sS$ the subset of stabilizers in $G$,
$$ \sS : = \{ g \in G | \exists x \in B, \ s.t. \ \  g x = x\}. $$

Let $D_i = \{ (x, g_i x)\}$. Then $$D_i \cap D_j \neq \emptyset \Leftrightarrow \exists   x \in B, \ s.t. g_i  x = g_j x \Leftrightarrow g_i^{-1} g_j \in \sS.$$

The group theoretical question that we have therefore in mind is a sort of sphere packing problem for groups $G$ acting on a curve.
We define the sphere with center $g$ to be the set $ g \sS $.

Our problem reduces to find a maximal number $r$ such that there are elements $g_1, \dots g_r$ so that $\forall i$ 
the sphere  $ g_i \sS $ contains only the element $g_i$, and no other $g_j \neq g_i$.

Moreover, by the proof of Hurwitz's theorem, and some easy arguments, one sees that the quotient  $ B/G = \PP^1$,
and then that  $G$ is a finite quotient of a  polygonal group
$$  T (n_1, n_2, \dots n_k) = < \ga_1, \dots, \ga_k| \ga_1 \dots \ga_k = 1, \ga_i^{n_i} = 1 \ , \ \forall i=1, \dots k > $$
which is hyperbolic, i.e. $\sum_{i=1}^k (1 - \frac{1}{n_i} ) > 2$.

We wonder whether Lubotzky' s theory of expander graphs might yield  a possible method to approach the question 
of finding a large $\al$, and also of giving a better  upper bound for $\al$ (\cite{lubo}).

\subsection{Moduli of Kodaira fibrations}
The study of moduli of Kodaira fibrations was initiated by Jost and Yau \cite{j-y83} (after that Kas \cite{kas} had  proven that the small deformations of Kodaira fibred surfaces are unobstructed),
showing that all deformations yield again a Kodaira fibred surface. 

Kodaira surfaces are a typical issue of the case where topology determines the moduli space, and the best characterization was obtained
by Kotschick \cite{kot}

\begin{theo}\label{kotschick}
Let $S$ be a complex surface.
A  Kodaira fibration on $S$ with fibres of genus $g$ and base curve of genus $b$ is equivalent to the datum of 

\begin{enumerate}
\item
an exact sequence
$$ 1 \ra \Pi_{g} \ra   \pi_1(S) \ra  
\Pi_{b} \ra  1 $$
 (here $\Pi_g$ denotes the
fundamental group of a compact curve of genus $g$) such that: 
\item

\[e(S)=4(b-1)(g-1),\]
\item
the monodromy homomorphism $ m : \Pi_b \ra Out (\Pi_g)$,
induced by conjugation in the previous exact sequence, has infinite image.
\end{enumerate}

\end{theo}

Kotschick used  the methods that we shall review in section \ref{fibred}, especially  the following  facts:

i) a fibration of a K\"ahler
manifold over a curve of genus $b \geq 2$ without multiple fibres
is determined by a surjection $\pi_1(X)  \ra \Pi_b$ with finitely generated kernel \cite{barlotti}, \cite{cime03}

2) The Zeuthen-Segre formula says that, in the case where $X$ is a surface $S$, 
$ e(S) \geq 4(b-1)(g-1)$, equality holding, when $g \geq 2$, iff  we have a differentiable bundle

3)  the fibration has constant moduli if and only if the image of the monodromy $m$ is finite.

We see then that every surface with the same fundamental group and Euler number as $S$ is again
Kodaira fibred with the same fibre genus $g$ and base genus $b$.

Of course, a major question which remains is:

\begin{question}\label{connected}
Let $S$ be a Kodaira fibred surface: do then the surfaces with the same
fundamental group and Euler number as $S$ form 
a connected component of the moduli space, or the union of a connected component
with  its  complex conjugate component?

\end{question}

One may moreover ask the following question.

\begin{question}\label{3fibred?}
How many Kodaira fibrations can a given algebraic surface possess?

Do there exist surfaces $S$ with three distinct Kodaira fibrations?
\end{question}

Kodaira's original example, whose generalization was explained  in the previous subsection, shows that $S$ can have
two distinct Kodaira fibrations, and one can indeed see that this is again a topological condition.

\begin{prop}\label{char} Let $S$ be a complex surface.
A double Kodaira fibration on $S$ is equivalent to the datum of
 two exact sequences
\[\xymatrix{1\ar[r] &\Pi_{g_i}\ar[r]& \pi_1(S) \ar[r]^{\bar\psi_i}&
\Pi_{b_i}\ar[r]& 1 & i=1,2},\] (here $\Pi_g$ denotes as before  the
fundamental group of a compact curve of genus $g$) such that:

\begin{enumerate}
\item
the monodromy homomorphisms $m_i : \Pi_{b_i} \ra Out (\Pi_{g_i})$
have infinite image;
\item
$b_i\geq2, g_i \geq 3$,
\item the composition homomorphism
\[\xymatrix{\Pi_{g_1}\ar[r]& \pi_1(S) \ar[r]^{\bar\psi_2}&
\Pi_{b_2}}\] is neither zero nor injective, and
\item  the Euler characteristic of $S$ satisfies
\[e(S)=4(b_1-1)(g_1-1)=4(b_2-1)(g_2-1).\]
\end{enumerate}
\end{prop}

The above result shows that surfaces admitting a double
Kodaira fibration form a closed and open subset
in the moduli spaces of surfaces of general type;
since for these one has a realization as a
 branched covering
$ S \ra B_1 \times B_2$, branched over a  divisor $D
\subset  B_1
\times B_2 $, it makes sense to distinguish
the \'etale case where $D$ is smooth and the two
projections $ D \ra B_i$ are \'etale.
It is not clear a priori that this property is also open and closed,
but in \cite{cat-rollenske} we were able to prove it.

\begin{theo}
Double \'etale Kodaira fibrations form a closed and open subset
in the moduli space of surfaces of general type.
\end{theo}

Concerning the previous question \ref{connected} there are only  partial results,  Jost and Yau studied the deformations of
the original example of Kodaira, while \cite{cat-rollenske} considered a more general question 
for an important  special class of double \'etale Kodaira fibrations.  To describe the latter results we recall  a
definition from \cite{cat-rollenske}.

\begin{defin}
A double \'etale Kodaira fibration is said to be {\bf standard} if the logarithmic double Kodaira
fibred surface $(B_1 \times B_2 , D)$ reduces, after \'etale base changes for $B_1$ and $B_2$,
to the very simple case of a logarithmic double Kodaira
fibred surface $(B \times B , D)$ for which $D$ is a union of graphs of automorphisms of $B$.
\end{defin}

\begin{theo}\label{standardmoduli}
The subset of  the moduli space corresponding to
\emph{standard} double \'etale Kodaira fibred
surfaces $S$ with a fixed fundamental group   consists  of at most two  irreducible connected components,
exchanged by complex conjugation, which
are isomorphic to the moduli
space of   pairs $(B, G)$,  where $B$ is a curve of genus $b$ at
least two and
$G$ is a group
   of biholomorphisms of $B$ of a given topological type.
\end{theo}

An interesting by-product of the study of the moduli space of standard double \'etale Kodaira fibrations is the following 
result which contradicted a spread belief (\cite{mostow-siu}), and shows that Kefeng Liu's theorem that a Kodaira fibration cannot exist 
 on a (free) ball quotient cannot be shown invoking non-rigidity.
 
\begin{theo}\label{Kodairarigid}
There are double Kodaira fibred surfaces $S$
which are rigid.
\end{theo}

Using these explicit descriptions 
Rollenske went further in \cite{soenkeCM} and showed that,  in the case where the branched cover has a cyclic Galois group,
then the closure of this irreducible component inside the Koll\'ar-Shepherd Barron-Alexeev  \cite{ksb} compactification
is again a connected component. One should observe that there are extremely few examples where the 
connected components of  the KSBA compact moduli space
 have been investigated, apart from the obvious case of rigid surfaces (see \cite{wenfei} for the case of surfaces isogenous to a product). 
 
 Finally, concerning the existence problem for Kodaira fibrations, we have the following
 
 \begin{question}\label{existence}
Given an exact sequence 
$$ 1 \ra \Pi_{g} \ra   \pi \ra  
\Pi_{b} \ra  1 $$
such that the image of $m : \Pi_{b} \ra \sM ap_g : = Out^+ (  \Pi_{g} )$
is infinite, when does there exist a Kodaira fibred surface $S$ with $\pi_1(S) \cong \pi$?
\end{question}

An obvious necessary condition is that the Abelianization $ \pi^{ab} = \pi / [ \pi, \pi]$
has even rank, by Hodge theory.

To the monodromy $m$ is associated a continuous map $ f : B \ra \sT_g / Im (m)$,
where $\sT_g$ denotes Teichm\"uller space, equivalently, an $m$-equivariant map of the universal cover $\tilde{B}$,
$\tilde{f} : \tilde{B} \ra \sT_g$. It is at present not clear to me if it is proven that $f$ can be deformed to a harmonic map;
the main difficulty seems however to be to show the holomorphicity of such a harmonic map, for which the condition on the Betti number
being even is the first obstruction.

Recently Arapura \footnote{Talk at the Conference for Bob Friedman's 60-th birthday, in May 2015, and \cite{arap}.}
observed that  there are other necessary conditions.   Let  $m_H$ be  the monodromy on $V: = \Pi_{g}^{ab} \otimes \QQ$,
i.e. the monodromy on cohomology, 
and let $G$ be the connected component of the identity in the Zariski closure of $ Im (m_H)$. Then we can first of all 
 replace the above condition on the parity of the first Betti number as the condition that 
 the space $V_G$ of $G$-coinvariants (the largest quotient on which $G$ acts trivially)
 has even dimension.
 
Then there are necessary conditions in special cases: for instance, in  the case where the space $V_G$ of coinvariants is zero,
   $G$ must be semisimple of classical Hermitian type.
In other words, if $V_G= 0$, then $ G(\RR)^0 / K$ must be  a Hermitian Symmetric Domain of classical type
(here $K$ is a maximal compact subgroup).

\section{Projective varieties which are classifying spaces}

\subsection{Generalities on projective varieties which are classifying spaces}

\begin{defin}

Define  $\sP \sC$ as  the class of projective varieties  which are classifying spaces for their fundamental 
 group $\pi_1(Z)$: equivalently, $\sP \sC$ is  the class of projective varieties $Z$ whose universal covering $\tilde{Z}$  is contractible.

\end{defin}

The class $\sP \sC$ is stable for Cartesian products, and for \'etale coverings, hence also for the relation of isogeny.

\begin{defin}
Two varieties $X,Y$ are said to be isogenous if there exist a third variety $Z$, and \'etale finite morphisms $ f_X : Z \ra X$,
$ f_Y : Z \ra Y$.

\end{defin}

The class $\sP \sC$ is however not stable for taking hyperplane sections: because for a compact manifold which is a classifying space
its real dimension is read off by the top nonzero cohomology group $H^m(\pi , \ZZ /2)$, where $\pi = \pi_1(Z)$: since if $Z$ is a classifying space,
or a $K(\pi,1)$ as one  also says in topology, then $H^m(\pi , R) \cong H^m(Z , R)$, for any ring of coefficients $R$.

Projective curves $C$ of genus $g$ are, by virtue of the uniformization theorem, in the class  $\sP \sC$ if and only if $ g \geq 1$.
For $g \geq 2$ they are  the quotients $C = \HHH /\Ga$, where $ \HHH : =\{ z \in \CC  | Im (z) > 0 \}$ and 
$ \Ga \subset \PP SL (2, \RR)$ is a discrete subgroup isomorphic to  $ \Pi_g $ (then necessarily the action is free and cocompact).

Primary examples of  projective varieties which are  $K(\pi,1)$'s are curves and Abelian varieties, and the varieties which are
isogenous to a product of these. Particularly interesting are the varieties isogenous to a product of curves of genera at least $2$.

It is interesting to observe (\cite{bms}, corollary 82) that Abelian varieties are exactly the projective $K(\pi, 1)$ varieties, for which $\pi$ is an abelian group.

\subsubsection{Locally symmetric manifolds of negative type}
A very interesting class of projective varieties which are  $K(\pi,1)$'s are the {\bf locally symmetric manifolds} $Z$  with ample canonical divisor $K_Z$.
These are in some sense a  generalization  in higher  dimension of curves of genus $g \geq2$: because the upper half plane
$ \HHH : =\{ z \in \CC  | Im (z) > 0 \}$   is biholomorphically equivalent to
the unit disk $\{ z \in \CC | |z| < 1\}$.

Locally symmetric manifolds $Z$  with ample canonical divisor $K_Z$ (also called locally Hermitian symmetric manifolds of negative curvature)
are the quotients of a {\bf bounded symmetric domain} $\sD$ by a cocompact discrete subgroup $\Ga \subset Aut (\sD)$
acting freely.  

Recall that a bounded symmetric domain $\sD$ is a bounded domain $\sD \subset \subset \CC^n$ such that its
group $Aut (\sD)$ of biholomorphisms contains, for each point $p \in \sD$, a holomorphic automorphism 
$ \s_p$ such that $ \s_p (p)= p$, and such that the derivative of $ \s_p $ at $ p$ is equal to $ - Id$. This property implies that $\sigma$ is an involution (i.e., it has order 2), and that
$Aut (\sD)^0$ (the connected component of the identity) is transitive on $\sD $ ; therefore one can write $\sD = G/K $, where $G$ is a 
connected Lie group, and $K$ is a maximal compact subgroup.

The classification of these bounded symmetric domains, done by Elie Cartan in
\cite{Cartan}, is based on the fact that such a 
 $\sD$ splits uniquely as the product of irreducible bounded symmetric domains.

 $\sD$ is a complete Riemannian manifold of negative sectional curvature,
hence it is contractible, by the Cartan-Hadamard theorem, and $ Z = \sD / \Ga$ 
 is a classifying space for the group $\Ga \cong \pi_1(X)$.

There are  four series of non sporadic bounded  irreducible domains, in their bounded realization, plus two exceptional types:

\begin{enumerate}
\item[(i)] 
   $I_{n,p} $ is the domain $ \sD = \{ Z \in Mat (n,p, \mathbb{C}) :
\mathrm{I}_p - ^tZ \cdot \overline{Z} > 0 \}$.\\
\item[(ii)] 
   $II_{n} $ is the intersection of the domain  $I_{n,n} $ with the
subspace of skew symmetric matrices.
 \item[(iii)] 
   $III_{n} $ is instead the intersection of the domain  $I_{n,n} $
with the subspace of  symmetric matrices.

\item[(iv)]   The Cartan - Harish Chandra realization of a domain of type $IV_{n}$ in $\mathbb{C}^{n}$ is the subset $\sD$ defined by the inequalities
(compare \cite{Helgason2}, page 527)

\vspace{.3cm}

\begin{tabular}{c}
 $|z_1^2 + z_2^2 + \cdots + z_n^2 | < 1$ , \\
  $1 + | z_1^2 + z_2^2 + \cdots + z_n^2 |^2 - 2\left( |z_1|^2 + |z_2|^2 + \cdots + |z_n|^2 \right) > 0 \, .$ \\
\end{tabular}\\

\item[(v)] $\sD_{16}$ is the exceptional domain of dimension $d=16$.

\item[(vi)] $\sD_{27}$ is the exceptional domain of dimension $d=27$.

\end{enumerate}
Each of these domains 
is contained in the so-called {\bf compact dual}, which is a Hermitian symmetric spaces of  compact type,
the easiest example being, for type I, the Grassmann manifold.

Among the   bounded symmetric domains are the so called {\bf bounded symmetric domains of tube
type}, those which are biholomorphic to a {\bf tube domain}, a generalized
Siegel
upper half-space
$$ T_{\sC} =   \mathbb{V} \oplus  \sqrt -1  \sC$$ where $\mathbb{V}$ is
a real vector
space and $\sC \subset \mathbb{V}$ is a { \em symmetric cone}, i.e., a self dual homogeneous convex cone
containing
no full lines.

In the case of type III domains, the tube domain is Siegel's  upper half space:
$$  \sH_g : = \{ \tau \in \Mat(g,g,\CC)| \tau = \ ^t\tau, \im (\tau) > 0 \},$$ 
a generalisation of the upper half-plane of Poincar\'e.

Borel proved in  \cite{Bo63} that for each bounded symmetric domain
$\sD$ there exists  a compact free
quotient $ X
= \sD / \Ga$,  called a compact Clifford-Klein form  of the
symmetric domain $\sD$.

 A classical result of J. Hano (see \cite{Hano} Theorem IV, page 886, and Lemma 6.2, page 317 of  \cite{milnorcurv}) asserts that a bounded homogeneous domain
that is the universal cover of  a compact complex manifold is symmetric.

\subsubsection{Kodaira fibrations}

The Kodaira fibrations $ f : S \ra B$ are a remarkable example of surfaces in the class $\sP \sC$. 

 Because $S$ is a smooth projective surface and it is known   that all the fibres of $f$ are smooth curves of genus $ g \geq 3$,
 whereas the base curve  has genus $b \geq 2$. 
 
 By the fundamental group exact sequence
 $$ 1 \ra \Pi_{g} \ra   \pi_1(S)  \ra  
\Pi_{b} \ra  1 $$
the universal cover $\tilde{S}$ is a differentiable fibre bundle over $\tilde{B}$ with fibre $\tilde{F}$,
hence it is diffeomorphic to a ball of real dimension $4$.

By simultaneous uniformization (\cite{Bers}) the universal covering $\tilde{S}$ of a Kodaira fibred surface $S$ is biholomorphic
to a bounded domain in $\CC^2$ (fibred over the unit disk $\De : = \{ z \in \CC | |z| < 1\}$ with fibres isomorphic to $\De$), which is not homogeneous.

 $\tilde{S}$ is not homogeneous by the Hirzebruch proportionality principle;
indeed there are only two bounded homogeneous domains in dimension $2$: the bidisk and the $2$-ball.
The bidisk is biholomorphic to $ \HHH \times  \HHH $, and its group of biholomorphisms is a semidirect product of 
$ Aut (\HHH )  \times  Aut (\HHH ) $ by $\ZZ/2$, in particular one has   Chern forms for the tangent bundle,
invariant by automorphisms, which can be written as 
$$c_1: = \phi(z_1) dz_1\wedge  \bar{dz_1} \oplus \phi(z_2)dz_2\wedge  \bar{dz_2},$$
$$c_2 : = \phi(z_1) \phi(z_2) dz_1\wedge  \bar{dz_1} \wedge dz_2\wedge  \bar{dz_2} = \frac{1}{2} c_1^2.$$

Hence the Chern index $\nu_C= 2$ if the universal cover of $S$  is the bidisk, i.e. $c_1(S)^2 = 2 c_2(S)$; whereas, if the universal cover is the ball,
$c_1(S)^2 = 3 c_2(S)$ by a similar argument (\cite{hirz}).

\subsubsection{Known projective classifying spaces in complex dimension two}

In complex dimension 2, we have the following list of projective classifying spaces:

\begin{enumerate}
\item
Abelian surfaces
\item
Hyperelliptic surfaces: these are the quotients of a complex torus of dimension 2  by a finite 
group $G$ acting freely, and in such a way that the quotient is not again a complex  torus. 

They have $p_g=0, q=1$.

These surfaces were classified by Bagnera and de Franchis (\cite{BdF}, see also \cite{Enr-Sev} and \cite{bpv}) and they are obtained as quotients $(E_1 \times E_2)/G$
where  $E_1, E_2$ are  two elliptic curves, and $G$ is an abelian group acting on $E_1$ by translations, and on $E_2$ effectively and in such a way
that $ E_2/G \cong \PP^1$. In other words, these are exactly the surfaces isogenous to a product of curves of genus $1$ which are not Abelian surfaces.

\item
Surfaces  isogenous to a product of curves $C_1 \times C_2$, where $C_1$ has genus $1$, $C_2$ has genus $g_2 \geq 2$.

These are quotients $ (C_1 \times C_2)/G $, where the finite group $G$ acts freely, and where
 we can assume that $G$ acts by a faithful diagonal action
$$ g(x,y) = (g(x), g(y))$$ (we regard $G$ as $G \subset Aut(C_1), G \subset Aut(C_2)$). These surfaces have Kodaira dimension 1.

\item
Ball quotients.

\item
Bidisk quotients, divided into the {\bf reducible case}, the case of surfaces $ (C_1 \times C_2)/G $, isogenous to a product of curves of genera $g_i \geq 2$, and
the {\bf  irreducible case} 
(for these   rigidity holds, as  proven by Jost and Yau \cite{J-Y85}).
\item
Kodaira fibred surfaces.
\item
Mostow-Siu surfaces of negative curvature (\cite{mostow-siu}): these are branched coverings of ball quotients, admitting a metric of negative scalar curvature.
Their Chern slopes $\nu_C(S) $ are very close to $3$, but strictly smaller than $2,96$. These surfaces are rigid, in particular they are not the 
Kodaira examples of Kodaira fibrations.

\item
More examples are gotten from  coverings of the plane branched over configurations of lines, as we shall discuss in a later section
(see  \cite{zheng} , and \cite{panov}, \cite{panov2} ).

\end{enumerate}

\begin{prop}Every surface in $\sP \sC$ is necessarily minimal, and indeed it contains no rational curve.
\end{prop}
\label{norat}

\Proof
In fact, if we have a rational curve $\phi : \PP^1 \ra S$, then the map  $\phi $ lifts to the universal cover, and is then null-homotopic
because $\tilde{S}$ is contractible. Hence also $\phi $ is null-homotopic, and the image rational curve $C : = \phi (\PP^1 )$
is homologous to zero; this is a contradiction, as soon as $S$ is a K\"ahler surface. A fortiori, if $S$ is projective.

\qed

\begin{cor}
Surfaces in $\sP \sC$ which are not of general type are exactly type 1), 2) and 3) above. Any surface in $\sP \sC$
which is of general type has ample canonical divisor.
\end{cor}
\Proof
The result follows by minimality and  by Enriques' classification, if the Kodaira dimension is $<1$.

If $S$ has general type, then it contains no rational curves, in particular no $(-2)$-curves, hence  $K_S$ is ample.

If the Kodaira dimension is one, then $S$ is minimal, propery elliptic, which means that a multiple of the canonical divisor
yields a morphism  $ f : S \ra B$ with general fibre an elliptic curve. By Kodaira's classification of singular fibres of  elliptic fibrations
 (\cite{kod1}) follows that all the fibres
are either smooth elliptic, or multiple of a smooth elliptic curve. We use now the orbifold fundamental group exact sequence
(see \cite{cko}, \cite{barlotti} or \cite{cime03})
$$ \pi_1 (F) \ra \pi_1(S) \ra \pi_1^{orb}(f) \ra 1 .$$
Here the orbifold fundamental group is defined as $\pi_1(B^*) / K$, $B^*$ being the set of non critical values of $f$,
and $K$ is normally generated by $\ga_i^{m_i}$, for each geometric loop going around a point $p_i \in B \setminus B^*$ 
for which the fibre $f^{-1} (p_i)$ is multiple of multiplicity $m_i$.
 Hence there is an intermediate covering $\hat{S} \ra \hat{B}$, possibly of infinite degree, such that all the fibres are smooth elliptic curves, and where $\hat{B}$
 is simply connected.
 
There is also a finite ramified covering $B' \ra B$ such that the pull back $f' : S' \ra B'$ has all the fibres which are smooth.
The $j$-invariant is constant on $B'$ since $B'$ is projective. Hence all the smooth fibres are isomorphic to a fixed elliptic curve $E$
and therefore
we obtain another finite cover $B'' \ra B'$ such that $S'' = E \times B'' $.
Since the Kodaira dimension of $S$ is one, we obtain that $B''$ has genus at least two, and there exists another \'etale covering $C \ra B''$ such that
$ C \ra B$ is Galois hence $ S = (E \times C) / G .$ 

\qed

\subsection{Galois conjugate of projective classifying spaces}

Let $X \subset \PP^n$ be a complex projective variety. Then, for each 
 $\s \in Aut (\CC)$,  the conjugate variety $X^{\s}$ is the set $\s(X)$: it is the projective variety 
 defined by the  ideal $I_X^{\s}$,   obtained from the ideal $I_X$ of $X$ applying the homomorphism $\s$ to the coefficients
 of the polynomials in $I_X$.
 
If $\s$ is complex conjugation, we get  $\overline{X}$, which is diffeomorphic  to $X$.

Observe that, by the theorem of Steiniz,  one has a surjection $ Aut (\CC) \ra Gal(\bar{\QQ} /\QQ)$,
and that  we have an action of the absolute Galois group  $ GAL : = Gal(\bar{\QQ} /\QQ)$
 on the set of   varieties $X$ defined over $\bar{\QQ}$,.
 
 $X$ and the conjugate variety $X^{\s}$  have the same Hodge numbers and Chern numbers. In particular,
for curves, the genus is preserved.

It is immediate also that Galois  conjugation by $\s \in Aut (\CC)$ preserves products and the equivalence relation given by isogeny,
indeed Galois conjugation does not change the algebraic fundamental group, as shown by Grothendieck \cite{sga1}.

\begin{theo}
Conjugate varieties $ X,  X^{\s} $ have isomorphic algebraic fundamental groups 
$$ \pi_1(X)^{alg} \cong  \pi_1(X^{\s} )^{alg} ,$$
($\pi_1(X)^{alg}$ is the profinite completion
of the topological fundamental group $G: = \pi_1(X)$, i.e. $\pi_1(X)^{alg}$ is the inverse limit
of the factor groups $G/K$, $K$ being a normal subgroup of finite index in $G$).

\end{theo}

It is easy to see (\cite{bms}, theorem 223) that 

i) If $X$ is an Abelian variety,  the same holds for any Galois conjugate  $X^{\s} $.

ii) If $S$ is a Kodaira fibred surface, then any  Galois conjugate $S^{\s} $ is also Kodaira fibred.

The following attempt of conjecture is based mainly on the fact that it holds for all known examples.

\begin{conj}
Assume that $X$ is a projective $ K (\pi, 1)$, and assume $\s \in Aut (\CC)$.

Is then the conjugate variety $  X^{\s} $ still a classifying space $ K (\pi', 1)$?

\end{conj}

We know since long, thanks to the result obtained by  J.P. Serre \cite{serre} in the 60's ,  that it is not true in general that $\pi_1 ( X^{\s} ) \cong \pi_1 (X)$.
 Serre showed   that there exists a field automorphism  $\s$  in the absolute Galois group $Gal(\bar{\QQ} / \QQ)$,
and a variety $X$ defined over a number field, such that $X$ and the Galois conjugate variety $X^{\s}$  have non isomorphic fundamental groups, in particular they are not homeomorphic.

This is also false for surfaces in the class $\sP \sC$, for instance in a joint paper with I. Bauer and F. Grunewald \cite{bcg2}
we obtained  the following theorem.

  \begin{theo} \label{fundamentalgroup}  If $\sigma \in Gal(\bar{\QQ} /\QQ)$ is not in the conjugacy class of complex conjugation,
then there exists a surface isogenous to a product $S$ such that $S$ and the Galois conjugate surface $S^{\sigma}$  have non-isomorphic fundamental groups. 
\end{theo}

 I. Bauer and F. Grunewald and the author  (\cite{almeria}, \cite{bcg2}) discovered  also many  explicit  examples of algebraic surfaces
isogenous to a product 
for which the same phenomenon holds (observe that the proof of the general theorem is, as it may be surmised, non constructive).

\begin{rem}
Gonzal\'ez-Diez and Jaikin-Zapirain \cite{gabino} later extended  theorem \ref{fundamentalgroup} to all 
automorphisms $\sigma$ different from
complex conjugation.

\end{rem}

\subsection{Some characterizations of locally symmetric varieties and Kazhdan's theorem in refined form}

Proceeding with other projective $K(\pi, 1)$'s, the question becomes more subtle and we have to appeal to a famous theorem
by Kazhdan on arithmetic varieties (see \cite{Kazh70}, \cite{Kazh83},  \cite{milne}, \cite{CaDS}, \cite{CaDS2}, \cite{V-Z}).

Here the result is in the end  much stronger: not only the conjugate variety is again locally symmetric, 
but the universal cover is indeed the same bounded symmetric domain!

\begin{theo} Assume that $X$ is a projective manifold with
$K_X$  ample, and that
the universal covering $\tilde{X}$ is a bounded symmetric domain.

 \noindent
 Let $\tau \in
\mathrm{Aut}(\mathbb{C})$ be an automorphism of $\mathbb{C}$.

 \noindent Then the conjugate variety  $X^{\tau}$ has universal covering
$\tilde{X^{\tau}} \cong \tilde{X}$.
\end{theo}

The above result rests in an essential way on the Aubin-Yau theorem (\cite{Yau78}, \cite{Aubin}) about the existence
of a K\"ahler Einstein metric for a projective manifold with ample canonical divisor $K_X$, and on  the results of Berger \cite{Berger}.

These results  allow precise algebro-geometric characterizations of such locally symmetric varieties. These results were pioneered
by Yau \cite{Yau88} \cite{Yau93},   his treatment of the non tube case is  however incorrect (he claims that (1) of theorem  \ref{tubes}, which characterizes the tube case,
 holds also in the non tube case).

Simpler proofs follow from recent results  obtained together with Antonio Di Scala. These results yield a simple and precise characterization of varieties possessing a  bounded symmetric
domain as universal cover, without having to resort to the existence of a finite \'etale covering where the holonomy splits.

For the tube case, we have the following theorem (see \cite{CaDS}), whose simple underlying idea is best illustrated in the case where the universal covering is a polydisk.

In this case one observes that, due to the nature of the automorphism group of $\HHH^n$ as semidirect product of $Aut(\HHH)^n$ with the symmetric group
in $n$ letters, the following tensor 
$$\Psi : = \frac{dz_1 \cdot \dots \cdot dz_n }{  dz_1 \wedge \dots \wedge  dz_n}  $$ 
is a semi-invariant for the group of automorphisms, it is multiplied by $\pm 1$ according to the sign of the corresponding
coordinates permutation.

It therefore descends to a section $0 \neq \phi  \in
H^0(S^n(\Omega^1_X)(-K_X) \otimes \eta)$, where $\eta$ is the  2-torsion invertible sheaf
associated to the sign character of the fundamental group $\pi_1(X) \cong \Ga  <  Aut(\HHH^n)$
(observe that, depending on the choice of $\Ga$, $\eta$ may be trivial).

The existence of such a tensor is  unfortunately, in complex dimension $n \geq 4$, not 
a characterizing property of polydisk quotients.

 Indeed, in dimension $4$, we have 
the bounded domain  $ \Omega \subset \CC^4 $,   $$ \Omega =     \{ Z \in  Mat (2,2,\CC) \  :  \  \operatorname{I_2} - ^{t} Z
 \cdot \overline{Z} >0\},$$
  the bounded (Harish-Chandra)  realization of the Hermitian
symmetric space 
  $SU(2,2) / S (U(2) \times U(2))$.
  
  Here the holonomy action of
$  (A, D) \in S (U(2) \times U(2))$ is given
by
$ Z \mapsto A Z D^{-1}$. Hence,  the square of the
determinant of $Z$ yields a section $\Psi$ which descends to a section $0 \neq \phi  \in
H^0(S^n(\Omega^1_X)(-K_X))$.

The main difference with respect to the polydisk quotient case is that the corresponding hypersurface in the projectivized 
tangent bundle of $X$ is nonreduced, 
since we start from 
an invariant hypersurface of degree $4$ which is twice a
smooth quadric.  With these examples in mind it should be easier to get the flavour of the following theorem.

\begin{theo}\label{tubes}

Let $X$ be a compact complex manifold of dimension $n$ with  $K_X$ ample.

Then the
following two conditions (1) and (1'), resp. (2) and (2') are equivalent:

\begin{itemize}

\item[(1)] $X$ admits a slope zero   tensor $0 \neq  \psi   \in
H^0(S^{mn}(\Omega^1_X)(-m K_X) )$, (for some  positive
integer $m$);

\item[(1')]$X \cong \Omega / \Gamma$ , where
$\Omega$ is a bounded symmetric domain of tube type and $\Gamma$ is a
cocompact
discrete subgroup of
$\mathrm{Aut}(\Omega)$ acting freely.
\item[(2)] $X$ admits a semi special tensor $0 \neq \phi  \in
H^0(S^n(\Omega^1_X)(-K_X) \otimes \eta)$, where $\eta$ is a 2-torsion invertible sheaf,
 such that there is a  point
$p\in X$ for which the
corresponding hypersurface $F_p : = 
\{\phi_p = 0 \}\subset \PP (TX_p)$ is  reduced.
\item[(2')]
The universal cover of $X$ is a polydisk.
\end{itemize}

Moreover, in case (1),  the degrees and the multiplicities of the irreducible
factors of the polynomial  $\psi_p$ determine uniquely the universal
covering
$\widetilde{X}=\Omega$.

\end{theo}

The crucial underlying  fact is the following discovery of Kor\'anyi-V\'agi.

Let $D \subset \CC^n$ be a homogeneous bounded symmetric
domain in its circle realization around the origin $0 \in
\CC^n$.

Let $K$ be the isotropy group of $D$ at the origin $0 \in
\CC^n$, so that we have $ D =  G / K$.

   A polynomial $f \in \mathbb{C}[X_1,\dots,X_n]$ is said to be $K$- semi-invariant  if there
    is a character $\chi \colon K \ra   \CC$ such that, for all $g \in K$,
$f(g X) =
\chi(g) f(X)$.

   Observe that, since $K$ is compact, we have:  $|\chi(g)| = 1$.

Let $D = D_1 \times D_2$ be the decomposition of $D$ as a product of
two domains where $D_1$ is of tube type and $D_2$ has no irreducible
factor of
tube type. \\

\begin{theorem}\cite[Kor\'anyi-V\'agi]{KoVa79} \label{KV}

Let $D = D_1 \times D_2$ be the above decomposition and let moreover
   $$D_1 = D_{1,1} \times D_{1,2} \times \cdots \times D_{1,p}$$ be the
decomposition of $D_1$ as a product of irreducible tube type domains
$D_{1,j},
\, \, \, (j=1,\cdots,p)$.

Then there exist, for each $j = 1, \dots p$, a unique
$K_j$-invariant polynomial $N_j (z_{1,j})$, where $K_j$ is the
isotropy subgroup of
$D_{1,j}$, such that:

   for all   $K$-invariant polynomial $f$   there exist a constant $c
\in \CC$ and exponents $k_j$ with
   \begin{enumerate}
\item \[ f =  c \prod_{j=1}^p N_{j}^{k_j}, \, \, \, \] hence in particular
\item \[ f(z_1,z_2) = f(z_1) \, \, ,\] where $z_1$ denotes a vector
in  the domain $D_1$ and $z_2 \in D_2$.
\end{enumerate}
\end{theorem}

The above theorem follows almost directly from \cite{KoVa79} by taking into account
that a $K$-invariant polynomial is, up to a multiple, an \emph{inner
function},
i.e., a function such that $ | f(z) | = 1$ on the Shilov boundary of $D$
(on which  $K$ acts transitively). 

Moreover  the polynomials $N_j$ have 
degree equal to the $rank(D_j)$ of the irreducible domain $D_j$,
($rank(D_j)$ denotes  the dimension $r$ of the maximal totally 
geodesic embedded polydisc $\HHH^r \subset D_j$,
  or, 
equivalently, if $ D = G/K$, with $ G = Aut (D)^0$, $rank (D) = 
rank 
(G^{\CC}) =$ the dimension of a maximal
  algebraic torus contained 
in the complexification $G^{\CC}$).

   The characterization is essentially  a   consequence of the unicity of these inner functions,
   and of the inequality $rank(D_j) \leq 
dim(D_j)$, where  equality holds if and 
only if $D_j = 
\HHH$.

In the case where there are non tube domains in the irreducible decomposition, one has to use
some ideas of Kobayashi and Ochiai \cite{KobOchi}, developed by Mok \cite{MokLibro} who introduced and studied  certain characteristic varieties
which generalize the hypersurfaces defined by the Kor\'anyi-V\'agi - polynomials.

\subsubsection{Algebraic curvature-type tensors and their First Mok characteristic varieties}

Consider the  situation where we are given a direct sum

$$T = T_1 \oplus  ... \oplus T_k$$
of irreducible representations $T_i$ of a group $H_i$
 ($T$  shall be the tangent space to a projective manifold
at one point, and $H = H_1 \times \dots \times H_k$ shall be the restricted holonomy group).

\begin{defin}
1)
An algebraic curvature-type tensor is a nonzero  element
$$ \s   \in End (T \otimes T^{\vee}).$$

2) Its first Mok characteristic cone  $\sC\sS \subset  T$ is defined as  the projection on the first factor
of the intersection of $ker (\s)$ with the set of rank 1 tensors, plus the origin:

$$ \sC\sS : = \{ t \in T | \exists t^{\vee} \in T^{\vee} \setminus \{0\}, (t \otimes  t^{\vee} ) \in ker (\s)  \}.$$

3) Its {\bf first Mok characteristic variety} is the subset $\sS : = \PP ( \sC\sS) \subset \PP (T)$.

4) More generally, for each integer $h$, consider
$$ \{  A \in T \otimes T^{\vee} | A  \in ker (\s), {\rm Rank} (A) \leq h  \},$$
and consider the algebraic cone which is its  projection on the first factor
$$ \sC\sS^h : = \{ t \in T | \exists A  \in ker (\s), {\rm Rank} (A) \leq h  , \exists t'  \in T: t = A   t'  \},$$
and define then $\sS^h : = \PP ( \sC\sS^h) \subset \PP (T)$ to be the {\bf h-th Mok characteristic variety.}

5) We define then the {\bf full characteristic sequence} as the sequence
$$ \sS = \sS^1   \subset \sS^2 \subset \dots \subset \sS^{k-1} \subset \sS^k =  \PP (T).$$
\end{defin}

\begin{rem}
In the case where $\s$ is the curvature tensor of an irreducible symmetric bounded domain $\sD$, Mok (\cite{Mok})
proved that the difference sets  $ \sS^h \setminus    \sS^{h-1}$ are exactly all the  orbits of the  parabolic subgroup $P$
 associated to the compact dual
$\sD^{\vee} = G/P$.
\end{rem}

The
concept of an algebraic curvature type
tensor $\s$ can be then used to prove the following theorem.

\smallskip
\begin{theo}\label{noball}
Let $X$ be a compact complex manifold of dimension $n$ with $K_X$  ample.

Then the universal covering $\tilde{X}$ is a bounded symmetric domain without factors isomorphic to  higher dimensional balls
if and only if there is a holomorphic tensor $\s \in H^0 ( End (T_X \otimes T_X^{\vee}))$
enjoying the following properties:

1)there is a  point  $p \in X$, and  a splitting of the tangent space $T = T_{X,p}$

$$T = T'_1 \oplus  ... \oplus T'_m$$

such that the first Mok characteristic cone  $\sC\sS$ is $\neq T$ and moreover $\sC\sS$ splits into m
irreducible components $\sC\sS' (j )$  with

2) $ \sC\sS' (j )  =  T'_1 \times  ... \times  \sC\sS'_j \times    ...   \times T'_m$

3)  $\sC\sS'_j  \subset T'_j $  is the cone over a smooth projective variety $\sS'_j$
unless  $ \sC\sS'_j = 0$ and dim $(T'_j) = 1.$

Moreover, we can recover the universal covering  of $\tilde{X}$ from the sequence of pairs
$ ( dim (\sC\sS'_j ) , dim (T'_j))$.

\end{theo}

 The case where there are ball factors is the case where the Yau inequality is used (see also \cite{V-Z}), and one can also in this case
 give a characterization which  is given in terms of $X$ and not of some unspecified \'etale cover $K'$ of $X$
 (Master Thesis of Daniel Mckenzie, 2013 \cite{daniel}).

A very interesting program, suggested by Yau in \cite{Yau93} is to extend these characterizations to the case
of quotients $ Y: = \sD / \Ga$, where $\sD$ is a bounded symmetric domain, but the quotient need not be compact,
and the action may  be non-free.

 This should be done using logarithmic sheaves $\Omega^1_X (log D)$,
where $(X,D)$ is a normal crossing compactification and resolution of $ (Y, Sing(Y))$.
These sheaves can 
possibly only be defined in orbifold sense (similarly to what is done  in the work of Campana et al. \cite{camp}),
otherwise it is not clear that the tensors which we  considered above, and which descend to $ Y \setminus Sing(Y)$,
do indeed  extend to $X$ logarithmically.

\subsection{Kodaira fibred surfaces and their conjugates}

The bulk of the previous subsection was to show examples where the universal cover of a projective variety
$X$ and of its Galois conjugates $X^{\s}$ are isomorphic.

Kodaira surfaces miracolously show us that we should not hope (by mere wishful thinking)  that the same result
should hold for all varieties in the class $\sP \sC$.

In fact, it was proven by Shabat (\cite{shabat1}, \cite{shabat2})

\begin{theo}
Let $f : S \ra B$ be a Kodaira fibration , and let $\tilde{S}$ be the universal covering of $S$, a bounded domain in $\CC^2$.
Then the fundamental group $\pi_1(S)$ has finite index inside $Aut (\tilde{S} )$.

\end{theo}

We have then the following consequence
\begin{theo}\label{kodconj}
There exist  families $S_t, t \in T$ of Kodaira fibrations whose  universal covers $\tilde{S_t}$ are not isomorphic.
In particular, there exist Kodaira fibred surfaces such that $S$ and some Galois conjugate $S^{\s}$
have non isomorphic universal covering.

\end{theo}

\Proof
It suffices to take a family $S_t, t \in T, \ dim (T) \geq 1$, of Kodaira fibrations where the fibres of the map of $T$  to the moduli space are
finite. Then, by Shabat's theorem, it follows that, for $t'  \in T$, the number of surfaces $S_t$ whose universal cover is isomorphic
to $\tilde{S_t'}$ is countable (finite?), since these surfaces correspond to conjugacy classes of  finite index subgroups of  
$ G' : = Aut (\tilde{S_t'})$ which are isomorphic to $\pi_1(S)$. And $G'$ is finitely presented, as well as $\pi_1(S)$.

For the second assertion, let $ \mathfrak M (S_t)$ be the irreducible component of the moduli space of surfaces of general type
(which is defined over $\QQ$) containing the image of $T$,
and let $S_t'$ be a surface whose moduli point is not algebraic. Then   the set of isomorphism classes of Kodaira fibred 
surfaces with universal covering isomorphic to $\tilde{S_t'}$ is countable, while 
the set of isomorphism classes of its conjugate surfaces 
${S_t'}^{\s}$ is uncountable (note that the corresponding moduli points  do not need a priori to belong to $ \mathfrak M (S_t)$, but this is irrelevant and can be
indeed arranged taking some explicit family of double Kodaira fibrations).

\qed

The result of Shabat was brought to attention by Gonz\'alez-Diez and Reyes-Carocca \cite{gabino-carocca}, who also  discuss the `arithmeticity' condition that a
Kodaira fibration is defined over a number field. They show, as a consequence of Arakelov's finiteness theorem,
that $S$ is arithmetic if and only if the base curve is so. And from this they deduce that   if two such have the same universal cover, and one is arithmetic,
then also the other is arithmetic.

Their work suggests the following possible extension of theorem \ref{kodconj}:

\begin{question}
Does there exist an arithmetic  Kodaira surface $S$ (i.e., $S$ is defined over $\overline{\QQ}$) and an automorphism $\s \in Aut (\CC)$
such that the universal coverings of $S$ and $S^{\s}$ are not isomorphic?
\end{question}

\subsection{Surfaces fibred onto curves of positive genus which are classifying spaces}

We consider now the following situation:  $f : S \ra B$ is a fibration of an algebraic surface $S$ onto a curve $B$
and we want to find a criterion for the universal covering $\tilde{S}$ to be  contractible.

We already remarked in proposition \ref{norat} that if $S \in \sP \sC$  then $S$ must be minimal, and indeed $S$ contains no rational curves. More
generally, the same argument shows that, for each curve $C \subset S$,  the fundamental group of the normalization
$C'$ of $C$ has infinite image in $\pi_1(S)$.

Use now the orbifold fundamental group exact sequence
( see \cite{cko}, \cite{barlotti} or \cite{cime03})
$$ \pi_1 (F) \ra \pi_1(S) \ra \pi_1^{orb}(f) \ra 1 ,$$
hence we obtain a fibration $\tilde{f} : \tilde{S} \ra \hat{B}$, where $\hat{B} \ra B$ is the ramified covering of $B$ corresponding
to the surjection $$\varphi : \pi_1 (B \setminus \{p_1, \dots, p_r\} ) \ra  \pi_1^{orb}(f)$$
(here $p_1, \dots, p_r$ are the points whose fibre is multiple, of multiplicity $m_i \geq 2$, and the kernel of $\varphi$ is normally generated by $\{ \ga_i^{m_i}\}$,
for $\ga_i$ a loop around $p_i$).

Notice that $ \hat{B}$ is compact if and only if $\pi_1^{orb}(f)$ is finite, and if and only if $ \hat{B} \cong \PP^1$,
since $ \hat{B}$ is simply connected. But even when $\pi_1^{orb}(f)$ is infinite
there is a finite ramified covering $B'$ of $B$ such that the pull back of the fibration $f$ has no multiple fibres
(it suffices that the ramification index at each $p_i$ is equal to the multiplicity $m_i$ of the multiple fibre).

By passing therefore to a finite \'etale covering (of the surface $S$), we may assume that $f : S \ra B$ has no multiple fibres,
and we shall assume in the sequel that the genus $b$ of $B$ is $\geq1$.

\begin{question} Are there  examples of  a fibration  $f : S \ra B$  of an algebraic surface $S$ onto a curve $B$, where $S$ is a projective classifying space,
 and $\pi_1^{orb}(f)$  is finite? Equivalently, with $ B \cong \PP^1$ and where $f$ has no multiple fibres?
\end{question}

Assume now that $f$ has no multiple fibres, so that $\hat{B}$ is the universal cover $\tilde{B}$ of $B$, and that the genus $b$ of $B$ is $\geq1$.

Next, by the necessary condition that for each curve $C$ contained in a fibre $\pi_1(C' ) \ra  \pi_1(S)$
has infinite image, it follows that all the fibres of $\tilde{f} : \tilde{S} \ra \tilde{B}$ are homotopy equivalent to CW complexes
of real dimension $1$.

Let $\sC \subset B$ be the set of critical values of $f$, $\sC = \{ p_1, \dots, p_h\}$. Choose a set $\Sigma$ of non intersecting paths joining
a fixed base point $p_0$ with the points of $\sC$, and similarly a set  $\Sigma'$  inside $\tilde{B}$ of non intersecting paths joining
a fixed base point $y_0 \in \tilde{S}$ mapping to $p_0$, with the points of the inverse image of $\sC$, which is a countable set $\{ y_n | n \in \NN , n \geq 1\}$.

\begin{lemma}
Let $f : S \ra B$ be a fibration without multiple fibres
over a curve $B$ of genus $b\geq1$, and assume that, for each irreducible curve  $C$ contained in a fibre of $f$, $C'$ denoting the normalization of $C$,
the homomorphism  $\pi_1(C' ) \ra  \pi_1(S)$
has infinite image.

Since $\Sigma'$ is  deformation retract of $\tilde{B}$, $\tilde{S}$ is homotopically equivalent to a CW complex of dimension $\leq 2$,
in particular  $\tilde{S}$ is contractible if and only if  $H_2(\tilde{S}, \QQ) = 0$. 

\end{lemma} 

\Proof
$\tilde{S}$ retracts onto $K' : = \tilde{f}^{-1} (\Sigma')$, and this shows the first assertion. 

The second follows since,  
$\tilde{S}$ being simply connected and homotopically equivalent to a CW complex of dimension $\leq 2$, it is contractible if and only if
$H_2(\tilde{S}, \ZZ) = 0$. In fact, 
by Hurewicz' theorem, the first non-zero homotopy group $\pi_m (\tilde{S}) $ is isomorphic to $ H_m(\tilde{S}, \ZZ) $,
which is obviously zero for $m \geq 3$. Hence if
$H_2(\tilde{S}, \ZZ) = 0$ all homotopy groups $\pi_i (\tilde{S}) = 0$ and $\tilde{S}$ is contractible (the converse is obvious).

Finally, by the universal coefficient formula, $ H_2(\tilde{S}, \ZZ) $ is torsion free, hence $ H_2(\tilde{S}, \ZZ) = 0$ 
if and only if $H_2(\tilde{S}, \QQ) = 0$. 

\qed

Let $\tilde{F_0}$ be the fibre over $y_0$, and let $\tilde{F_n}$ be the fibre over $y_n$.

We can write $K' = \cup_{n \in \NN} K_n$, where $K_{n}$ is the inverse image of the union $\Sigma'_n$ of the segments joining $y_0$ with $y_i$, for $ i \leq n$.

Clearly $$H_2(\tilde{S}, \QQ) = 0 \Leftrightarrow H_2(K_n, \QQ) = 0 \ \forall n \in \NN.$$
We can now use the theorem of Mayer-Vietoris, using the following notation: $K_n^*$ shall be the inverse image of $\Sigma^*_n$, the union of the open segments, together
with the point $y_0$, whereas $\partial K_n$ shall be the inverse image of $\{ y_1, \dots, y_n\}$.

 Now, $K_n$ is the union of two open subsets, the first is $K_n^*$ which is homotopically equivalent to $\tilde{F_0}$, and the second which is homotopically equivalent
 to $\partial K_n$, the disjoint union of the fibres $\tilde{F_i}, \ 1 \leq i \leq n$.
 Moreover, the intersection of the two open sets  is homotopically equivalent to $n$ disjoint copies of $\tilde{F_0}$.
 
 Define $H_i : = H_1 (\tilde{F_i}, \QQ)$. Then we have the Mayer-Vietoris exact sequence
 
 $$ 0 \ra  H_2(K_n, \QQ)  \ra H_0^n \ra H_0 \oplus H_1 \oplus \dots \oplus H_n \ra H_1(K_n, \QQ) \ra 0.$$ 

Here, we have a  homomorphism $ r_i : H_0 \ra H_i$ which is obtained by the fact that a neighbourhood of the fibre $\tilde{F_i}$ retracts onto $\tilde{F_i}$:
$r_i$ is  a surjection whose kernel is  the group $V_i$ of vanishing cycles.
The homomorphism of the i-th summand $ H_0$ inside $ H_0 \oplus H_1 \oplus \dots \oplus H_n$ has first component which is the identity,
and all the other components equal to zero with the exception of  one which is indeed $ r_i : H_0 \ra H_i$.

We obtain the following theorem.

\begin{theo}
Let $f : S \ra B$ be a fibration without multiple fibres
over a curve $B$ of genus $b\geq1$, and assume that, for each irreducible curve  $C$ contained in a fibre of $f$, $C'$ denoting the normalization of $C$,
the homomorphism  $\pi_1(C' ) \ra  \pi_1(S)$
has infinite image.

Then  $\tilde{S}$ is contractible if and only if, for each $n \in \NN$, the subgroups of vanishing cycles form a direct sum
$V_1  \oplus \dots \oplus V_n$ inside $H_0$. 
\end{theo}

\Proof
The kernel of $H_0^n \ra H_0 \oplus H_1 \oplus \dots \oplus H_n $ are the elements $(x_1, \dots , x_n)$ such that
$(x_1+  \dots + x_n, r_1(x_1), \dots, r_n(x_n)) = (0, 0. \dots, 0)$. Hence these are the elements of 
$V_1  \oplus \dots \oplus V_n$ which map to $0$ inside $H_0$. 

\qed

We observe that, even if the criterion is essentially a characterization of fibred algebraic surfaces which are
projective classifying spaces, the condition on the vanishing cycles is not so easy to verify, and, up to now, good examples are 
still missing (cf. however later sections for candidates).

Observe also the following relation with the well known Shafarevich conjecture. Keep the assumption that
$f : S \ra B$ is a fibration without multiple fibres
over a curve $B$ of genus $b\geq1$, and assume that, for a general fibre $F$, the homomorphism  $\pi_1(F ) \ra  \pi_1(S)$
has infinite image.
 If there is an  irreducible curve  $C$ contained in a fibre of $f$ such that,  $C'$ denoting the normalization of $C$,
the homomorphism  $\pi_1(C' ) \ra  \pi_1(S)$
has finite image, while the  homomorphism  $\pi_1(C ) \ra  \pi_1(S)$ has infinite image, then $S$ would be a counterexample to
the Shafarevich conjecture that $\tilde{S}$ is holomorphically convex (cf. \cite{bog-katz}). Since then $\tilde{S}$ would contain a 
connected infinite chain 
of compact curves.

\section{Surfaces fibred onto curves}\label{fibred}

\subsection{The Zeuthen-Segre formula}

 The Zeuthen Segre formula is a  beautiful formula, valid for any smooth algebraic surface $S$:
  
 \begin{theo}{\bf (Zeuthen-Segre, classical)}
Let $S$  be a smooth projective surface, and let $ C_{\la} , \  \la  \in \PP^1$, be a linear pencil of curves of genus $g$ which meet transversally 
in $\de $ distinct points. If $\mu$ is the number of singular curves in the pencil (counted with multiplicity), then
$$ \mu -  \de - 2 (2g-2) = I + 4, $$ where the integer $I$ is an invariant of the algebraic surface, called Zeuthen-Segre invariant.
 
 \end{theo}
 
 Here, the integer $\de$ equals the self-intersection number $C^2$ of the curve $C$, while in modern terms  the integer $ I + 4$ is not only an algebraic invariant, but is indeed a topological invariant: the topological Euler-Poincar\'e characteristic $e(S)$. 
 
 Observe that $ I + 4 + C^2$ is then the  topological Euler-Poincar\'e characteristic
 of the surface blown up in the $\de = C^2$ base points of the pencil. In this way the Zeuthen-Segre formula generalizes 
 for every fibration $f : S \ra B$ of an algebraic surface onto a curve $B$, and the formula measures the deviation
 from the case of a topological bundle, for which one would have $e(S) = 4 (b-1) (g-1)$, where $g$ is 
 the genus of a fibre.
 
 The importance of the formula is that, as suggested by the interpretation of $\mu$ as a number of points,
 the difference $\mu : = e(S) - 4 (b-1) (g-1)$ is always non negative.
 
 The formula is well known using topology (see \cite{bpv}), but it is very convenient to have an algebraic formula, which is proven using 
 sheaves and exact sequences (see the  lecture notes \cite{cb}).

\begin{defin}
Let $ f : S \ra B$ be a fibration of a smooth algebraic surface $S$ onto a curve of genus $b$, and consider a singular  fibre 
$F_t  = \sum n_i C_i$, where the $C_i$ are distinct  irreducible curves.

Then the divisorial singular locus of the fibre is defined as the divisorial part of the critical scheme,
$ D_t : =  \sum (n_i -1) C_i$, and the Segre number of the fibre is
defined as 
$$\mu_t : = deg \sF + D_t K_S  - D_t^2,$$
where the sheaf  $\sF$  is concentrated in the singular points of the reduction of the
fibre, and is the quotient of $\hol_S$ by  the ideal sheaf  generated by the components of the vector 
$d \tau / s$, where $ s = 0 $ is the equation of $D_t$, and where $\tau$ is the pull-back of a 
local parameter
at the point $t \in B$.

More concretely, $$ \tau  = \Pi_j   f_j^{n_j} , s = \tau / ( \Pi_j   f_j),$$ and the logarithmic derivative yields

$$ d \tau = s [  \sum_j n_j (df_j \Pi_{h \neq j} f_h)].$$

\end{defin} 

The following is the refined Zeuthen-Segre formula
\begin{theo}{\bf (Zeuthen-Segre, modern)}
Let $ f : S \ra B$ be a fibration of a smooth algebraic surface $S$ onto a curve of genus $b$, and with fibres of genus $g$.

Then $$c_2(S) = 4 (g-1)(b-1) + \mu,$$ where $\mu = \sum_{t\in B} \mu_t$,
and $\mu_t \geq 0$ is defined as above. Moreover, $\mu_t$ is strictly positive, except if the fibre is smooth or a multiple of a smooth
curve of genus $g=1$. 

\end{theo}

Most of the times, the formula is used in its non refined form: if $g >1$, then either $\mu>0$, or $\mu=0$ and we have
a differentiable fibre bundle. In this case there are two alternatives: either we have a Kodaira fibration, or all the smooth
fibres are isomorphic, and we have a holomorphic fibre bundle (\cite{fg65}).

\subsection{The positivity results of Arakelov,  Fujita and Kawamata}

\subsubsection{Arakelov's theorem}

\begin{defin}$
A \  fibration \ f : S \ra B$  of a smooth algebraic surface $S$ onto a curve of genus $b$ is said to be (relatively) minimal, if there is no $(-1)$-curve contained in a fibre.
Moreover, it is said to be {\bf isotrivial}, or with constant moduli,  if all the smooth fibres are isomorphic.
\end{defin}

Isotriviality is equivalent to the condition that the moduli morphism $ \psi : B \ra \mathfrak M_g$ is constant, and it implies that
there exists a finite Galois base change $B' \ra B$ such that the pull back $S' \ra B'$ is birational to a product.

The theorem of Arakelov (\cite{ara}) gives a numerical criterion for isotriviality.

\begin{theo}\label{ara}
Let  $f : S \ra B$  be a minimal fibration  of a  smooth algebraic surface $S$ onto a curve of genus $b$, where the genus $g$ of the fibres is strictly positive.
 Define $K_{S|B}$, the relative 
canonical divisor, as $ K_{S|B} : = K_S - f^* (K_B)$. Then $K_{S|B}$ is nef, and big unless the fibration is isotrivial.
In particular $K_{S|B}^2 \geq  0 $ and $K_{S|B}^2 > 0 $ if the fibration is not isotrivial.
\end{theo}

\Proof(Idea)
(I) The minimality of the fibration and $g \geq 1$ ensure that $K_{S|B} \cdot C = K_S \cdot C \geq 0$ for each curve $C$
contained in a fibre.

(II) For the case where $C$ is not contained on a fibre, we use that the line bundle $\hol_C(K_{S|B})$ is generically a quotient of the pull-back
of $V: = f_*\hol_S(K_{S|B})$.  Indeed, if  $ p : C \ra B$ is induced by $f$, there is a non zero morphism
$p^* V \ra  \hol_C(K_{S|B})$ and one applies the theorem of Fujita \ref{Fuji1} (that we shall soon describe) stating that $p^*V$ is nef, hence
$\hol_C(K_{S|B})$ has  non negative degree.

(III) Since the divisor $K_{S|B}$ is nef, either $K_{S|B}$ is big or the graded ring associated to it has Iitaka dimension $1$,
and yields a map to a curve $C$, $\varphi : S \ra C$. We consider $F : =  f \times \varphi : S \ra B \times C$,
and we consider the Hurwitz formula: $K_S = F^* (K_B + K_C) + R$, which proves that
$K_{S|B} = F^* ( K_C) + R$. Since all the sections of multiples of $K_{S|B} $ pull back from $C$, it follows that $R$ is horizontal, hence
all the fibres are ramified covers of the same curve $C$, and branched in the same set of points.
From this it follows that all the smooth fibres are isomorphic.

\qed

The following corollary contains an observation by Beauville \ref{beachi}and the fact, already mentioned several times,
that a Kodaira fibred surface has positive index.

\begin{cor}
Under the same assumptions as in Arakelov's theorem \ref{ara}, but assuming $g \geq 2$: then $\chi(S) \geq (g-1) (b-1)$, equality if and only if 
we have a holomorphic bundle. 

In particular, a Kodaira fibred surface $S$  has a strictly positive index, i.e.  $S$ has $c_1^2 (S) = K^2_S > 2 e (S) = 2 c_2(S)$.
\end{cor}
\Proof
The Arakelov inequality $K_{S|B}^2 = K_S^2 - 8 (b-1)(g-1) \geq 0 $ and the Zeuthen-Segre inequality $e(S) \geq  4 (b-1)(g-1)$ add up,
in view of Noether's theorem (the first equality in the next formula)
to yield the new inequality $$12\chi(S) = K_S^2  + e(S) \geq 12 (b-1)(g-1).$$ Moreover, equality holds if both Arakelov's and Zeuthen-Segre's inequality 
are equalities, implying that we have an isotrivial fibration and a differentiable bundle: hence, as we already observed, $f$ is a holomorphic bundle.

If we have a Kodaira fibration, then $e(S) =  4 (b-1)(g-1)$ and the fibration is not isotrivial,
hence $K_S^2  > 8 (b-1)(g-1) = 2 e(S)$.

\qed 

\subsubsection{Fujita's direct image theorems}

An important progress in classification theory was stimulated by  a theorem of Fujita, who showed (\cite{fuj1}) that the direct image of
the relative dualizing sheaf $$ \om_{X|B}  =  \hol_X (K_{X|B}) : =  \hol_X (K_X - f^* K_B)$$ is numerically positive, i.e. every quotient bundle
has non negative degree
(a fact that was used in the previous subsection to give a different proof of Arakelov's theorem).

\begin{theo}{\bf Fujita's first theorem}\label{Fuji1}
If $X$ is a compact K\"ahler
manifold and $ f : X \ra B$ is a fibration onto a projective curve $B$ (i.e., $f$ has connected fibres), then the direct image sheaf
$$  V : = f_* \om_{X|B} $$ is a  nef   vector bundle on $B$, equivalently $V$ is `numerically semipositive',
meaning that each quotient bundle $Q$
of $V$ has degree $\deg (Q) \geq 0$.
\end{theo}

In particular , if $X$ is an algebraic surface $S$, then $deg (V) \geq 0$, which is the inequality 
observed previously 
$$ deg (V) = \chi(S) - (g-1) (b-1) \geq 0$$
(except that the  characterization of the case of equality does not follow right away from theorem \ref{Fuji1}, one needs theorem  \ref{fuj2})

In the note \cite{fuj2} Fujita announced the following quite  stronger result:

\begin{theo}{\bf (Fujita's second theorem, \cite{fuj2})}\label{fuj2}

Let $f : X \ra B $ be  a fibration of a compact K\"ahler manifold $X$ over a projective curve $B$, and consider 
the direct image  sheaf $$ V : = f_* \om_{X|B} = f_* ( \hol_X (K_X - f^* K_B)).$$
Then $V$ splits as a direct sum  $ V  = A \oplus Q$, where $A$ is an ample vector bundle and $Q$ is a unitary  flat bundle.
\end{theo}

Fujita sketched the proof, but referred to a forthcoming article concerning the positivity  of the so-called local exponents; the result 
was used in the meantime, for instance I used it  in my joint work with Pignatelli, but since  Fujita's detailed  article was never written, 
there was some objection
(see \cite{Barja}) to the use of this beautiful result:
 this was a motivation for Dettweiler and myself  to write down a complete proof, which is going to appear 
  in the article contributed to Kawamata's 60-th birthday volume
 (\cite{cd}).
 
 \subsubsection{Kawamata's positivity theorems}
 
 In the meantime the idea of the proof had become more transparent, through Kawamata's use of 
 Griffihts' results on Variation of Hodge Structures ( the relation being that
 the fibre of $V: = f_* \om_{X|B}$ over a point $b \in B$, such that $X_b : = f^{-1} (b)$ is smooth,  is  the vector space 
$  V_b   = H^0( X_b, \Omega^{n-1}_{X_b}).$ 

   Kawamata (\cite{kaw0} \cite{kaw1}) improved on Fujita's result, solving a long standing problem
   and proving  the subadditivity of Kodaira dimension for
 such fibrations, $$ Kod (X) \geq Kod(B) + Kod (F),$$ (here $F$ is a general fibre) showing the semipositivity also for 
the direct image of  higher powers of the relative dualizing sheaf 
$$W_m : =  f_* ( \om_{X|B}^{\otimes m}) = f_* ( \hol_X (m(K_X - f^* K_B))).$$

  Kawamata also extended his result to the case where the dimension of the base variety $B$  is $> 1$  in \cite{kaw0},
giving later a simpler proof of semipositivity  in \cite{kaw2}. There has been a lot of literature on the subject ever since,
(see the references we cited in \cite{cd},  see \cite{e-v} for the ampleness of $W_m$ when the fibration is not birationally  isotrivial,
and see \cite{ff14} and \cite{ffs14}). 

Kawamata also  introduced  a simple lemma, concerning
the degree of line bundles on  a curve whose metric grows at most  logarithmically around a finite number of singular points, which played a crucial role for the proof
of Fujita's second theorem.

\medskip

Indeed the  missing details concerning the proof of the second theorem of Fujita, using Kawamata's lemma and some crucial estimates
given by Zucker (\cite{zucker}) for the growth of the norm of sections of the $L^2$-extension of Hodge bundles, were provided in \cite{cd}, 
whose main contribution was   a negative answer  to a question posed by Fujita in 1982 (Problem 5, page 600 of \cite{katata},
Proceedings of the 1982 Taniguchi  Conference).

\subsection{Fujita' s semiampleness question}

To understand the  question posed by Fujita  it  is not only  important to have in mind Fujita's second theorem, but it is also very convenient to recall the following 
 classical definition used by Fujita in \cite{fuj1}, \cite{fuj2}.

 Let $V$ be a holomorphic vector bundle  over a projective curve $B$.

\begin{defin} 
 Let $p :  \PP : = Proj (V) = \PP (V^{\vee}) \ra B$ be the associated projective bundle, and  let $H$ be a  hyperplane divisor (s.t.
 $p_* (\hol_{\PP} (H)) = V$).

Then $V$ is said to be:

(NP)  numerically semi-positive if and only if every quotient bundle $Q$ of $V$ has  degree $deg(Q) \geq 0$,

(NEF) nef if and only if $H$ is nef   on $\PP$,

(A) ample  if and only if $H$ is ample  on $\PP$ 

(SA) semi-ample  if and only $H$ is semi-ample  on $\PP$ (there is a positive multiple  $mH$ such that
the linear system $|mH|$ is base-point free).
\end{defin}
\begin{rem}

Recall that  (A) $\Rightarrow$ (SA) $\Rightarrow$ (NEF) $\Leftrightarrow$(NP), 
the last  follows from the following result  due to Hartshorne.

\end{rem}
\begin{prop}\label{ample}
A vector bundle $V$ on a curve   is  nef if and only it is numerically semi-positive, i.e.,  if and only if every quotient bundle $Q$ of $V$ has  degree $deg(Q) \geq 0$,
and $V$ is ample  if and only if every quotient bundle $Q$ of $V$ has  degree $deg(Q) > 0$.

\end{prop}
Recall the following standard definition, used in the statement of Fujita's second theorem. 

\begin{defin}
A flat holomorphic vector bundle  on a complex manifold $M$ is  a holomorphic vector bundle $\sH : = \hol_M \otimes_{\CC} \HH$,
where $\HH$ is a local system of complex vector spaces associated to a representation $\rho : \pi_1(M) \ra GL  (r, \CC)$,
$$  \HH : =( \tilde{M} \times \CC^r )/  \pi_1(M) ,$$
$ \tilde{M} $ being the universal cover of $M$ (so that $ M =  \tilde{M} /  \pi_1(M)$).

 We say that $\sH$ is  unitary flat if it is associated to a  representation $\rho : \pi_1(M) \ra U (r, \CC)$.

\end{defin}
\bigskip

We come now to Fujita's question:
 
 \begin{question}\label{fujitaquestion} {\bf (Fujita)}
 Is the direct image  $V : = f_* \om_{X|B}$ semi-ample ?
 
 \end{question}

In \cite{cd} we established    a  technical result  which clarifies how Fujita's question is very closely related to Fujita's  II theorem
\begin{theo}\label{semiample}
Let $\sH$ be a unitary flat vector bundle on a projective manifold $M$, associated to a representation $\rho : \pi_1(M) \ra U (r, \CC)$.
Then $\sH$ is nef;  moreover $\sH$ is semi-ample if and only if $ Im (\rho )$ is finite.
\end{theo}

The idea of the proof is to reduce to the case where we have such a bundle over  a curve (obtained by taking successive hyperplane sections of $M$), 
decompose the representation as a direct sum of irreducible unitary representations,
and then  use the theorem of Narasimhan and Seshadri \cite{NS} implying that a unitary flat holomorphic bundle
over a curve is holomorphically trivial if and only if the representation is trivial.

Hence in  our particular case, where $ V  = A \oplus Q$ with $A$ ample and $Q$ unitary  flat, the semiampleness of $V$  simply means that the flat bundle has finite monodromy
(this is another way to wording the fact that  the representation of the fundamental group 
 $ \rho : \pi_1(B) \ra  U (r, \CC)$ associated to the flat unitary
 rank-r bundle $Q$  has finite image).

The main new result in our joint work \cite{cd} was, as already said, to provide  a negative answer to Fujita's question in general:
\begin{theo}\label{maintheorem}
There exist  surfaces $X$ of general type endowed with  a fibration $ f : X \ra B$ onto a curve $B$ of genus $\geq 3$, and with fibres of genus $6$,  such that $V : = f_* \om_{X|B}$ 
splits as a direct sum $ V = A  \oplus Q_1 \oplus Q_2$, where $A$ is an ample   rank-2  vector bundle, and the flat unitary rank-2 summands $Q_1, Q_2$
have infinite monodromy group (i.e., the image of $\rho_j$ is infinite). In particular, $V$ is not semi-ample.

 \end{theo}
 
 Recently (\cite{cd3}) we have found an infinite series of counterexamples, which are based on the same ideas,
 but are quite simpler; we shall report on them in a later subsection (they shall  be called here BCDH-surfaces).
 
 Notice that Fujita's second theorem follows right away from the first  in the case where the base curve is $\PP^1$,
 since then every vector bundle splits as a direct sum of line bundles, and we can separate the summands with strictly positive degree from  the trivial summands.

 Also, one  can say something more precise in the case where the base curve $B$ has genus $1$, or under other assumptions,
 which imply that indeed $V$ is semiample.

 \begin{cor}

Let $f : X \ra B $ be  a fibration of a compact K\"ahler manifold $X$ over a projective curve $B$.
Then  $ V : = f_* \om_{X|B} $
is  a direct sum  $ V  = A \bigoplus (\oplus_{i=1}^h Q_i)$, 
with $A$   ample  and each $Q_i$  unitary  flat 
without any nontrivial degree zero quotient. 
Moreover,

(I) if $Q_i$ has rank equal to 1, then it is a torsion bundle ($\exists \  m$ such that $Q_i^{\otimes m}$
is trivial) (Deligne)

(II) if the  curve $B$ has genus 1, then rank  $(Q_i) = 1 , \ \forall i.$ 

(III)  In particular,  if  $B$ has genus  at most  1, then $V$ is semi-ample.

\end{cor}

{\em Idea of proof:}

(I) was proven by Deligne \cite{TH2}   (and by Simpson \cite{simpsonTors} using the theorem of Gelfond-Schneider), while

(II) Follows since  $\pi_1(B)$ is abelian, if $B$ has genus $1$: 
hence every representation splits as a direct sum of 1-dimensional ones.

\qed

  \subsubsection{How do the flat bundles appear}
  In order to get a fibration $f: S \ra B$ where we have a splitting $V = f_* \om_{S|B}  \bigoplus Q$
  our idea is to use symmetry, for instance the fibres are cyclic  coverings $C$ of $\PP^1$ 
  with group $\ZZ/n$ and branched in $4$ points.
  
  Then we get  a curve $C=C_x$ birational to the curve
described by an equation of the form:
\begin{equation}\label{eq51} z_1^n = y_0^{m_0} y_1^{m_1} (y_1 - y_0 )^{m_2}  (y_1 - x y_0 )^{m_3} , \ \  x \in \CC \setminus \{0,1  \}.
\end{equation}

 Here, we shall make the restrictive assumption that $$ 0 < m_j \leq n-3, \ {\rm  and} \ \  m_0 + m_1 + m_2 + m_3= n.$$ 
  
  Then $C$   admits a Galois cover  $ \phi : C \ra \PP^1$ with cyclic 
Galois group equal to the group of n-th roots of unity in $\CC$, $$ G = \{ \zeta \in \CC^* | \zeta^n = 1\},$$
acting by scalar multiplication on $z_1$.
The choice of a generator in $G$ yields an isomorphism $ G \cong \ZZ/n$.

Now, the vector space $V_x : = H^0 (\Omega^1_{C_x})$ splits according to the characters $\chi: G \ra \CC^*$,
$$  V_x = \bigoplus_{\chi} V^{\chi}$$
and each  character can be written as $\zeta \mapsto \zeta^i$.
  
There is an easy formula (see \cite{DeligneMostow}) for $ dim (  V^{\chi} )$:
 $$ dim (  V^{\chi} ) = \frac{1}{n} ( [im_0] +   [im_1] + [im_2] + [im_3] - n ),  \ i \neq 0 $$
 where $0 \leq [m] \leq n-1$ denotes the remainder of division by $n$
 (a standard representative of the residue class).
 
For the dual character then 
$ dim (  V^{\bar{\chi}} ) = \frac{1}{n} (n - [im_0] +  n-  [im_1] +n- [im_2] + n- [im_3] - n )$
hence $$ dim (  V^{\chi} ) + dim (  V^{\bar{\chi}} ) = 2.$$
Hence, by Hodge Theory, if  $H^{\chi} = V^{\chi}  \oplus V^{\bar{\chi}}$ is the corresponding eigenspace for the action of $G$ on $H^1(C_x, \CC)$,
its dimension equals $2$.

By our choice for the integers $m_i$, for $i=1$ we get a character $\chi_1$ such  that $V^{\chi_1} = 0$,
hence 
$$ H^{\chi_1} = V^{\chi_{n-1}},$$
hence we get a flat summand for $V$.

The next question is: how can we assert that the monodromy group for the flat summand is infinite?

The idea can be explained in simple terms like this: if it were finite, then there would be
a monodromy invariant positive definite scalar product. This would also hold when we conjugate the characters via the action of the absolute Galois group
(for $n$ prime, they are all conjugate). If the monodromy is irreducibile, there is only one monodromy invariant  scalar product, hence
for all conjugate characters we would have a positive definite scalar product, which amounts to the condition that one never has
$ dim (  V^{\chi} ) =  dim (  V^{\bar{\chi}} ) = 1$.

On the other hand, if we assume that  $m_0+m_3$ is    invertible in $\ZZ/n\ZZ,$  there is a $j$ such that $j (m_0+m_3) \equiv -1 \ (mod \ n)$.
For instance $m_0+m_3$ is    invertible in $\ZZ/n\ZZ$ if we take $m_0 =1, m_1= 1, m_2 = 1, m_3 = n-3$
and we assume that $n$ is an odd number.

Define now $m'_i : = [ m_i j]$, where $[a]$ denotes the remainder of division by $n$.  Hence $m'_0+m'_3= n-1$.

We have 
the obvious inequalities $  2\leq m'_1+m'_2\leq 2n-2.$ Hence 
$$ n+1\leq m'_0+\cdots+m'_3\leq 3n-3$$ and therefore 
$$  m'_0+\cdots+m'_3= 2n.$$ 
Hence the underlying unitary form is indefinite for the character $\chi_j$, whereas it is definite for the conjugate character $\chi_{n-j}$.

Moreover, (see \cite{cd3}) once the rank one summand $V^{\chi_{j}} $ is shown to be ample, the irreducibility of the flat bundle associated to $H^{\chi_j} $
(hence of the monodromy) follows for the following chain of  arguments:

\begin{itemize}
\item
If the  flat vector bundle were reducible there would be an exact sequence of flat vector bundles
$$0 \ra \sH' \ra H^{\chi_j}  \ra \sH'' \ra 0 $$
where both $\sH' , \  \sH''  $ have rank 1.
 \item
$V^{\chi_{j}} $ is a holomorphic subbundle of $H^{\chi_j} $, hence it has a nontrivial homomorphism to it
\item
 every homomorphism of $V^{\chi_{j}} $ to a flat rank one bundle $\sH' , \  \sH''$ must be zero
 \item
 contradiction!  
 \end{itemize}

The final conclusion is  that the  monodromy of the flat factor $V^{\chi_{n-1}}$  is infinite, and that $V$ is not semiample; the concrete
 examples of this situation shall be described in section 5.

\subsection{Castelnuovo-de Franchis and morphisms onto curves}

In order to address  the study of moduli spaces of Kodaira fibred surfaces,
or of surfaces which are isotrivially fibred, it is convenient to review some simple results, pioneered by Siu,
and which concern topological conditions implying the existence of holomorphic maps of a compact K\"ahler manifold $X$
onto an algebraic curve.

Siu  used harmonic theory in order to construct holomorphic maps from K\"ahler manifolds to projective curves. The first result in this direction
was the  theorem of \cite{siucurves}, also obtained by Jost and Yau
(see \cite{JostYau}  and also \cite{j-y83} for other results).

\begin{theo}\label{siucurves}{\bf (Siu)}
Assume that a compact K\"ahler manifold $X$ is such that there is a surjection  $\phi : \pi_1(X) \ra \pi_g$, where $ g\geq 2$ and, as usual, $\pi_g$ is the fundamental group of a projective curve of genus $g$.
Then there is a projective curve $C$ of genus $g' \geq g$ and a fibration $f : X \ra C$ (i.e., the fibres of $f$ are connected)  such that $\phi$ factors through $\pi_1 (f)$.
\end{theo}

In this case the homomorphism leads to a harmonic map to a curve, and one has to show that the Stein factorization yields a map to some Riemann surface
which is holomorphic for some complex structure on the target.

It can be seen more directly  how the K\"ahler assumption  is used, because this assumption guarantees that
holomorphic forms are closed, i.e., $ \eta \in H^0 (X, \Omega^p_X) \Rightarrow d \eta = 0$.

At the turn of last century this fact was  used by Castelnuovo and de Franchis (\cite{CdF}, \cite{deFranchis}):

\begin{theo} {\bf (Castelnuovo-de Franchis)}
Assume that $X$ is a compact K\"ahler manifold, $ \eta_1 , \eta_2  \in H^0 (X, \Omega^1_X)$ are $\CC$-linearly independent, and  the wedge product $ \eta_1 \wedge \eta_2$ is d-exact.
Then $ \eta_1 \wedge \eta_2 \equiv 0 $ and there exists a fibration $f : X \ra C$ such that $\eta_1, \eta_2 \in f^* H^0 ( C, \Omega^1_C)$. In particular,
$C$ has genus $g \geq 2$.
\end{theo}

Even if the proof is well known, let us point out that the first assertion follows from the Hodge-K\"ahler decomposition, while $ \eta_1 \wedge \eta_2 \equiv 0 $ implies the existence
of a non constant rational function $\fie$ such that $\eta_2 = \fie \eta_1$. This shows that the foliation defined by the two holomorphic forms has Zariski closed leaves,
and the rest follows then rather directly taking the Stein factorization of the rational map $\fie : X \ra \PP^1$.

Now, the above result, which is holomorphic in nature, combined with the Hodge decomposition, produces results which are topological in nature
(they actually only depend on the cohomology algebra structure of $H^*(X, \CC)$). 

To explain this in the most elementary case, we start from the  following simple observation.  If two linear independent vectors in the first cohomology group $H^1 (X, \CC)$ of a K\"ahler manifold have wedge product which is trivial
in cohomology, and  we represent them as $ \eta_1 + \overline{\om_1},  \eta_2 + \overline{\om_2},$ for $ \eta_1 , \eta_2 , \om_1, \om_2  \in H^0 (X, \Omega^1_X)$,
then by the Hodge decomposition and the first assertion of the theorem of  Castelnuovo-de Franchis 
$$ (\eta_1 + \overline{\om_1}) \wedge (  \eta_2 + \overline{\om_2}) = 0 \in H^2 (X, \CC)$$
implies $$ \eta_1  \wedge   \eta_2 \equiv 0 ,  \  \om_1  \wedge   \om_2 \equiv 0.$$

We can apply Castelnuovo-de Franchis unless $\eta_1  ,   \eta_2$   are $\CC$-linearly dependent, and similarly $\om_1  ,   \om_2$.
W.l.o.g. we may assume then $\eta_2  \equiv 0$ and $\om_1 \equiv 0$. But then $\eta_1 \wedge  \overline{\om_2} = 0$ implies that the Hodge norm
 $$ \int_X (\eta_1 \wedge  \overline{\om_2}) \wedge \overline{(\eta_1 \wedge  \overline{\om_2}) } \wedge \xi^{n-2}= 0 ,$$
 where $\xi$ is here the K\"ahler form. A simple trick is to observe that
  $$ 0 = \int_X (\eta_1 \wedge  \overline{\om_2} ) \wedge \overline{(\eta_1 \wedge  \overline{\om_2}) } \wedge \xi^{n-2}=   - \int_X (\eta_1 \wedge  \om_2) \wedge \overline{(\eta_1 \wedge  \om_2}) \wedge \xi^{n-2} ,$$ therefore   the same integral yields that the Hodge norm of $\eta_1 \wedge \om_2$ is zero,
 hence $\eta_1 \wedge \om_2 \equiv 0;$ the final conclusion is that  we can in any case apply Castelnuovo-de Franchis and find a map to a projective curve $C$
 of genus $ g \geq 2$.
 
 More precisely, one gets  the following theorem (\cite{albanese}):
 
 \begin{theo}{\bf (Isotropic subspace theorem)}
 On a compact K\"ahler manifold $X$ there is a bijection between isomorphism classes of fibrations $ f : X \ra C$ to a projective curve of genus $g\geq 2$,
 and real  subspaces $V \subset H^1 (X, \CC)$ (`real' means that $V$ is self conjugate, $\overline{V} = V$) which have dimension $2g$ and are of the form $ V = U \oplus \bar{U}$, where $U$ is a maximal isotropic subspace for 
 the wedge product $$H^1 (X, \CC) \times H^1 (X, \CC) \ra H^2 (X, \CC).$$
 \end{theo}
 
 The above result, as simple as it may be,  implies the  few relations  theorem of Gromov (\cite{Gromov}),
 which  in turn  implies theorem \ref{siucurves} of Siu (see e.g. \cite{bms} for an ampler discussion).
 
There is another result (\cite{cime03})  which again, like the isotropic subspace theorem, determines explicitly the genus of the target curve
 (a result which is clearly useful for classification and moduli problems).
 
 \begin{theo}\label{orb-fibr}
 Let $X$ be a compact K\"ahler manifold, and let $f : X \ra C$ be a fibration   onto a projective  curve $C$, of genus $g$,
 and assume that there are exactly $r$ fibres which are multiple with multiplicities $ m_1, \dots m_r \geq 2$. Then $f$ induces
an orbifold
fundamental group exact sequence
$$ \pi_1 (F) \rightarrow \pi_1 (X) \rightarrow  \pi_1 (g; m_1, \dots m_r) 
\rightarrow 0,$$  where $F$ is a smooth fibre of $f$, and where the orbifold fundamental group
$$  \pi_1 (g; m_1, \dots m_r) $$
is defined as 
$$   \ \langle  \al_1, \be_1, \dots, \al_g, \be_g, \ga_1, \dots \ga_r | \ \  \Pi_1^g [\al_j, \be_j] \Pi_1^r \ga_i  = \ga_1^{m_1}= \dots = \ga_r^{m_r} = 1\rangle.
$$ 
 Conversely, let $X$ be a compact K\"ahler manifold and let $(g, m_1, \dots m_r)$
be a hyperbolic type,  i.e., assume that $  2g-2 + \Sigma_i (1 - \frac{1}{m_i} ) > 0.$ 

Then
 each  epimorphism
$\phi : \pi_1 (X) \rightarrow  \pi_1 (g;
m_1, \dots m_r)$ with finitely generated kernel is obtained from a fibration $f : X \rightarrow C$ of
type $(g; m_1, \dots m_r)$.
 \end{theo}
 
 With these results and the Zeuthen-Segre formula, it is easy to explain on the one hand the characterization of surfaces isogenous to a product,
 and on the other hand Kotschick's result on moduli of Kodaira surfaces.

\begin{theo}\label{isogenous}

a) A  projective smooth surface $S$ is isogenous   to a product of two curves of respective genera $g_1, g_2 \geq 2$ ,  if and only if
the following two conditions are satisfied:

1) there is an exact sequence
$$
1 \rightarrow \pi_{g_1} \times \pi_{g_2} \rightarrow \pi = \pi_1(S)
\rightarrow G \rightarrow 1,
$$
where $G$ is a finite group and where $\pi_{g_i}$ denotes the fundamental
group of a projective curve of genus $g_i \geq 2$;

2) $e(S) (= c_2(S)) = \frac{4}{|G|} (g_1-1)(g_2-1)$.

\noindent
b) Write $ S = (C_1 \times C_2) / G$. Any surface $X$ with the
same topological Euler number and the same fundamental group as $S$
is diffeomorphic to $S$ and is
also isogenous to a product. 

c) The corresponding subset of the moduli space of surfaces of general type
$\mathfrak{M}^{top}_S = \mathfrak{M}^{diff}_S$, corresponding to
surfaces orientedly homeomorphic,
resp. orientedly diffeomorphic to $S$, is either
irreducible and connected or it contains
two connected components which are exchanged by complex
conjugation.

\end{theo}

{\em Idea of the proof of b)}
 $\Ga : = \pi_1 (S) $  admits a subgroup $\Ga'$ of index $d$ such that 
$ \Ga' \cong  (\pi_{g_1} \times \pi_{g_2})$.  Let $S'$ be the associated unramified covering of $S$. Then application of the isotropic subspace theorem or of theorem \ref{orb-fibr} yields 
a pair of holomorphic maps $f_j : S' \ra C_j$, hence a holomorphic map $$ F : =  f_1 \times f_2 : S' \ra C'_1 \times C'_2.$$

 Then the fibres of $f_1$ have genus $h_2 \geq g_2$,
hence by the Zeuthen Segre formula  $ e(S') \geq  4 (h_2 -1 ) (g_1 - 1)$, equality holding if and only if all the fibres are smooth.

But   $ e(S') = 4 (g_1-1)(g_2-1) \leq 4 (h_2 -1 ) (g_1 - 1)$, so $h_2= g_2$, all the fibres are smooth hence isomorphic to $C'_2$; therefore $F$ is an isomorphism.

\qed

The previous theorem generalizes also to varieties isogenous to a product of curves of dimension $n\geq 3$(\cite{isogenous}, \cite{bms}).

We try now to compare the cited theorem of Kotschick (theorem \ref{kotschick}) with the previous one.

\begin{theo}
Assume that $S$ is a compact K\"ahler surface, and that

(i)  its fundamental group sits into an exact  sequence, where $ g,b \geq 2$:

$$ 1 \ra \pi_g \ra \pi_1(S) \ra \pi_b \ra 1   $$

(ii) $ e(S) = 4 (b-1) ( g-1)$.

Then $S$  has a smooth holomorphic fibration $ f : S \ra B$, where $B$ is a projective curve of genus $b$, and where all the fibres are smooth projective curves
of genus $g$.  $f$ is a Kodaira fibration if and only if the associated homomorphism $ \rho : \pi_b \ra \sM ap_g$ has image of infinite order, else it
is a surface isogenous to a product of unmixed type and where the action on the first curve is free.
\end{theo}

\Proof
 By theorem \ref{orb-fibr} the above exact sequence yields a fibration $ f : S \ra B$ such that there is a surjection $\pi_1(F) \ra \pi_g$, where $F$ is a smooth fibre.
 Hence, denoting by $h$ the genus of $F$, we conclude that $ h \geq g$, and again we can use the Zeuthen-Segre formula to conclude that $h = g$ and that all fibres
 are smooth. So $F$ is a smooth fibration. Let $C' \ra C$ be the unramified covering associated to $ \ker (\rho)$: then the pull back family $ S' \ra C'$ has
 a topological trivialization, hence is a pull back of the universal family $\sC_g \ra \sT_g$ for an appropriate holomorphic map $\fie : C' \ra \sT_g$.
 
 If $ \ker (\rho)$ has finite index, then $C'$ is compact and, since Teichm\"uller space is a bounded domain in $\CC^{3g-3}$, the holomorphic map is constant.
 Therefore $S'$ is a product $ C' \times C_2$ and, denoting by $G : = \im (\rho)$, $S = (C' \times C_2)$, and we get exactly the surfaces 
 isogenous to a product such that  the action of $G$ on the curve $  C'$ is free. 
 
 If instead $G : = \im (\rho)$ is infinite, then the map of $C'$ into Teichm\"uller space is not constant, since the isotropy group of a point corresponding to a curve $F$ is, 
 as we saw, equal to the group
 of automorphisms of $F$ (which is finite). Therefore, in this case, we have a Kodaira fibration.
 
 \qed

\subsection{Singular fibres and mapping class group monodromy}
Let us consider again the situation where we have a fibration $f : S \ra B$ of an algebraic surface $S$ onto a curve of genus $b$,
and where the fibres have genus $g \geq 2$. We assume throughout that the fibration is relatively minimal and not a product.

We have devoted most of our attention to the consideration of the case where all the fibres $F$ of $f$ are smooth curves.
This can occur whenever the genus $b \geq 1$: just take  a surface isogenous to a product of unmixed type,
i.e., $ S = (C_1 \times C_2) / G$, where $G$ acts diagonally and freely on $C_1$, with quotient $C_1 / G = B$.

Then $f :  S = (C_1 \times C_2) / G \ra C _1 / G = B$ is the desired holomorphic bundle.

For Kodaira fibrations the situation is slightly more complicated, it is possible that the genus $b=2$, as we saw, but the fibre genus $g$  should be at least $3$.

At any rate, if we have a differentiable fibre bundle over a curve $B$ of genus $1$, then, $\tilde{B}$ being the universal covering of $B$, 
$\tilde{B} \cong \CC$,  the Torelli map $\tilde{B} \ra \HHH_g$
into Siegel's half space must be constant, therefore we
have a holomorphic bundle, which is then isogenous to a product, as described above.

The conclusion is that there must be at least one singular fibre if $b=1$ and we have nonconstant moduli.

On the other hand, when $b=0$, and the fibration is not a product, the number of singular fibres must be at least three.

Indeed, if $B^*= B \setminus \sC$ is the set of noncritical values of $f$, we want that there is a nonconstant map  $B^* \ra \HHH_g$;
if $B= \PP^1$, therefore the number of singular fibres, i.e. $|\sC|$, must at least three.

If we assume however that all the fibres are moduli semi-stable curves, then the number of singular fibres
must be at least  $4$ for $g \geq 1$ and at least $5$ for $g \geq 2$, as proven by Beauville and  Tan  \cite{beass} and \cite{tanss}
(\cite{zamora} gave an improvement for $g \geq 3$, that the number must be at least $6$).

Tan's inequality is $$ \frac{g}{2} ( 2b -2 +s) >  deg (f_* (\omega_{S|B})) = \chi(S) - (b-1)(g-1),$$
so that for $b=1$ the number $s$ of singular fibres is $> \frac{2}{g} $, which gives $s \geq 2$ only if $g=2$.

Can one  obtain better estimates on the number of singular fibres also when the genus of the base curve is $b=1$?

A very interesting result was obtained by Ishida \cite{ishida} who found, in the moduli space of surfaces with $q=p_g=1, K^2=3$,
described in joint work with Ciliberto (\cite{CaCi})
some surfaces whose Albanese fibration (onto a curve of genus $b=1$) has a unique singular fibre, and fibre genus $g=3$; but the singular fibre is not reduced.

As far as I know, there are no known examples of a non isotrivial fibration over a curve $B$ of genus $b=1$, and possessing
only one singular fibre, irreducible and with only nodes as singularities.

If such a fibration were to exist, the local monodromy around the singular fibres, a product of commuting Dehn twists, would be a commutator
in the mapping class group. This raises a general question 

\begin{question}
Which products of commuting Dehn twists are a commutator in the mapping class group $Map_g$, for $ g \geq 2$?
\end{question}

The question in genus $g=1$ is clear, since then the mapping class group is $SL(2, \ZZ)$,
and there the transvection (Picard-Lefschetz transformation) 
$$ T_1: e_1 \mapsto e_1, e_2 \mapsto e_1 + e_2$$
is not a commutator, and indeed no parabolic transformation in 
$$\PP SL(2, \ZZ) = \ZZ/2 * \ZZ/3 = \langle A,B | A^2 = B^3 = 1\rangle$$
is a commutator. The easiest way to see this is to express $T_1$ as a product of the two generators,
$T_1^{-1}  = AB$, so that the  image of $T_1$ in the Abelianization  $\ZZ/2 \times  \ZZ/3$
is the element $(1,-1)$. In fact, no power of $T_1$ is a commutator in  $SL(2, \ZZ)$.

In higher genus there is a surjection $Map_g \to Sp(2g, \ZZ)$ and there the obstructions seem to cease to exist, as shown by Corvaja and Zannier
(personal communication). Writing $T_2$ for the transvection such that
$$ T_2: e_1 \mapsto e_1, e_2 \mapsto e_1 + e_2, e_3 \mapsto e_3, e_4 \mapsto e_3 + e_4,$$
they show that $T_2$ is a commutator in $Sp(4, \ZZ)$, and that $T_1$ and the similarly defined $T_3$ are commutators in $Sp(6, \ZZ)$.

 In the mapping class group the question becomes more subtle: first of all, there are
 Dehn twists  on nonseparating and on separating curves; in the latter case the homology class of the curve is trivial, 
 and the image of the Dehn twist in $Sp(2g, \ZZ)$ is trivial.
 
 Endo and Kotschick \cite{endokot} showed that in the separating case the Dehn twist cannot be a commutator:
 the idea is that otherwise one would have a symplectic fibration over a torus with only one singular fibre with
 the Dehn twist as local monodromy.
 
 Cartwright and Steger (see \cite{4})  recently constructed via computer calculations a surface of general type with $ q=p_g=1$ and $K^2 =9$
 (hence a ball quotient). In this case the Albanese map \cite{catkeum} has either three irreducible singular fibres with one node,
 or one irreducible singular fibre having either  three nodes, or  a tacnode.
 
 The existence of this surface shows that the product of three Dehn twists can be a commutator in $Map_g$,
 and contradicts \cite{kot2}. Recently Stipsicz and Yun \cite{SY} announced the following result: a product of at most two Dehn twists cannot be a commutator in 
$Map_g$  for $g \geq 1$.

\section{Covers branched over line configurations}
\subsection{Generalities on Abelian coverings of the plane}

Assume that in  the projective plane  $\PP : = \PP^2$  are given distinct  irreducible curves $C_j$ of respective degrees $m_j$; we shall denote by $C$ the union of the curves $C_j$.

Then an irreducible  Abelian cover branched over the union $C$ of these curves and with group $G$ is determined by: 

\begin{itemize}
\item
elements 
$g_j \in G$   such that   they generate $G$ and such that   

\item the sum of the $ m_j g_j $ is equal to zero in $G$:

$$\sum_j  m_j g_j = 0 \in G.$$ 
\end{itemize}

In fact, one has an exact sequence $$H^2 (\PP, \ZZ) \ra H^2 (C, \ZZ) \ra H^3 (\PP,C,  \ZZ) \ra 0$$
and the last group is isomorphic to $H_1(\PP \setminus C, \ZZ) $
via Lefschetz duality, while  $ H^2 (C, \ZZ) \cong  \oplus_j \ZZ C_j$. This means that  $H_1(\PP \setminus C, \ZZ) $ is generated by loops around the curves $C_j$,
but these satisfy one relation.

The   relation has the above form because, denoting  $ r : H^2 (\PP, \ZZ) \ra \oplus_j \ZZ C_j$, the image of 
a line $L$  in $coker (r) = H_1(\PP \setminus C, \ZZ) $  is zero.

If we blow up some points $P_i $  which are singular for $C$, obtaining a birational morphism $\pi : Y \ra \PP$, then the homology of the complement of the preimage of $C$,
$H_1(Y \setminus \pi^{-1} C, \ZZ) $,  is generated by the loops around the curves
$D_j $ which are  strict transform of the curves $C_j$, and by  the loops around the exceptional curves $E_i$.

The monodromy $\mu$ of the Abelian cover takes these generators   to respective elements   $g_j$ and $\e_i $ of the group.

 Recall the relation 
$\sum_j  m_j g_j = 0$. 

If we write $$   D_j = m_j L  - \sum_{i}  a_{j,i} E_i  = 0$$

we can repeat the same argument with Lefschetz duality, and 
the relation saying that the image of $E_i$  in $coker (r_Y)$  is zero (here $ r_Y : H^2 (Y, \ZZ) \ra \oplus_j \ZZ D_j  \bigoplus \oplus_i \ZZ E_i $) yields

$$\sum_{j}  a_{j,i} g_j   =   \e_i.$$

This  third formula  determines the image $\e_i$  of a loop around $E_i$ under the monodromy homomorphism $\mu$ .

\medskip

\begin{defin}
{\bf (DEFINITION OF A MAXIMAL COVER )}

Let $d_j $ be the order of the element $g_j $ and consider the following group $G''$ (relating to the terminology used e.g. in the  lecture notes \cite{cime03},
$G'' $ is the abelianization of the orbifold fundamental group of the cover), defined as:

$G'' $ is  the quotient of the direct sum of  $\ZZ / d_j \ZZ $ by the relation  $$\sum_j  m_j g''_j = 0$$

where $g_j'' $ is the standard generator of the summand  $G''_j = \ZZ / d_j \ZZ$ .

Clearly the monodromy $\mu$ factors through $\mu'' $ , which sends the loop around $C_j$ (resp. $D_j$) to $g_j'',$
and the obvious surjection of $G'' \ra G$.

We get corresponding (irreducible) normal coverings   $ Z'' , Z $ of the plane $ \PP $ such that   $$ \PP =  Z'' / G''  = Z / G \ , \ Z = Z'' / H,$$
where $H$  is the kernel of the surjection  $G'' \ra G.$

A cover is said to be {\bf maximal} if $ G '' = G $, i.e., $Z=Z''$.
\end{defin}

\begin{rem}.  The quotient  $Z = Z'' / H$  is only ramified in a finite set.
\end{rem}

The importance of the concept of maximal covering is the following: let $Y$ be the surface obtained  by  blowing up  the points $P_i $
where $C$ is not a normal crossing divisor,
and assume that  the divisor $D$ in $Y$,  union of the $E_i $'s and the $D_j $'s,  is a normal crossing divisor  (this happens if
and only if  $C$ has only ordinary singularities).

In this case 
the local monodromies are  the elements  $g_j'' $, respectively 

$$ \e''_i =   \sum_{j}  a_{j,i} g''_j .$$ 

Now, we can write our covering $X \ra Y$ (the normalized  fibre product $Z \times _{\PP^2} Y$) as a GLOBAL QUOTIENT    $X = X'' / H$, 
and since the covering $X'' \ra X$ is only ramified in a finite set,
we get that, if $X''$ is smooth, then   $X$ has only cyclic quotient singularities (this approach has the advantage of making the description of the singularities
shorter).

It is therefore important to see general conditions which ensure the smoothness of $X''$: this is however technical, so we skip this analysis here.

As done by Hirzebruch (\cite{hirzlines}, \cite{hirzbook}) formulae simplify drastically if we require all the curves $C_j$ to be lines $L_j$, and we let all the
orders $d_j$ of the elements $g_j$ to be equal to the same integer $n$. 

In this case the maximal cover is called  the Kummer cover of exponent $n$ of the plane branched on the $r$ lines $L_1, \dots , L_r$ and its Galois group is the 
group $$(\ZZ/n)^{r-1} = (\ZZ/n)^{r} / \ZZ e, $$ where $e : = \sum_i e_i$.

\begin{defin}\label{HK}
Let $\sC$ be a configuration of distinct lines 
$$L_1 = \{ l_1 (x_1,x_2,x_3) = 0\}, \dots, L_r =  \{ l_r (x_1,x_2,x_3) = 0\}\subset \PP^2,$$
where we assume that there is no point in the plane $P$ belonging to each line $L_i$. Hence, without loss of generality, we 
can assume, after changing the numbering of the lines, and after a projective change of coordinates, that $L_i = \{x_i=0\}$
for $i=1,2,3$.

The linear forms $(l_1, \dots, l_r)$ yield an embedding $ l : \PP^2 \ra \PP^{r-1}$, and let  $(y_1, \dots, y_r)$ be coordinates in $\PP^{r-1}$.
Consider next the Galois covering $$\psi_n : \PP^{r-1} \ra \PP^{r-1}, \ \   \psi_n  ((z_1, \dots, z_r)) = (z_1^n, \dots, z_r^n),$$
with Galois group $(\ZZ/n)^{r-1} $. 

Let $Y$ be the fibre product of $L$ and $\psi$: 
$$ Y = \{ (x,z) \in  \PP^2 \times \PP^{r-1} | l (x) = \psi (z) \}.$$

Under our assumption on the linear forms, $Y$ indeed embeds in $ \PP^{r-1}$ as the complete intersection of $r-3$ hypersurfaces:
$$  Y  = \{ (z) |  z_j^n = l_j (z_1^n, z_2^n, z_3^n), \ j = 4, \dots, r\}.$$ 

The minimal resolution of singularities $X$ of $Y$ is called the {\bf Hirzebruch-Kummer covering}  of $\PP^2$ of exponent $n$
branched on the configuration $\sC$ of lines, and denoted $ HK_{\sC} (n)$.

Quite similarly one defines the {\bf Hirzebruch-Kummer covering}  of $\PP^m$ of exponent $n$
branched on a configuration $\sC$ of hyperplanes.
\end{defin}

\begin{rem}\label{Fermat}
a) When $r = m+2$ one gets a hypersurface, the Fermat hypersurface.

b) In the case of $\PP^1$ one gets curves which are also called {\bf generalized Fermat curves}.

c) If the configuration is a normal crossing configuration, $Y$ is a smooth complete intersection, so it has a lot of deformations.

\end{rem}

\subsection{Invariants}

Assume that, as in the previous subsection, we have an Abelian cover $S \ra Y$ branched on a normal crossing divisor $D' : =  \sum_j D_j + \sum_iE_i$.
To each divisor is associated the cyclic subgroup  (inertia subgroup) generated by  the local monodromy: $<g_j>$ in the case of $D_j$, and 
$<\e_i>$ in the case of $E_i$. Let $d_j$ be the order of $g_j$ and $d'_i$ the order of $<\e_i>$.

The covering $S$  is smooth if and only if, at each intersection point of two branch curves, the corresponding two inertia subgroups yield a direct sum.
In this case also the maximal cover $S''$ is smooth, at least provided that  the order of $\e''_i$ in $G''$ equals $d'_i$;
 and $S''$ is an \'etale covering of $S$ with group $H$, such that $G''/ H \cong G$. 
 
 The Chern numbers of $S$ can be easily calculated, since $K_S$ is the pull back of a divisor with rational coefficients
 $$ K_S = p^* (K_Y + \sum_j (1- \frac{1}{d_j}) D_j + \sum_i  (1- \frac{1}{d'_i}) E_i ) = $$
 $$= p^* ( - 3 L  + \sum_j (1- \frac{1}{d_j}) (m_j L - \sum_i a_{j,i} E_i) + \sum_i  (2- \frac{1}{d'_i}) E_i ) .$$ 
 
 Whereas, for the Euler number one uses the fact that it is additive for a stratification with strata which are orientable (non compact) manifolds
 (of different dimensions).
 
The simplest case is the case of a Kummer covering
\footnote{ One can define in the utmost generality  the Kummer covering of exponent $n$ of a normal variety $Y$ branched on $B$ as the normal finite
covering associated to
the epimorphism $\pi_1( Y \setminus B) \ra H_1 ( Y \setminus B, \ZZ) \otimes_{\ZZ} \ZZ/n$ } of exponent $n$ branched on $r$ lines, where $m_j = 1$, $d_j = n$, and where we assume that also
$d'_i = n$, for all $ i=1, \dots, k$. Then, observing that $a_{j,i} \in \{0,1\}$ and writing $v_i : = \sum_j a_{j,i}$ for the valency of the point $p_i$,
we get:
 $$K_S  = p^* ( - 3 L  + \sum_j (1- \frac{1}{n}) ( L - \sum_i a_{j,i} E_i) + \sum_i  (2- \frac{1}{n}) E_i ) = $$
 $$p^* ( (- 3   + r (1- \frac{1}{n})  ) L  + \sum_i  (1 + (1- \frac{1}{n}) (1 - v_i) ) E_i ) .$$
 Whence $$K_S ^2 = n^{r-1} [(- 3   + r (1- \frac{1}{n})  )^2 - \sum_i  (1 + (1- \frac{1}{n}) (1 - v_i) ) 
 ^2]$$
 Observe that we only needed to blow up the points of valency $v_i \geq 3$.
 
 Whereas, to calculate the Euler number we write $$k + 3 = e (Y) =  e ( Y - D') + \sum_i e (E_i^*) + \sum_j e(D_j ^*) + N,$$
 where $^*$ denotes the  intersection of a component of $D'$ with the smooth locus of $D'$, and $N$ is the number of singular points of $D'$.
 
 Let  $\de$ be  the number of double points of $C= \sum_j L_j$, $$\de = \frac{1}{2} r(r-1)-  \frac{1}{2} \sum_i v_i(v_i-1). $$
 Then 
 $$k + 3 = e (Y) =  e ( Y - D') + (2 k - \sum_i v_i  +  2 r  - \sum_i v_i  - 2 \de) + ( \sum_i v_i + \de).$$
hence, writing $v : = \sum_i v_i$
$$k + 3 = e (Y) =  e ( Y - D') + (2 k - 2 v  +  2 r    - 2 \de) + ( v + \de).$$ 
 the Euler number of $S$  is then just equal to $$ e(S) =   n^{r-1} [ k + 3 - (1- \frac{1}{n}) (2 k - 2 v  +  2 r    - 2 \de)  - (1 -  \frac{1}{n^2})( v + \de)]. $$
 
 In order to calculate the irregularity of $S$, and in general also $q = h^1 (\hol_S) ,  p_g =  h^2 (\hol_S)$, the best method is (see \cite{bc}) to calculate explicitly the
 decomposition of $p_* \hol_S$ into eigensheaves, 
 $$  p_* \hol_S = \hol_Y \oplus ( \bigoplus _{\chi} \hol_Y (- L_{\chi})) .$$

\subsection{Covers branched on lines in general position}

Assume that we have an Abelian covering with $d_j = n$ $\forall j=1, \dots, r$, branched on the $r$ lines $L_1, \dots, L_r$.
We assume that the  $r$ lines are  in general position, this means that $C$ has only double points, hence we get exactly  $\de = \frac{1}{2} r(r-1)$ intersection points.
In this case the fundamental group of  $\PP^2 \setminus C$ is Abelian, free of rank $(r-1)$.

Then any such covering is a quotient of the Hirzebruch-Kummer covering of exponent $n$, and we can assume that $r \geq 4$,
otherwise our surface is either singular or equal to $\PP^2$.

We saw in remark \ref{Fermat} that for $r=4$ the Hirzebruch-Kummer covering yields the Fermat surface
of degree $n$  in $\PP^3$.

The simplest case is  $n=2$ and  $r=4$: the Hirzebruch-Kummer cover
is the smooth quadric 
$$\{ u^2 = x^2 + y ^2 + z^2\} \subset \PP^3,$$
 hence the Del Pezzo surface  $\PP^1 \times \PP^1$.
 
 The next case of five lines in general position yields as Hirzebruch-Kummer cover   a complete intersection $Y$ of type $(2,2)$ in $\PP^4$,
 hence again a Del Pezzo surface, of degree $4$. There is only one intermediate covering branched on the five lines: it is a 
 singular Del Pezzo surface of degree $2$.

The next case of six lines in general position is interesting, the Kummer covering is a smooth K3 surface, with a group $(\ZZ/2)^5$ of covering automorphisms,
whereas the double cover branched on the six lines is a K3 surface with $15$ nodes coming from the $\de=15$ double points;
 in the case where the six lines are tangent to a conic $Q$,  it is a 16-nodal Kummer surface   blown up in one node (the blown up node corresponds
  to one of the two components in the inverse image of $Q$).
  
  In general, in order to obtain a surface of general type, we must have at least $4$ lines, and moreover we need, in the case where $d_j = n \ \forall j$, that
  $ r (n-1) > 3 n \Leftrightarrow (r-3)n > r $.
    
  Hence the smallest such case is for $r=4$, $n=5$. 
   If the $4$ lines are in general position, as we already mentioned,  the Hirzebruch-Kummer covering is the Fermat quintic
  $$\{ u^5 = x^5 + y ^5 + z^5\} \subset \PP^3,$$
  and any intermediate free $\ZZ/5$ quotient yields a classical Godeaux surface. 
  
  \subsubsection{Pardini's surfaces}
  Even if the surfaces we have been playing  with right now  may a priori look rather uninteresting, and more fun than serious mathematics, a remarkable example in the 
  class of abelian coverings branched on lines in general position was found
  by Pardini (\cite{Pardinican}); it belongs to  the next case, with $r= 5$ lines, $ n=5$, and group $G = (\ZZ/5)^2$.
  
  The covering surfaces have, by the above formulae, $K_S^2 = 25, e(S) = 35$ hence, by Noether's formula $\chi(S) = \frac{1}{12} (25 + 35) = 5$.
  
  Since the fundamental group of the surface is $(\ZZ/5)^2$ (the Kummer cover being simply connected,  see \ref{HK}), the conclusion is
  that $q(S)=0, p_g(S) = 4$, and the canonical map maps to $\PP^3$.
  
  To make a long story short, recall that, when the canonical map of a surface of general type has  degree $\geq 2$, then only two alternatives
  are possible (see \cite{beacan}):
  \begin{enumerate}
  \item
  the canonical image of $S$, $\Sigma : = \Phi_{K_S} (S)$ has $p_g(S) = 0$, or 
  \item
    the canonical image of $S$, $\Sigma : = \Phi_{K_S} (S)$ is canonically embedded.
  
  \end{enumerate}
  Babbage claimed that only the first case should occur, then three essentially equivalent examples of case 2) were found by 
  Beauville, van der Geer-Zagier, and  the present author, (\cite{babbage}, \cite{vdgz},\cite{beacan}), where the degree
  of the canonical map  $ \Phi_{K_S}$ equals $2$.
  
  It is still an open question which is the highest possible value for $ deg (\Phi_{K_S})$ in case 2), but the world record is $5$,
  obtained by Pardini (and later, independently, by Tan \cite{tan}).
  
  To give a smooth covering with group $G = (\ZZ/5)^2$, branched over five lines in general position, it is equivalent
  to give the $5$ monodromy vectors $g_1, \dots, g_5 \in  (\ZZ/5)^2$, with the property that:
  
  i) $\sum_j g_j = 0$, and 
  
  ii) two distinct vectors $g_i$ and $g_j$ are  linearly independent (this is the condition for the smoothness of $S$).
  
  Pardini's choice can best be explained in the following abstract way: we take the five vectors to be the five points of  an affine line in $(\ZZ/5)^2$
  not passing through the origin!
  
  Now, conditions i) and ii) are satisfied trivially, and there remains to see the advantage of this clever choice.
  
  This rests on the fact that there is a homomorphism $\psi : G \ra \ZZ/5$ mapping all the five elements $g_j$ to $ 1 \in \ZZ/5$.
  To the homomorphism $\psi$ there corresponds an intermediate $\ZZ/5$-covering $Z$ of $\PP^2$,
  i.e. we have a factorization $ S \ra Z \ra \PP^2$ of the original covering.
  
  If $l_j(x_0,x_1,x_2) = 0$ is the equation of the line $L_j$, then we see right away that
   $$Z = \{ x_3^5 = \prod_{j=1}^5 l_j(x_0,x_1,x_2) \} \subset \PP^3.$$ 
  
   Now, $Z$ is a quintic in $\PP^3$, whose only singularities are the 10 points $x_3= l_i(x_0,x_1,x_2)= l_j(x_0,x_1,x_2) = 0$.
   These singularities are rational double points of type $A_4$, hence $Z$ is a canonical model with $p_g(Z)=4$.
   
   Since  $p_g(S)= p_g(Z)=4$, the conclusion is:
   
   \begin{theo}{\bf (Pardini)}
   There exist surfaces $S$ with $K^2_S = 25$, $K_S$ ample, $p_g(S) = 4, q(S) =0$, such that the 
   canonical map of $S$ maps with degree $5$ onto a canonically embedded surface  $Z$.
   \end{theo}
   
   One can ask what happens if the five monodromy vectors  $g_1, \dots, g_5$ are affinely independent,
   yet satisfy i) and ii). This is interesting, since it shows the usefulness of calculating the character sheaves.
   
   We write a character $\chi : G \ra \ZZ/5$ as a pair $(a,b), 0 \leq a,b, < 5$, so that, for $(x,y) \in  (\ZZ/5)^2$,   $\chi (x,y) = ax + by \in \ZZ/5$.
   
   The character sheaves are of the form $\hol_{\PP^2} ( - L_{\chi})$, where the $L_{\chi}$ are calculated applying $\chi$ to the monodromy vectors,
   taking  the remainder for division by $5$    ( we denote by $[d]$   the remainder for division of $d$ by $5$,  hence $0 \leq [d] < 5$),
   and then summing all these remainders.
   
   Just in order to have a concrete example, take the five vectors $$(1,0), (0,1) , (1,1), (2, -1), (1, -1):$$ then 
   $$  5  L_{\chi} = a + b + [a+b] + [2a-b] + [a-b].$$

   This is important, because one can write 
   $$  H^0 (S, \hol_S (K_S) = \oplus_{\chi} z_{\chi} H^0 (\PP^2, \hol_{\PP^2} ( -3 +  L_{\chi}),$$ 
   where $z_j = 0$ is the equation of the ramification divisor $R_j$ corresponding to the j-th line $L_j$,
   and $ z_{\chi} : = \prod_{\chi(j) = 0} z_j^{5-1 -[\chi(g_j)]}$.

   In our example, $p_g=4$ and there are exactly $4$ characters such that  $L_{\chi} = 3$,
   namely, $a=b=4 , a=4, b=0, a=3, b=4, a=1, b=3$.
   
   One sees therefore  that the base locus of the canonical system $K_S$ equals the subscheme intersection of the four divisors
   $$ R_3 + 4 R_5, R_1 + 2 R_3 + 2 R_4, 4 R_2 + R_4, 3 R_1 + R_2 + R_5, $$ 
   that is, we get  the points $R_2 \cap R_3$ and $R_4 \cap R_5$ with multiplicity $1$.
   
   Hence
   \begin{prop}
   A $(\ZZ/5)^2$-covering of the plane branched on $5$ lines in general position, and with monodromy vectors 
   $$(1,0), (0,1), (1,1), (2, -1), (1, -1)$$ is smooth with $K^2_S = 25$,  $p_g(S) = 4, q(S) =0$,  has a canonical system with 2 simple base points, 
   and its canonical map is  birational onto a surface $\Sigma$ of degree $23$ in $\PP^3$.
   \end{prop}
   
  We remark that there are other  Abelian covers branched on $5$ lines in general position, which are {\bf uniform} (i.e. with $d_j = n \ \forall j=1, \dots, r$)
  and with exponent  $n=5$, but we do not pursue their classification here.
  
\subsection{Hirzebruch's and other ball quotients}

One obtains an infinite fundamental group whenever we take a  {\bf uniform} (i.e. with $d_j = n \ \forall j=1, \dots, r$) Abelian cover
branched on a union of lines, such that there exists a 
point $p_i$  of valency $v_i = : w \geq 3$. 

This is because projection with centre  the point $p_i$ leads to a fibration  $S'' \ra B$, 
where $B$ is an Abelian covering of $\PP^1$ branched on $w$ points,
and with group $(\ZZ/n)^{w-1}$. 

The genus of $B$ satisfies  $2 (b-1) = n^{w-1} (-2 + w (1-\frac{1}{n})),$
hence $b > 1$ as soon as $(w-2)(n-1) > 2$, which we shall now assume, and which holds for $n \geq 4$.

Since the maximal covering $S''$ maps onto a curve with genus $b \geq 2$, the fundamental group of $S''$ is infinite,
as well as the one of $S$, of which $S''$ is a finite \'etale covering.

Yet, it may still happen that the first Betti number of $S$ is zero, as we shall see.

Hirzebruch \cite{hirzlines} found explicit examples of ball quotients by taking Kummer coverings branched on a union of lines.

The one on which we shall concentrate mostly  is the one obtained for $n=5$  taking $6$ lines which are the sides of a complete quadrangle
(which can be also visualized as the sides of a triangle plus its three medians, which   meet in  the barycentre).

Other important examples were the Hirzebruch-Kummer coverings (cf. \ref{HK})  $HK_{\sH}(5)$ and $ HK_{\sD\sH}(3)$
associated to the Hesse configuration of lines  $\sH$, a configuration of type $(9_4,12_3)$ formed by the $9$ flexes of a plane cubic curve,
and the twelve lines joining pairs of flexpoints, respectively to  the dual Hesse configuration $\sD\sH$ of type $(12_3, 9_4)$
of the nine lines dual to the flexpoints.

Recall that every smooth plane cubic curve is isomorphic to one in the Hesse pencil:

$$ C_{\la} : = \{  x^3 + y^3 + z^3 + 6 \la xyz = 0 \}. $$ 
$C_{\la}$ is smooth, except for $\la = \infty$, or $ 8 \la^3= -1$; its flexes are the intersection of $C_{\la}$
with its Hessian cubic curve, which is precisely 
$$ H_{\la} : = \{  \la^2 (x^3 + y^3 + z^3) - (1 + 2 \la^3)  xyz = 0 \} = C_{\mu}, \ \mu = - \frac{1 + 2 \la^3}{6 \la^2}. $$

Hence the nine points are the base points of the pencil, 
$$ \{ xyz= x^3 + y^3 + z^3 = 0 \},$$
while $C_{\la} = H_{\la} $ exactly when $ 8 \la^3= - 1$: hence the $4$ singular curves of the pencil are
$4$ triangles (every point is a flexpoint!), and these $4$ triangles produce the $12$ lines.

More examples of ball quotients can be found in \cite{hirzRUSS} and \cite{hirzbook}
(they are also related to hypergeometric integrals, see \cite{DeligneMostow}).

\subsection{Symmetries of the del Pezzo surface of degree $5$}

The blow-up $Y$ of  $4$ points $P_1, \dots, P_4 \in \PP^2$  in general position 
is the del Pezzo surface of degree $5$.

This surface is the moduli space of ordered quintuples of points in $\PP^1$,
as we shall now see.
The six lines can be labelled $L_{i,j}$, with $i\neq j \in \{1, 2,3,4\}$ ($L_{i,j}$ is the line $\overline{P_i P_j}$).

The following is well known

\begin{theo}
The automorphism group of the Del Pezzo surface $Y$  of degree $5$ is isomorphic to $\mathfrak S_5$.
\end{theo}

\Proof
There is an obvious action of the symmetric group $\mathfrak S_4$ permuting the $4$ points, but indeed there is more
(hidden) symmetry, by the  symmetric group $\mathfrak S_5$. This can be seen denoting by $E_{i,5}$ the exceptional curve lying over
the point $p_i$, and denoting, for $i\neq j \in \{1, 2,3,4\}$, by $E_{i,j}$ the line $L_{h,k}$, if $\{1, 2,3,4\}=  \{i, j, h,k\}$.

For each choice of $3$ of the four points, $\{1, 2,3,4\} \setminus \{h\}$, consider the standard Cremona transformation $\s_h$
based on these three points. To $\s_h$ we associate the transposition $(h,5) \in \mathfrak S_5$, and the upshot is that
$\s_h$ transforms the $10$ $(-1)$ curves $E_{i,j}$ via the action of $(h,5) $ on pairs of elements in $\{1, 2,3,4,5\}$. 

Indeed there are  five geometric objects permuted by $\mathfrak S_5$: namely, $5$ fibrations $\varphi_i : Y \ra \PP^1$,
induced, for $1 \leq i \leq 4$, by the projection with centre $P_i$, and, for $i=5$, by the pencil of conics through
the $4$ points. Each fibration  is a conic bundle, with exactly three singular fibres, corresponding to the possible partitions of
type $(2,2)$ of the set $\{1, 2,3,4,5\} \setminus \{i\}$.

To conclude that $\mathfrak S_5 = \Aut (Y)$,  we observe that $Y$ contains exactly ten lines, i.e. 
irreducible curves $E$ with $ E^2 = E K_Y = -1$.  We have the following easy result:

\begin{lemma}
The curves $E_{i,j}$, which generate the Picard group, have an intersection behaviour which
is dictated by the simple rule (recall that $ E_{i,j}^2 = -1, \ \forall i \neq j$)
$$ E_{i,j} \cdot E_{h,k} = 1 \ \Leftrightarrow \{i,j\} \cap \{h,k\}= \emptyset, \ \ E_{i,j} \cdot E_{h,k} = 0   \Leftrightarrow \{i,j\} \cap \{h,k\}\neq  \emptyset.$$
In this picture the three singular fibres of $\varphi_1$ are $$E_{3,4} +  E_{2,5}, \ E_{2,4} +  E_{3,5}, \ E_{2,3} +  E_{4,5}.$$

The relations among the $E_{i,j}$'s in the Picard group come from the linear equivalences 
$E_{3,4} +  E_{2,5} \equiv E_{2,4} +  E_{3,5} \equiv  E_{2,3} +  E_{4,5}$ and their $\mathfrak S_5$-orbits.

\end{lemma}

Therefore each automorphism $\psi$ of $Y$ permutes the 10 lines, preserving the incidence relation.
Up to multiplying $\psi$ with an element of the subgroup $\mathfrak S_5$, we may assume that $\psi$ fixes $E_{1,2}$, hence that $\psi$ 
permutes the three curves $E_{3,4}, E_{3,5}, E_{4,5}$. By the same trick we may assume that 
$\psi$ fixes $E_{1,2}$,$E_{3,4}, E_{3,5}, E_{4,5}$. Hence $\psi (E_{1,j}), 3 \leq j \leq 5$ is either 
$E_{1,j}$ or $E_{2,j}$. Multiplying $\psi$ by the transpostion $(1,2)$ we may assume that also
 $\psi (E_{1,3}) = E_{1,3}$. The incidence relation now says that $\psi$ fixes all the ten lines,
 and there remains to show that $\psi$ is the identity. But blowing down the four curve $E_{i,5}$
 we see that $\psi$ indudes an automorphism of $\PP^2$ fixing each of the points $P_1, P_2, P_3, P_4$.
 Hence $\psi$ is the identity and we are done.

\qed

\begin{rem}
We have the following correspondences:
\begin{itemize}
\item
The lines $E_{i,j}$ correspond to the transpositions in $\mathfrak S_5$.
\item
The 15 intersection points $ E_{i,j} \cdot E_{h,k} $, for $| \{i.j.h.k\}| = 4$
correspond to isomorphisms of $(\ZZ/2)^2$ with a subgroup of $\mathfrak S_5$.
\item
The five  2-Sylow subgroups correspond to the five conic bundles $\varphi_i$, $i=1, \dots, 5$,
and the triples of singular fibres correspond to the three different embeddings 
$(\ZZ/2)^2 \ra \mathfrak S_5$ with the same image.
\item
The six 5-Sylow subgroups correspond combinatorially to pairs of opposite {\bf pentagons}.
Here a pentagon is the equivalence class of a bijection 
$ \sP : \ZZ/5 \ra \{1,2,3,4,5\}$ for the action of the dihedral group $D_5$ 
on the source ( $ n \in \ZZ/5 \mapsto \pm n + b, b \in \ZZ/5$). Whereas a pair of opposite 
pentagons is the equivalence class for the action of the affine group $A(1, \ZZ/5)$
on the source.
\item
To a combinatorial pentagon corresponds   a {\bf geometric pentagon}, 
i.e. a union of lines $E_{\sP(i), \sP(i+1)}, i \in \ZZ/5$ each meeting the 
following line $E_{\sP(i+2), \sP(i+3)}$. The divisor of a geometric pentagon
is an anticanonical divisor $D_{\sP}$, and the sum   of two opposite 
geometric pentagons $D_{\sP}+ D_{\sP^o}$
is just the sum of the 12 lines $E_{i,j}$.
\item
If $D_{\sP} = div (s_{\sP}),  s_{\sP} \in H^0(\hol_Y (-K_Y))$, we obtain 
five independent quadratic equations for the anticanonical embedding 
$$ Y \ra \PP (H^0(\hol_Y (-K_Y))^{\vee},$$
from the six equations $s_{\sP} s_{\sP^o} = \delta$,
where $ div(\de) = \sum E_{i,j}$.
\item
The above symmetry is the projective icosahedral symmetry,
i.e., the symmetry of the image of the icosahedron in $\PP^2(\RR)$:
the 10 lines correspond to pairs of opposite faces,
the 15 points to pairs of opposite edges,
the 6 pairs of opposite pentagons correspond to pairs 
of opposite vertices.

\end{itemize}

\end{rem}

A basis for $ H^0(\hol_Y (-K_Y))$ is given by the six sections 
corresponding to the pentagons where $4,5$ are never neighbours.

Written as sections of $ H^0(\hol_{\PP^2} ( 3))$ vanishing at the points $P_1, \dots, P_4$
which we assume to be the coordinate points and the point $(1,1,1)$, they are:

$$ s_{i,j} = x_i x_j (x_j - x_k), \  \{i,j,k\} = \{ 1,2,3\}. $$ 

All this leads to beautiful Pfaffian equations for $Y \subset \PP^5$, which are 
$\mathfrak S_5$ equivariant.
Indeed
$ H^0(\hol_Y (-K_Y))$ is the unique irreducible representation $\chi_6$ of $\mathfrak S_5$
of dimension 6, whereas the representation $\rho$ on the set of pairs of opposite pentagons splits as the 
direct sum of the trivial representation with an irreducible representation $\chi_5$
of dimension $5$, the representation $\tilde {\rho}$ on the pentagons
is the direct sum of $\rho$ with the tensor product of $\rho$ with the signature character $\chi_1$.
The permutation representation of  $\mathfrak S_5$ on $ \{ 1,2,3,4,5\}$ splits as the trivial representation
direct sum with an irreducible representation $\chi_4$ of dimension $4$.
$\chi_6$ is the only  irreducible representation such that $\chi_6 \cong \chi_6 \otimes \chi_1 $.

In this way one obtains all the irreducible representations of $\mathfrak S_5$, see \cite{JL}, page 201.
We shall not go further here with the equations of $Y$, for the anticanonical embedding 
$Y \subset \PP^5$, and for the embedding $Y \subset (\PP^1)^5$ via $\varphi_1 \times \dots \varphi_5$.

But we shall now prove the classical:

\begin{theo}
The Del Pezzo surface $Y$ of degree 5 is the moduli space for ordered 5-uples of points in $\PP^1$,
i.e., the GIT quotient of $(\PP^1)^5$ by $\PP GL(2, \CC)$.

\end{theo}
\Proof

Another model for $Y$ is the blow up of $\PP^1 \times \PP^1$ in the three diagonal points $(\infty, \infty), (0,0), (1,1)$,
and the ten lines come from the three blown up points,
plus the strict transforms of  the diagonal and of the vertical and horizontal lines 
$x=0, x=1, x= \infty, y=0, y=1, y= \infty$.

Removing these seven lines in $\PP^1 \times \PP^1$ we obtain a point $(u,v) $ such that
the five points $\infty, 0,1, u,v$ are all distinct.

If we approach a smooth point in the diagonal line, say $(u,u) $ we obtain the 5-uple
$\infty, 0,1, u,u$ where  the  fourth and the fifth  points are equal $P_4 = P_5$, and the other three
are different (so that set theoretically we have four distinct points).
By $\mathfrak S_5$-symmetry, the same occurs whenever we get a smooth point
of the divisor $ \sum E_{i,j} \subset Y$: $P_i = P_j$ and the other three are different from $P_i$
and pairwise different.

If we tend to the point $(0,1) $, we get the 5-uple $\infty, 0,1, 0,1$,
where $P_2 = P_4$, $P_3 = P_5$; again by symmetry, to the point
$E_{i,j} \cap E_{h,k} $ corresponds a quintuple with $P_i = P_j$, $P_h = P_k$,
and where set theoretically we have three distinct points.

Now, as shown in Mumford's book \cite{git}, especially proposition 3.4, page 73 , in this case the semistable 5-uples
are stable, and a quintuple is unstable if and only if three points are equal.
One can easily conclude that we have an isomorphism of $Y$ with the GIT-quotient 
$(\PP^1)^5$ // $ \PP GL(2, \CC)$.

\qed

\subsubsection{Hirzebruch-Kummer covers of a del Pezzo surface of degree $5$}

 Let us now take as branch locus $D: = \sum_{i,j} E_{i,j} $, and let us notice that $D$  is linearly equivalent to twice 
  the anticanonical divisor $-K_Y$. 
 
 Therefore, for any $n$- uniform covering $S$ of $Y$ branched on the $10$ lines, and which is smooth, 
 $$ K_S = p^* ( K_Y -2  (1-\frac{1}{n}) K_Y)= p^* ( - (\frac{n-2}{n}) K_Y ) .$$
 
 In particular, for the Kummer covering, one has
 $$ K_{S''}^2 = 5 (n-2)^2 n^3. $$
 
 Whereas the Euler number can be calculated  (in \cite{cd3} it is calculated in an elegant way) as:
 $$ e(S'') = n^3(  3 + 2(n-2)(n-3)) = n^3 (2 n^2 - 10 n + 15). $$
 
 The Chern slope equals then
 $$ \nu_C (S'') =  \frac{5 (n-2)^2}{ 3 + 2(n-2)(n-3)} = \frac{5}{2} \frac{n-2}{n-3 +  \frac{3}{2(n-2)}}  > \frac{5}{2}.$$
 The same formula shows that the slope is a decreasing function of $n$ for $n \geq 5$, and for $n=5$ we obtain 
 $ \nu_C (S'') = \frac{45}{15} = 3.$ Hence the theorem of Hirzebruch \cite{hirzlines}:
 
 \begin{theo}
 A smooth $n$-uniform Abelian cover $S$ with $n \geq 5$ branched on the $10$ lines of a del Pezzo surface of degree $5$
 is a surface with ample canonical divisor, with  positive index, and indeed a ball quotient if and only if $n=5$.
 \end{theo}
 
 We observe, as a side remark,   that the Kummer coverings above embed into $C_n^5$, where 
 $C_n$ is the Fermat curve of degree $n$, $$C_n = \{ x^n + y^n + z^n = 0\} \subset \PP^2,$$
 via the Cartesion product of the Stein factorizations of the maps induced by the $\varphi_i,\  i=1, \dots, 5$.
 
 We shall next discuss some of the surfaces $S$ mentioned in the previous theorem, addressing
  the question of their irregularity $q(S)$.
 
 \subsubsection{A problem posed by Enriques, and its partial solution}
 
 Enriques posed in his book \cite{enr} the following problems:
 
 \begin{question}
 Given a surface $S$ with $p_g (S) = 4$ and with birational canonical map onto its image $\Sigma \subset \PP^3$,
 
  I) what is the maximum value  for its canonical degree $K_S^2$?
  
  II) what is the maximum value for $deg (\Sigma ) $?
 
 \end{question}
 He indeed suspected that there should be an  upper bound equal to $24$, for which counterexamples were given in \cite{bidouble}.
 
 Now, without loss of generality, we may assume that $S$ is minimal, because otherwise $K_S^2$ decreases, and we observe that, since $K_S$ is nef, we have
 $$   45 \geq K_S^2 \geq deg(\Phi_{K_S}) \cdot  deg (\Sigma ),$$
 where the first inequality is a consequence of the Miyaoka-Yau inequality $ K_S^2 \leq 9 (1-q(S) + p_g(S) ) = 9 (5 - q(S))$.
 
 Hence, in order to achieve the equality $K_S^2 = 45$, when $p_g (S) = 4$, we must have a ball quotient which has $q(S)=0$.
 
 With I. Bauer \cite{bc} we showed that there exists such a surface.
 
 \begin{theo}
 There exists an Abelian cover of the del Pezzo surface $Y$ of degree $5$, with group $(\ZZ/5)^2$, and branched on the $10$ lines, which is regular,
 i.e. $q(S)=0$,  has $p_g (S) = 4$, $K_S^2 = 45$, and has a  birational canonical map onto a surface $\Sigma$ of degree $19$.

 \end{theo}
 
 Indeed in the course of the search we classified all such coverings of the del Pezzo surface $Y$ of degree $5$, with group $(\ZZ/5)^2$, and branched on the $10$
 lines.

We considered  the group $\mathcal{G}$,  generated by $\mathfrak{S}_5$ and $Gl(2, 
\mathbb{Z}/ 5 \mathbb{Z})$, acting on the set of 
admissible  monodromy vectors.

  A MAGMA computation
showed that $\mathcal{G}$ has four orbits, and 
representatives for these orbits could be taken as:
$$
\mathfrak{U}_1 = (
\begin{pmatrix}
1 \\
0
\end{pmatrix}, \begin{pmatrix}
1 \\
0
\end{pmatrix}, \begin{pmatrix}
0 \\
1
\end{pmatrix}, \begin{pmatrix}
2 \\
1
\end{pmatrix}, \begin{pmatrix}
2 \\
1
\end{pmatrix}, \begin{pmatrix}
4 \\
2
\end{pmatrix});
$$

$$
\mathfrak{U}_2 = (
\begin{pmatrix}
1 \\
0
\end{pmatrix}, \begin{pmatrix}
1 \\
0
\end{pmatrix}, \begin{pmatrix}
0 \\
1
\end{pmatrix}, \begin{pmatrix}
2 \\
1
\end{pmatrix}, \begin{pmatrix}
4 \\
2
\end{pmatrix}, \begin{pmatrix}
2 \\
1
\end{pmatrix});
$$

$$
\mathfrak{U}_3 = (
\begin{pmatrix}
1 \\
0
\end{pmatrix}, \begin{pmatrix}
1 \\
0
\end{pmatrix}, \begin{pmatrix}
0 \\
1
\end{pmatrix}, \begin{pmatrix}
4 \\
1
\end{pmatrix}, \begin{pmatrix}
3 \\
2
\end{pmatrix}, \begin{pmatrix}
1 \\
1
\end{pmatrix});
$$

$$
\mathfrak{U}_4 = (
\begin{pmatrix}
1 \\
0
\end{pmatrix}, \begin{pmatrix}
1 \\
0
\end{pmatrix}, \begin{pmatrix}
0 \\
1
\end{pmatrix}, \begin{pmatrix}
1 \\
1
\end{pmatrix}, \begin{pmatrix}
0 \\
3
\end{pmatrix}, \begin{pmatrix}
2 \\
0
\end{pmatrix}).
$$

In particular we saw that
$\mathcal{G} \cong
Gl(2, \mathbb{Z}/ 5 \mathbb{Z}) \times \mathfrak{S}_5$,
and concluded the classification with the following result:

\medskip

\begin{theo}\label{four}
Let $S_i$ be the minimal smooth surface of general type with $K^2 
=45$ and $\chi =5$ obtained
  from the covering induced by the admissible six - tuple 
$\mathfrak{U}_i$, where $i \in \{1,2,3,4\}$. Then we
have that $S_3$ is regular (i.e., $q(S_3) = 0$), whereas $q(S_i) = 2$ 
for $i \neq 3$.
\end{theo}

The story concerning Enriques' question  is not yet completely finished, because it is not clear whether one can achieve $deg (\Sigma) = 45$
(the current published record \footnote{Recently we achieved $deg (\Sigma) = 32$.} in this direction is $deg (\Sigma) = 28$, \cite{bidouble}).

\subsection{Bogomolov-Miyaoka-Yau fails in positive characteristic}

Even if the BMY inequality was proven by Miyaoka with purely algebraic methods, still the proof uses characteristic $0$ arguments in an essential way.

Robert Easton \cite{easton} gave easy examples showing that indeed the BMY inequality does not hold in positive characteristic.
These examples are related to the Hirzebruch-Kummer coverings of the plane, and the main idea  is that in characteristic $p > 0$
there are configurations of lines which cannot exist in characteristic zero.

These configurations are just the projective planes $\PP^2_{\ZZ/p} \subset \PP^2_K$, for each algebraically closed field $K$ of characteristic $p>0$.

The easiest case is the so called {\bf Fano plane}  $\PP^2_{\ZZ/2} \subset \PP^2_K = : \PP^2.$
We have a configuration $\sC$ of type $7_3 7_3$, seven lines passing through seven points, each triple for the configuration.

For each odd number $n \geq 5$ we consider the Hirzebruch-Kummer covering of exponent $n$, $S_n$, a finite Galois cover of 
the blow up $Z$ of $\PP^2$ in the seven points, with Galois group $(\ZZ/n)^6$.

The canonical divisor of $Z$ is $ K_Z = - 3 H + \sum_1^7 E_i$, where $E_i $ is the inverse image of the point $P_i$. Denoting by
$L_i$ the proper transform of the line $L_i$, since $\sum_1^7 L_i = 7 H - 3  \sum_1^7 E_i $,
the canonical divisor of $S_n$ is the pull back of 
$$  K_Z  +  (1 -  \frac{1}{n} )( \sum_1^7 L_i +  \sum_1^7 E_i ) = (-3 + 7 (1 -  \frac{1}{n} ) ) H +  ( 1  - 2 (1 -  \frac{1}{n} ) ) \sum_1^7 E_i  $$

Hence   
$$   K^2_{S_n} =  n^6  [( 4  -  \frac{7}{n} )^2 - 7 ( (1 -  \frac{2}{n} )^2 ] = n^6 [ 9 - \frac{28}{n} + \frac{21}{n^2}].$$ 

The formulae we illustrated earlier yield

$$ c_2 (S_n) = n^6 [ 10  + 14   (1 -  \frac{1}{n} ) - 21  (1 -  \frac{1}{n^2} )] = n^6 [ 3  - 14   \frac{1}{n}  + 21 \frac{1}{n^2} ] .$$ 

Hence 
$$  K^2_{S_n} =  3 c_2 (S_n)  + n^6   ( 14  \frac{1}{n}  - 42  \frac{1}{n^2} ) = 3 c_2 (S_n)  + n^4   ( 14  n   - 42   ) .$$

Hence a particular case of Easton 's  theorem

\begin{theo}
In characteristic equal to $2$ the Hirzebruch-Kummer coverings of the plane branched on the Fano of configuration of lines have Chern slope
$$ \nu (S_n) =  \frac{ K^2_{S_n}} {c_2 (S_n) }=  3 ( 1 +  \frac{14 n - 42} {9 n^2 - 42 n + 63}) > 3,$$ 
violating the Bogomolov Miyaoka Yau inequality.
The maximum of the slope is attained for $n=5$, $ \nu (S_5) = 4 + \frac{3}{39}$.
\end{theo}

\section{Counterexamples to Fujita's semiampleness question, rigid manifolds}

\subsection{BCDH surfaces, counterexamples  to Fujita's question}
Recently, we found new counterexamples to Fujita's question \ref{fujitaquestion}(we say counterexamples since we heard through the grapevine
that experts were expecting a positive answer; however  this counterexample does not inficiate the abundance conjecture).

  \begin{theo}\label{surfaces}
  There exists an infinite series of surfaces with ample canonical bundle,  whose Albanese map is a  fibration $ f : S \ra B$ onto a curve $B$ of genus $b= \frac{1}{2}(n-1)$, and with fibres of genus $g = 2b= n-1$,  where $n$ is any integer relatively prime with $6$.
  
  These Albanese  fibrations  yield negative answers to Fujita's question about the semiampleness of $V : = f_* \om_{S|B}$, since
  here $V : = f_* \om_{S|B}$
splits as a direct sum 
$ V = A  \oplus Q$, where $A$ is an ample    vector bundle, and $Q$ is a  unitary  flat  bundle with
infinite monodromy group.

The fibration $f$ is semistable: indeed all the fibres are smooth, with the exception of three fibres which are the union of two smooth curves of genus $b$
which meet transversally in one point.

For $n=5$ we get   three surfaces which are rigid, and are quotient of the unit ball in $\CC^2$ by a torsion free cocompact lattice $\Ga$.
We shall call them BCD-surfaces (cf. theorem \ref{four}).
The rank of $A$, respectively $Q$, is in this case equal to $2$.
 \end{theo}

The easiest way to describe these surfaces, which are Abelian covers of the del Pezzo surface $Y$ of degree $5$ with group $(\ZZ/n)^2$,
branched over the $10$ lines of $Y$, is to look at a birational model which is an Abelian covering of $\PP^1 \times \PP^1$
branched over  the diagonal and of the vertical and horizontal lines 
$x=0, x=1, x= \infty, y=0, y=1, y= \infty$.

 We consider again the equation 
$$ z_1^n =  y_0 ^{m_0} y_1^{m_1}  (y_1 - y_0 )^{m_2}   (y_1 - x y_0 )^{m_3}  , \ \  x \in \CC \setminus \{0,1  \}$$
but we homogenize it to obtain the equation

$$ z_1^n = y_0 ^{m_0} y_1^{m_1}  (y_1 - y_0 )^{m_2}   (x_0 y_1 - x_1 y_0 )^{m_3}     x_0^{n-m_3}. \  $$

The above equation describes a singular surface $\Sigma'$ which is a cyclic covering of 
 $\PP^1 \times \PP^1$ with group $ G : = \ZZ / n$; $\Sigma'$ 
is contained inside the line bundle $\LL_1$  
over $\PP^1 \times \PP^1$  whose sheaf of holomorphic sections $\sL_1$ equals $ \hol_{\PP^1 \times \PP^1}(1,1)$.

The first  projection $\PP^1 \times \PP^1 \ra  \PP^1$ induces a morphism $p : \Sigma' \ra  \PP^1$
and we consider the curve $B$, normalization of the covering of $\PP^1$ given by
$$ w_1^n =  x_0 ^{n_0} x_1^{n_1}  (x_1 - x_0 )^{n_2}  .$$ 

We consider the normalization $\Sigma$ of the  fibre product $\Sigma' \times_{\PP^1} B$. 

$\Sigma$ is an abelian covering of $\PP^1 \times \PP^1$ with group $(\ZZ / n)^2$, and, if $GCD (n,6) = 1$, 
 for a convenient choice of
the integers $m_j , n_i$, for instance for 
$$  m_0= m_1=  m_2= 1, \ m_3 =  n-3,  n_0= n_1 = 1 , n_2 = n-2$$
we obtain a smooth $(\ZZ / n)^2$- Abelian covering of the blow up of $\PP^1 \times \PP^1$ in the three diagonal points $(\infty, \infty), (0,0), (1,1)$,
which is the del Pezzo surface $Y$ of degree $5$.

The corresponding singular fibres are only $3$, and come from one of the $5$ conic bundle structures on $Y$,
here given by the  first  projection $\PP^1 \times \PP^1 \ra  \PP^1$: hence one sees right away that the singular fibres are reducible, 
and that they are the union of 
two smooth curves of genus $b$ intersecting transversally in one point.

The surfaces in theorem \ref{surfaces}, which we shall call BCDH-surfaces, have an \'etale unramified covering given by
the Hirzebruch Kummer coverings  of the del Pezzo surface $Y$ of degree $5$ with group $(\ZZ/n)^2$,
branched over the $10$ lines of $Y$, which we shall denote HK-surfaces.

In the case $n=5$ we get ball quotients, and, for $5 | n$, we get therefore that the universal covering is a branched covering of the ball.
This motivates the following questions, which are clearly satisfied for $n=5$:

\begin{question}
\begin{enumerate}
\item
Are the BCDH-surfaces rigid?
\item
Is their universal covering $\tilde{S}$, or more generally the universal covering of the  HK-surfaces,  a Stein manifold? 
\item
Is the universal covering $\tilde{S}$ of the  BCDH- and HK-surfaces contractible?
\item
Do the  BCDH-surfaces admit a metric of negative curvature?

\end{enumerate}

\end{question}

The last question is motivated by analogy with the examples considered by Mostow-Siu: namely,
one has a K\"ahler-Einstein metric on $Y$ deprived of the $10$ lines, and
on the covering one should interpolate with another  metric localized on the ramification divisor.

Of course, an interesting question is whether the K\"ahler-Einstein metric $\omega_n$ on $S_n$ is negatively curved. 
To analyse the question one can observe that $\omega_n$  is the pull-back of a metric on $Y$ with given cone angles along the branch divisors.

Now, a positive answer to (4) would imply (3); however, for (3),  the vanishing cycles  criterion which was given earlier may be applied.

Question  (2) should have a positive answer, while for  the answer to question (1), it is positive,
as proven in joint work with I. Bauer \cite{bc-rigid}:

\begin{theo}\label{HKrigid}
The HK-surfaces are rigid for $ n \geq4$, hence also the BCDH-surfaces are rigid.

\end{theo}

Observe that, since the  BCDH-surfaces have the HK surfaces as \'etale covers, the non rigidity of the former would imply
non rigidity of the latter (see \ref{etale}).  Moreover, for $n=3$ rigidity does not hold, this surface was studied by Roulleau in \cite{Roulleau}.

The proof for  HK surfaces uses many tools, first of all  the special geometry of the del Pezzo surface of degree $5$, and rather complicated arguments
using logarithmic sheaves in order to control the deformations of the Abelian coverings.

\begin{cor}
 The Albanese fibrations of BCDH-surfaces yield rigid curves in the moduli spaces $\overline{\mathfrak M_{n-1}}$.
\end{cor}

We noticed that all the fibres of the Albanese map have compact Jacobian, hence the following:

\begin{question}
Do the  Albanese fibrations of BCDH-surfaces yield rigid curves in the moduli spaces $\mathfrak A_{n-1}$?

\end{question}

Theorem \ref{HKrigid} raises  several questions: here, given a configuration of lines $\sC$, we denote by
$HK _{\sC } (n)$ the Hirzebruch-Kummer covering of exponent $n$ ramified on the lines of the configuration $\sC$.

Natural questions are (see \cite{bc-rigid}:
\begin{question}
{\bf I)} For which rigid configuration $\sC$ of lines in $\PP^2$ is the associated Hirzebruch Kummer covering $ HK_{\sC}(n)$ rigid for $ n > > 0$? 

{\bf II)} For which rigid configuration $\sC$ of lines in $\PP^2$ is the associated Hirzebruch Kummer covering $ HK_{\sC}(n)$ a $K(\pi,1)$ for $ n > > 0$? 

{\bf III)} For which rigid configuration $\sC$ of lines in $\PP^2$ does the associated Hirzebruch Kummer covering $ HK_{\sC}(n)$ possess a K\"ahler metric
of negative sectional curvature for $ n > > 0$? 

\end{question}

Observe that if III) has a positive answer, then also II), by the Cartan-Hadamard theorem. 

Abelian coverings branched over configurations were also used by Vakil in \cite{murphy}, who showed  
  that for $n >> 0$
the local deformations of $ HK_{\sC}(n)$ correspond to the product of the deformations of the configuration 
$\sC$   with a smooth manifold. Vakil used a result of Mnev \cite{mnev} 
in order to show that, up to a product with a smooth manifold, one obtains all possible singularity types.

Let me end this subsection by commenting that not only BCD- and BCDH-surfaces are quite interesting from
many algebro-geometric and complex analytic points of view,
but that the features of their Albanese fibration  have also found  very interesting applications for the construction of remarkable
symplectic manifolds, in the work of Akhmedov and coworkers
(see {\cite{akh} and literature cited  therein).

\subsection{Rigid compact complex manifolds}

Recall the following notions of rigidity (see \cite{bc-rigid} for more details)

\begin{definition}\label{rigid} \
\begin{enumerate}
\item  Two compact complex manifolds $X$ and $X'$ are said to be {\em  deformation equivalent} if and only if there is a 
proper smooth holomorphic map  $$f \colon \mathfrak{X}  \rightarrow \sB 
$$ 
where $\sB$ is a connected (possibly not reduced) complex space and there are points $b_0, b_0' \in \sB$ such that the fibres $X_{b_0} : = f^{-1} (b_0),
X_{b'_0} : = f^{-1} (b'_0)$ are respectively isomorphic to $X, X'$ ($X_{b_0} \cong X, X_{b'_0} \cong X'$).

\item A compact complex manifold $X$ is said to be {\em globally rigid} if for any compact complex manifold $X'$, which is deformation equivalent to $X$, we have  an isomorphism $X \cong X'$.
\item  A compact complex manifold $X$ is instead said to be  {\em (locally) rigid} (or just {\em rigid}) if for each deformation of $X$,
$$f \colon (\mathfrak{X},X)  \rightarrow (\sB, b_0)
$$ 
there is an open neighbourhood $U \subset \sB$ of $b_0$ such that $X_t := f^{-1}(t) \cong X$ for all $t \in U$.
\item  A compact complex manifold $X$ is said to be  {\em infinitesimally rigid} if 
$$H^1(X, \Theta_X) = 0,$$
where $\Theta_X$ is the sheaf of holomorphic vector fields on $X$.
\item
$X$ is said to be  {\em strongly rigid} if the set of compact complex manifolds $Y$ which are homotopically equivalent to $X$,
$\{ Y |  Y  \sim_{h.e.} X\}$ consists of a finite set of isomorphism classes of globally rigid varieties.
\item
$X$ is said to be  {\em \'etale  rigid} if every \'etale (finite unramified) cover $Y$ of $X$ is rigid.
\end{enumerate}
\end{definition}

\begin{remark}
1) If $X$ is infinitesimally rigid, then $X$ is also locally rigid. This follows by the Kodaira-Spencer-Kuranishi theory, since $H^1(X, \Theta_X)$ is the Zariski tangent space of the germ of analytic space which is the base $\Def(X)$ of the Kuranishi semiuniversal deformation of $X$.
If  $H^1(X, \Theta_X) =0$, $\Def(X)$ is a reduced point and all deformations are locally trivial.

2) Obviously strong rigidity implies global rigidity; both global  rigidity and  \'etale  rigidity imply local rigidity.

6) The simplest example  illustrating the difference between global and infinitesimal rigidity
is the Del Pezzo surface $Z_6$ of degree $6$, blow up of the plane $\PP^2$ in three non collinear points.
It is infinitesimally rigid, but it deforms to  the weak del Pezzo surface of degree $6$,  the blow up $Z'_6$ of the plane $\PP^2$ in three  collinear points.
$Z'_6$ is not isomorphic to $Z_6$ because for the second surface the anticanonical divisor is not ample.

\end{remark}

The following  useful general result is established in \cite{bc-rigid} using many earlier results (and the Riemann Roch theorem for the second statement):

\begin{theorem}\label{m(S)}
A  compact complex manifold $X$ is rigid if and only if  the Kuranishi space $\Def(X)$ (base of the Kuranishi family of deformations)
is  $0$-dimensional. 

In particular, if $X=S$ is a smooth compact complex surface and 
$$
10 \chi(\mathcal{O}_S) - 2 K_S^2 + h^0(X, \Theta_S) >  0,
$$
then $S$ is not rigid.

\end{theorem}

We have moreover (ibidem)

\begin{prop}\label{etale}
If $p : Z \ra X$ is \'etale, i.e. a finite unramified holomorphic map between compact complex manifolds,
 then  the  infinitesimal  rigidity of $Z$ implies the  infinitesimal  rigidity of 
$X$.
Moreover,  if $Z$ is rigid, then also $X$ is rigid.

\end{prop}
{\em Idea of proof.}
For the first assertion, one observes that $  H^1(Z, \Theta_Z) =  H^1(X, p_* (\Theta_Z) ) = 0$,
and that $p_* (\Theta_Z) = p_* (p^* \Theta_X)) =  \Theta_X \otimes (p_* \hol_Z)$
has $\Theta_X $ as a direct summand.

For the second assertion one reduces to the Galois case $X = Z/G$, where,
as shown in  \cite{montecatini}: 
$$  \Def (X) = \Def(Z)^G  \subset \Def(Z).$$ 

Hence, if $ \Def (Z)$ has dimension $0$,  a fortiori  also $\Def (X)$.

\qed

While the only rigid curve is $\PP^1$, in the case of surfaces the list is only known for surfaces which are not of general type.
Indeed, using surface classification, it is shown in \cite{bc-rigid}:

\begin{theorem}
Let $S$ be a smooth compact complex  surface, which is (locally) rigid. Then either
\begin{enumerate}
\item $S$ is a minimal surface of general type, or
\item $S$ is a Del Pezzo surface of degree $d \geq 5$ (i.e., $\PP^2 $, $\PP^1 \times \PP^1$,
$S_8, S_7, S_6, S_5$, where $S_{9-r}$ is the blow-up of $\PP^2 $
in $r$ points which are in general linear position.)
\item $S$ is an Inoue surface of type $S_M$ or $S_{N,p,q,r}^{(-)}$ (cf. \cite{inoue}).
\item
Rigid surfaces in class (1) are also globally rigid, surfaces in class  (3) are  infinitesimally and globally rigid,
 surfaces in class (2) are infinitesimally rigid, but the only rigid surface in class (2)  is the projective plane $\PP^2$.
\end{enumerate}
\end{theorem}

In particular,  rigid surfaces have Kodaira dimension either $2$ (general type), or $- \infty$.
In higher dimension $ n \geq 3$  it shown in \cite{bc-rigid} that there are rigid compact complex manifolds for each possible Kodaira 
dimension, except possibly Kodaira 
dimension $= 1$. Probably this exception does not really occur, at least  for large $n$.

\begin{remark}
For surfaces of general type it is expected  to find examples which are  rigid,
but not infinitesimally rigid: such an  example would be the one of  a minimal surface $S$ such that its canonical model $X(S)$ is infinitesimally rigid
and singular (see \cite{b-w}).

\end{remark}

The intriguing part of the story is that all the old examples of globally rigid surfaces, except $\PP^2$, are projective classifying spaces.
Indeed, before the examples which we denote here by  HK(n)-surfaces,
all  known examples of rigid surfaces of general type were the following:
\begin{enumerate}
\item the  {\em ball quotients}, which  are infinitesimally rigid, strongly rigid  and \'etale rigid (\cite{siuannals}, \cite{mostow}).
\item {irreducible bi-disk quotients}, i.e. those surfaces whose  
universal covering of $S$ is $\BB_1\times \BB_1 \cong \HH \times \HH$, where $\HH$ is the upper half plane,
and moreover if we write $S = \HH \times \HH / \Gamma$ the fundamental group $\Gamma$ has dense image for any 
of the two projections $\Gamma \rightarrow PSL(2, \RR)$; they are infinitesimally rigid, strongly rigid  and \'etale rigid (\cite{J-Y85}, \cite{Mok2}).
\item {\em Beauville surfaces}; they are  infinitesimally rigid, strongly rigid  but not  \'etale rigid (\cite{isogenous}).
\item {\em Mostow-Siu surfaces}, \cite{mostow-siu}; these are infinitesimally rigid, strongly rigid  and \'etale rigid.
\item the  rigid  {\em Kodaira fibrations} constructed  by the author and Rollenske, \cite{cat-rollenske}, and mentioned in theorem \ref{Kodairarigid}; these are rigid and strongly rigid, infinitesimal rigidity
and  \'etale rigidity is not proven in \cite{cat-rollenske} but could be true.
\end{enumerate}

\begin{remark}
Examples (1)-(3), and (5) are strongly rigid.
\end{remark} 

A natural question is therefore:

\begin{question}\label{rigidnonpc}
Do there exist infinitesimally  rigid surfaces of general type which are not projective classifying spaces ?
\end{question}

 We believe that the answer should be yes, not only because wishful thinking in the case of
 surfaces almost invariably turns out to be contradicted, but for the following reason,
 which relates to the later section on Inoue -type varieties.
 
 \begin{rem}
 (i) Assume that $Z$ is a projective classifying space of dimension $n \geq 3$, and let $Y$ be
 a smooth hyperplane section of $Z$: then $Y$ is not a projective classifying space.
 
 This is a consequence of the Lefschetz hyperplane theorem: $\pi_1(Y) \cong \pi_1(Z)$.
 
 If  $Y$ were a classifying space for $\pi_1(Y) \cong \pi_1(Z)$, then 
 $$H^* (Y, \ZZ) \cong H^* (\pi_1(Y), \ZZ) \cong H^* (\pi_1(Z), \ZZ)\cong H^* (Z, \ZZ),$$
 in particular $H^{2n}(Y, \ZZ) = \ZZ$, against the fact that the real dimension of $Y$ is $ 2n-2$,
 which implies that $H^{2n}(Y, \ZZ) = 0$.
 
 (ii) The same occurs if $Y$ is an iterated hyperplane section of $Z$, with $dim (Y) \geq 2$.
 
 (iii) In particular, even if $Z$ admits a  metric of negative curvature, $Y$ cannot have one such negative metric, by the Cartan-Hadamard 
 theorem. In fact, concerning the metric inherited from $Z$, observe that only the Hermitian curvature decreases in subbundles.
 
 \end{rem}

Concerning question \ref{rigidnonpc}, the case of BCDH-surfaces, and $HK(n)$-surfaces is not completely settled. 

For the surfaces $HK(n)$ the answer is positive,  in case that $5$ divides $n$, due to  the work of Fangyang Zheng \cite{zheng}
who extended the Mostow-Siu technique  to the case of normal crossings; from this also  strong and \'etale rigidity follow in this case.

For other values of $ n \geq 4$ the work of Panov \cite{panov} gives  a positive answer for  $n > >0$ (but unspecified): his method consists in finding
polyhedral metrics of negative curvature. 

\section{Surfaces isogenous to a product and their use}

Even if the topic of surfaces isogenous to a product does appear at first sight skew to the topic of surfaces
which are ramified coverings branched on a union of lines, this idea is deceptive.

In fact, a simple way to construct a curve $C_i$ with $G$-symmetry is to construct a Galois branched covering 
of $\PP^1 \setminus \sB_i $, with group $G$, and where we assume  that the branch locus is exactly the finite set $ \sB_i $.

Hence, the surface isogenous to a product 
$  S : = (C_1 \times C_2)/G$ is a ramified covering of $\PP^1 \times \PP^1$
branched on the horizontal and vertical lines 
$$(\PP^1 \times \sB_2 )\cup (\sB_1 \times \PP^1) .$$

Note that we have $\De_G \subset G \times G$, the diagonal subgroup, and a Galois diagram
$$  C_1 \times C_2 \ra  S : = (C_1 \times C_2)/\De_G \ra \PP^1 \times \PP^1 = (C_1 \times C_2)/(G \times G).$$ 

The covering $S \ra \PP^1 \times \PP^1 $ is Galois if and only if $G$ is Abelian. 

By the main theorem on surfaces which are isogenous to a product, such a surface $S$ is (strongly) rigid if $|\sB_1 |=|\sB_2|= 3$,
and in this case in \cite{isogenous} we called these surfaces ` Beauville surfaces'.

The reason for this is that Beauville surfaces with Abelian group occur only when $G = (\ZZ/n)^2$ with $ GCD (n,6)=1$,
as shown  in \cite{isogenous}, and the original example by Beauville in \cite{beabook} was exactly the case of $G = (\ZZ/5)^2$.
For these surfaces $C_1 = C_2$ is the Fermat curve of degree $n$, and the only  difficulty consists
in finding actions such that $\De_G$ acts freely on the product.

In the article \cite{bcg} the existence problem for Beauville surfaces was translated into group theoretical terms:
because each covering $C_1 \ra \PP^1$ is determined, in view of the Riemann existence theorem, by its
branch locus (here fixed!)  and its  monodromy;
and in this case the monodromy means the datum of three elements $a,b,c \in G$ which generate $G$ and satisfy $ abc=1_G$.

Then one gets two triples $(a,b,c) (a',b',c')$ and the condition that the action on $C_1 \times C_2$ is free amounts
to the disjointness of the stabilizers, $\sS_1 \cap \sS_2= \{1\}$, here the stabilizer set   $\sS_1$ is  the union
of the conjugates of the powers of the respective elements $a,b,c$, and similarly for $\sS_2$.

In particular the question: which (nonabelian) simple groups except $\mathfrak A_5$ occur as $G$ for a Beauville surface?
It was solved for many groups in loc. cit.,  has  then 
attracted the attention of group theorists and was solved in the affirmative,
even if there are still open questions concerning whether  this can also be achieved via non real structures, via real structures,
via  strongly real structures (see \cite{gll}, \cite{gm}, \cite{fmp}, \cite{bbpv}, see also the book \cite{beasurfbook} and references
 therein for a partial account).

Now, our contention here is that Beauville surfaces not only create links between algebraic geometry and group theory,
but that they continue to yield very interesting algebro-geometric examples.

For instance, Beauville surfaces with Abelian group $G = (\ZZ /n)^2$ were used in \cite{catlaz} to answer a question posed by Jonathan Wahl,
hence giving rise to new examples of threefolds $Z$ (obtained as cones over such a  surface $S$)  which fulfill the following properties:

1)  $Z$ is Cohen-Macaulay, 

2)  the dualizing sheaf  $\omega_Z$ is torsion, 

3) the index $1$ -cover $Z'$ is not Cohen-Macaulay, in particular $Z$ is not $\QQ$-Gorenstein.

This is the  technical result answering the question by Wahl, showing in particular the existence of  regular  surfaces with subcanonical ring which
is not Cohen-Macaulay.

\begin{theo}\label{main}
For each $r = n-3$, where $n \geq 7$ is  relatively prime to $30$, and for each 
$  m,  \ 1 \leq m \leq r-1$, there are Beauville type surfaces $S$  with  $q(S) =0$ ($q(S) : =  \dim H^1 (S, \hol_S) $)
 s.t. $ K_S = r L$, and $ \  H^1 (m L) \neq  0 \ $.

\end{theo}

\subsection{Automorphisms acting trivially on cohomology}

Another application is in the direction of giving examples of surfaces admitting automorphisms which act trivially on integral cohomology
(but are not isotopic to the identity).

In this context, Cai, Liu and Zhang \cite{cai-wenfei} have proven the following theorem

\begin{theo}

Let $S$ be a minimal smooth surface of general type with $q(S) \geq 2$. Then either $S$ is rational cohomologically 
rigidified,
i.e. every automorphism acts trivially on 
$H^* (S, \QQ)$, or the subgroup $Aut(S)_{\QQ}$ of automorphisms acting trivially on the rational cohomology
algebra  is isomorphic to $ \ZZ/2$,
and $S$ satisfies 

I) $K^2_S = 8 \chi, q(S)=2$,

II)  the Albanese map is surjective,
and 

III) $S$ has a pencil of genus $b=1$.

\end{theo}

To show that the estimate is effective, they classified the surfaces isogenous to a product such that  $q(S) \geq 2$
and $Aut(S)_{\QQ} \neq 0$.

 In joint work  with Gromadzki we have shown that indeed, for one of these examples (and probably for both), the action is even trivial on $H^* (S, \ZZ)$;
 but it is   not trivial on the fundamental group. To show triviality of the action
  on $H^* (S, \ZZ)$ it suffices to show that the action is trivial on the torsion group $H_1(S, \ZZ)$;
 because  then the result follows from the universal coefficients
theorem and by Poincar\'e duality.

 The following is still an open question
 
 \begin{question}
 Do there exist surfaces of general type with nontrivial automorphisms which are isotopic to the identity?
 
 \end{question}

The  question (see \cite{handbook}) is crucial in order to compare the Kuranishi and the Teichm\"uller space of
surfaces (and higher dimensional varieties).

\section{Topological methods for moduli}

In this very short   section, devoted to  concrete moduli theory in the tradition of Kodaira and Horikawa, i.e. the  fine 
 classification  of complex projective varieties (see e.g. \cite{quintics}), we shall try to show  how in some lucky cases,
  with big fundamental group, topology helps to achieve the fine classification,
allowing explicit descriptions  of  the structure of moduli spaces. 

This was done quite effectively  in several papers (\cite{burniat1}, \cite{keumnaie},
\cite{burniat2}, \cite{burniat3}, \cite{bc-inoue}, \cite{bc-CMP}, \cite{bcf}), and in the article \cite{bms} we already amply reported on this direction of research.

For this reason the exposition in this section  shall be rather brief, we refer to \cite{bms}  for several preparatory results,
and for other related topics,  such as 
orbifold fundamental groups,
Teichm\"uller spaces, moduli spaces of curves with symmetry,  and also for an account of the results on 
 the regularity of  classifying maps, such as harmonicity, addressed by Eells and Sampson, and their complex analyticity, addressed by Siu,
 which are key ingredients for the study of moduli through topological methods (see especially \cite{5book} on this topic);
 and which lead to   rigidity and quasi-rigidity properties of projective varieties which are classifying spaces
 (meaning that their moduli spaces are completely determined by their  topology).

We shall focus here instead on  a few concrete problems in moduli theory, in particular  new constructions of surfaces with $p_g=q=0$
or $p_g=q=1$.

\subsection{Burniat surfaces and Inoue type varieties}

Among the algebraic surfaces obtained as Abelian coverings of the plane branched on interesting configurations of lines,
the oldest examples were the so-called Burniat surfaces, some surfaces with $q=p_g=0$ and $K^2_S = 6,5,4,3,2$. 

Again,  I shall skip their description, especially since I already reported on them at the Kinosaki Conference in the Fall of 2011; even if there are still 
interesting open questions concerning the connected component of the moduli space containing Burniat surfaces with
 $K_S^2= 3$.

Following a suggestion of Miles Reid, Masahisa Inoue \cite{inoue} gave another description of the Burniat surfaces, as quotients of a hypersurface
in the product of three elliptic curves. Using this method, he went further with his construction, and obtained new (minimal) surfaces
of general type with $q=p_g=0$ and $K^2_S = 7$, which are now  called {\bf (algebraic) Inoue surfaces}.

These were given as quotients of a complete intersection of two surfaces inside the product of four elliptic curves,
but a closer inspection showed that indeed they are quotients of a surface inside the product
of a curve of genus $5$ with two elliptic curves. More precisely, an Inoue surface $S$ admits an unramified
$(\ZZ / 2\ZZ)^5$ - Galois covering $\hat{S}$ which is an
ample divisor in $E_1 \times E_2 \times D$, where
$E_1, E_2$ are elliptic curves and $D$ is a projective  curve of genus
$5$.

Hence  the fundamental group of an Inoue
surface with $p_g = 0$ and $K_S^2 =7$   sits
in an extension ($\Pi_g$ being as usual  the fundamental group of a projective
curve of genus $g$):
$$ 1 \rightarrow \Pi_{5} \times \mathbb{Z}^4 \rightarrow \pi_1(S)
\rightarrow (\mathbb{Z}/2\mathbb{Z})^5 \rightarrow 1.
$$

It turned out that the ideas needed to treat the moduli space of this special family of Inoue surfaces could be put in a rather general framework,
valid in all dimensions,  and together with I. Bauer 
 we proposed the study,  obtaining  several results, of what we called 
Inoue-type varieties. 

\begin{defin} {\bf (\cite{bc-inoue})}
Define a complex projective manifold $X$ to be an {\bf Inoue-type manifold} if
\begin{enumerate}
\item
$ dim (X) \geq 2$;
\item
there is a finite group $G$ and an   unramified $G$-covering $
\hat{X} \ra X$,
(hence $ X = \hat{X} / G$) such  that
\item
$ \hat{X}$ is an ample divisor inside a $K(\Ga, 1)$-projective manifold
$Z$, (hence by the theorem of  Lefschetz  $\pi_1 ( \hat{X}) \cong \pi_1 (Z)
\cong \Ga$) and moreover
\item
the action of $G$ on $ \hat{X}$   yields a faithful action on $\pi_1
( \hat{X}) \cong \Ga$:
in other words the exact sequence
$$ 1 \ra \Ga  \cong \pi_1 ( \hat{X}) \ra \pi_1 ( X) \ra G \ra 1$$
gives an injection $ G \ra \Out (\Ga)$, defined by conjugation by lifts of elements of $G$;
\item
the action of $G$ on $ \hat{X}$   is induced by an action of $G$ on $Z$.
\end{enumerate}
 We  say  that an  Inoue-type manifold $X$  is 
  a 
  {\bf special Inoue type manifold}
if  moreover
$$ Z = (A_1 \times \dots \times A_r) \times  (C_1 \times \dots \times
C_h) \times (M_1 \times \dots \times
M_s)$$ where each
$A_i$ is an Abelian variety,  each $C_j$ is a curve of genus $ g_j \geq 2$,
and $M_i$ is  a
compact quotient of an irreducible bounded symmetric domain of 
dimension at least 2 by a
torsion free subgroup;
and   a {\bf 
classical Inoue type manifold}
if  instead
$ Z = (A_1 \times \dots \times A_r) \times  (C_1 \times \dots \times
C_h) $ where as above each
$A_i$ is an Abelian variety,  each $C_j$ is a curve of genus $ g_j \geq 2$.
\end{defin}

The main idea underlying  the importance of this notion is the fact that, under suitable assumptions,
one can say that if $X$ is an Inoue type manifold, and $Y$ is homotopically equivalent to $X$
(indeed even some weaker cohomological conditions are sufficient),
then also $Y$ is an  Inoue type manifold, and in special cases one can conclude that $Y$ and $X$ belong to the same irreducible connected 
component of the moduli space. We omit to state the general results, referring to \cite{bc} and \cite{bms}, and here we shall just treat a concrete case,
in the next subsection.

\subsection{Bagnera-de Franchis varieties and applications to moduli}

In our article \cite{bcf2}, appeared  in the volume dedicated to Kodaira, we treated a special case of the theory
of Inoue type varieties, the one where the action of $G$ happens to be free also on $Z$, and $Z$ is the simplest projective classifying space, an Abelian variety.

Define the 
{\bf Generalized Hyperelliptic Varieties}  (GHV) as the quotients $A/G$ of an Abelian Variety $A$ by a finite group $G$ acting freely,
and with the property  that $G$ is not a subgroup of the group of  translations. Without loss of generality one can then assume
that $G$ contains no translations, since the subgroup $G_T$ of translations in $G$ would be a normal subgroup, and
if we denote $G' = G/G_T$, then $A/G = A' / G'$, where $A'$ is the Abelian variety $ A' : = A/G_T$.

A smaller class is the class of  {\bf Bagnera-de Franchis (BdF) Varieties}: these are the  quotients $X = A/G$ were $G$ contains no translations,
and $G$ is a cyclic group of order $m$, with generator $g$ (observe that, when $A$ has  dimension $n=2$, the two notions coincide, thanks
to the classification result of Bagnera-de Franchis in \cite{BdF}).

Bagnera-de Franchis varieties have a simple description as quotients of Bagnera-de Franchis varieties of product type,
according to the following definition:

\begin{defin}
A  Bagnera-de Franchis manifold (resp.: variety)  {\bf  of product type } is a quotient $ X= A/G $ where $A = A_1 \times A_2$,  $A_1, A_2$ are complex tori
(resp.: Abelian Varieties), 
and $G \cong \ZZ/m$ is a cyclic group operating freely on $A$, generated by an automorphism of the form
$$ g (a_1, a_2 ) = ( a_1 + \be_1, \al_2 (a_2)) ,$$
where $\be_1 \in A_1[m]$ is an element of order exactly $m$, and similarly $\al_2 : A_2 \ra A_2$ is a linear automorphism
of order exactly $m$ without $1$ as eigenvalue (these conditions guarantee that the action is free). 

\end{defin}

This is then the  characterization of general Bagnera- de Franchis varieties.

\begin{prop}\label{BdF}
Every Bagnera-de Franchis variety $ X= A/G $, where  $G \cong \ZZ/m$ contains no translations, is the quotient of a  Bagnera-de Franchis variety
of product type,  $(A_1 \times A_2)/ G$ by any finite subgroup $T$ of $A_1 \times A_2$ which satisfies the following properties:

1) $T$  is the graph of an isomorphism between two respective 
subgroups $ T_1 \subset A_1,  T_2 \subset A_2,$

2) $(\al_2 - \Id) T_2 = 0$

3) if $ g (a_1, a_2 ) = ( a_1 + \be_1, \al_2 (a_2)) ,$ then the subgroup of order $m$ generated by $\be_1$ intersects $T_1$ only in $\{0\}$.

In particular, we may write $X$ as the quotient $ X =  (A_1 \times A_2)/ (G \times T)$ by the abelian group $G \times T$. 

\end{prop}

This notion was then  used in \cite{bcf} to make a construction that we briefly describe.

Let $A_1$ be an elliptic curve, and let $A_2$ be an Abelian surface with a line bundle $L_2$ yielding a polarization of type $(1,2)$.
Take as $L_1$ the line bundle $\hol_{A_1} ( 2 O)$,
and let $L$ be the line bundle on $A' : = A_1 \times A_2$ obtained as the exterior tensor product of $L_1$ and $L_2$,
so that
$$ H^0 (A', L) =  H^0 (A_1, L_1)  \otimes  H^0 (A_2, L_2) .$$

Moreover, choose the origin in $A_2$ so that the space of sections $H^0 (A_2, L_2) $ consists only of {\em even} sections.

We  take then, using properties of the Stone-von Neumann representation of the Heisenberg group,
a Bagnera-de Franchis threefold $ X: = A/ G$, where $A = (A_1 \times A_2) / T$, and $ G \cong T \cong \ZZ/2$,
and  a surface $S \subset X$ which is the quotient of a $(G\times T)$-invariant $D \in |L|$, so that $S^2 = \frac{1}{4} D^2 = 6$.

We could then prove the following.

\begin{theo}\label{BCF}
Let $S$ be a surface of general type with invariants $K_S^2 = 6$, $p_g =q = 1$ such that there exists an unramified double cover
$ \hat{S} \ra S$ with $ q (  \hat{S} ) = 3$, and such that the Albanese morphism $ \hat{\al} :  \hat{S}  \ra A$ is birational onto its image $Z$,
a divisor in $A$ with  $ Z^3 = 12$.

Then the canonical model of $\hat{S}$ is isomorphic to $Z$, and the canonical model of $S$ is isomorphic to $Y = Z / (\ZZ/2)$, 
which a divisor in a Bagnera-De Franchis threefold $ X: = A/ G$, where $A = (A_1 \times A_2) / T$, $ G \cong T \cong \ZZ/2$,
and where the action is given by

\begin{equation}  G : = \{ Id, g \}, \ \ g (a_1 + a_2 ) : = a_1 + \tau /2 - a_2 + \la_2 / 2, \\\forall a_1 \in A_1, a_2 \in A_2\\
  \end{equation} 
  \begin{equation} 
  T : = ( \ZZ/2 ) ( 1/2 + \la_4 / 2) \subset A = (A_1 \times A_2) .
  \end{equation} 

These surfaces exist, have an irreducible four dimensional moduli space, and their Albanese map $\al : S \ra A_1 = A_1/ A_1[2]$ has 
general fibre a non hyperelliptic curve of genus $g=3$.

\end{theo}

\bigskip

\bigskip

{\bf Acknowledgement: I would like to thank the organizers for their kind invitation and for giving me the opportunity to meditate
on several interesting questions, Valery Alexeev for reminding me of Easton's construction, 
and some members of our algebraic geometry seminar in Bayreuth  for pointing out misprints and inaccuracies
in a preliminary version.}

\end{document}